\documentclass[10pt]{amsart}   	% use "amsart" instead of "article" for AMSLaTeX format
\usepackage[hmargin=1in,vmargin=1in]{geometry}                  		% See geometry.pdf to learn the layout options. There are lots.
\usepackage{graphicx}				% Use pdf, png, jpg, or eps§ with pdflatex; use eps in DVI mode
\usepackage{amssymb,amsthm,amsmath,mathrsfs,verbatim}
\usepackage[utf8]{inputenc}
\usepackage[english]{babel}
\usepackage{fancyhdr,color}
\usepackage{graphicx,enumerate, fancyhdr, microtype}
\usepackage{xcolor,bm,color,tikz,euscript,sansmath,ytableau}
\usepackage{tikz-cd}

\usepackage[all]{xy}

\usetikzlibrary{3d}
\usetikzlibrary{positioning,matrix,arrows,decorations.pathmorphing}

\usetikzlibrary{positioning,arrows}
\usetikzlibrary{decorations.markings}
\tikzstyle{circledvertex} = [draw, black, shape=circle, minimum size=8pt, inner sep=1pt]
\tikzstyle{invisivertex} = [black, shape=rectangle, minimum size=0pt, inner sep=2pt]
\tikzstyle{point}=[draw, black, fill,shape=circle, minimum size=4pt, inner sep=0pt]
\tikzstyle{magpoint}=[draw, magenta, fill,shape=circle, minimum size=4pt, inner sep=0pt]

\tikzstyle{origVertex} = [draw, shape=circle, minimum size=0pt, inner sep=0pt]
\tikzstyle{boundaryEdge} = [magenta, densely dashed]

\newcommand{\R}{\mathbb{R}}
\newcommand{\C}{\mathbb{C}}
\newcommand{\Z}{\mathbb{Z}}

\newcommand{\A}{\mathcal{A}}
\newcommand{\K}{\mathcal{K}}
\newcommand{\cham}[1]{\mathcal{C}({#1})} % Use \cham{<<blah>>} to denote the chambers of a cone or arrangement <<blah>>.
\newcommand{\inter}[1]{\mathcal{L}^{\operatorname{int}}(#1)} % Use \inter{<<blah>>} to denote the intersection poset of a cone <<blah>>.

\newcommand{\hypinter}[1]{\mathcal{L}(#1)} % Use \hypinter{<<blah>>} to denote the intersection poset of an arrangement <<blah>>.
\newcommand{\codim}[1]{\operatorname{codim}(#1)} % Use \codim{<<blah>>} to denote the codimension of a subspace.

\newcommand{\whit}[2]{\thinspace c_{#1}(#2)\thinspace} % Use \whit{k}{<<blah>>} to denote the kth Whitney number of of a cone or arrangement <<blah>>.

\newcommand{\linext}[1]{\operatorname{LinExt}(#1)}

\newcommand{\antichain}{{\sf Antichain}}
\newcommand{\Poin}{{\sf Poin}}
\newcommand{\pairs}{{\sf pairs}}

\newcommand{\PTP} %{\mathfrak{S}_{\operatorname{tran}}}
 {\mathfrak{S}^{\pitchfork}}
 \newcommand{\PTPar}%{\Pi_{\operatorname{tran}}}
  {\Pi^{\pitchfork}}
\newcommand{\LRmax}{{\sf LRmax}}
\newcommand{\MLRmax}{{\sf MLRmax}}
\newcommand{\Sym}{\mathfrak{S}}
\newcommand{\Symp}{\mathfrak{S}^{\operatorname{(pcyc)}}}
\newcommand{\Ma}{{M(\overline{a})}}
\newcommand{\Pa}{{P_{\overline{a}}}}
\newcommand{\Phia}{{\Phi^*}}
\newcommand{\Psia}{{\Psi^*}}
\newcommand{\edit}[1]{\textcolor{black}{#1}}

\newcommand{\cint}[1]{\mathsf{int}(#1)} % Use \cint{<<blah>>} to denote the interior of a cone.
\newcommand{\relint}[1]{\mathsf{relint}(#1)} % Use \relint{<<blah>>} to denote the relative interior of a cone.

\makeatletter
\renewcommand*\env@matrix[1][*\c@MaxMatrixCols c]{%
  \hskip -\arraycolsep
  \let\@ifnextchar\new@ifnextchar
  \array{#1}}
\makeatother

\newcommand{\des}{\mathsf{des}}
\newcommand{\opp}{\mathsf{opp}}
\newcommand{\Des}{\mathsf{Des}}
\newcommand{\chain}[1]{\mathsf{#1}}
\newcommand{\stat}{\mathsf{stat}}
\newcommand{\blocks}{\mathsf{blocks}}

\newcommand{\supp}[1]{\mathsf{supp}(#1)} % support of a multiset permutation #1
\newcommand{\foata}[1]{\mathsf{pcyc}(#1)} %Number of prime Cycles
\newcommand{\fc}{\foata{\sigma}} %Number of Foata Cycles
\newcommand{\cyc}{\mathsf{cyc}} %Number of cycles
\newcommand{\multi}[1]{M(\overline{#1})} % \multi{A} is the multiset defined by a vector A
\newcommand{\poset}[1]{P_{\overline{#1}}} % \poset{A} is the disjoint union of chains of lengths given by the vector A
\newcommand{\total}[1]{|\overline{#1}|} % \total{A} is the sum of the entries of a vector A

\newcommand{\SN}{\mathfrak{S}_n}
\newcommand{\SH}{\mathfrak{S}_H}

\newcommand{\Q}{\operatorname{Int}_{\ell}}
\newcommand{\bb}[1]{\underline{b}_{#1}}

\theoremstyle{plain}
\newtheorem{thm}{Theorem}[section]
\newtheorem{lem}[thm]{Lemma}
\newtheorem{coro}[thm]{Corollary}

\newtheorem{prop}[thm]{Proposition}

\newtheorem{prob}[thm]{Open Problem}

\theoremstyle{definition}
\newtheorem{defn}[thm]{Definition}
\newtheorem{ex}[thm]{Example}

\newtheorem{question}[thm]{Question}

\theoremstyle{remark}
\newtheorem{rem}[thm]{Remark}

\title{Whitney Numbers for Poset Cones}
\author{Galen Dorpalen-Barry, Jang Soo Kim and Victor Reiner}
\date{\today}
\address{School of Mathematics, University of Minnesota, Minneapolis, MN 55455, USA}
\email{dorpa003@umn.edu}

\address{Department of Mathematics, Sungkyunkwan University, Suwon, Gyeonggi-do 16419, South Korea}
\email{jangsookim@skku.edu}

\address{School of Mathematics, University of Minnesota, Minneapolis, MN 55455, USA}
\email{reiner@math.umn.edu}

\thanks{The first and third authors were partially supported by NSF grant DMS-1601961.
The second author was supported by NRF grants \#2019R1F1A1059081 and \#2016R1A5A1008055.}


\begin{document}

%%%%%%%%%%%%%%%%%%%%%%%%%%%%%%%%%%%%%%%%%%%%%%%%%%%%%%%%%%%%%%%%%%%%%%%%%%%%%%

\begin{abstract}
Hyperplane arrangements dissect $\R^n$ into connected components called chambers, and a well-known theorem of Zaslavsky counts chambers as a sum of nonnegative integers called
Whitney numbers of the first kind.  His theorem generalizes to count chambers within any cone defined as the intersection of a collection of halfspaces from the arrangement, leading to a notion of Whitney numbers for each cone.
This paper focuses on cones within the braid arrangement, consisting of the reflecting hyperplanes $x_i=x_j$ inside $\R^n$ for the symmetric group, thought of as the type $A_{n-1}$ reflection group.  Here,
\begin{itemize} \item cones correspond to posets,
\item chambers within the cone correspond to linear extensions of the poset,
\item the Whitney numbers of the cone interestingly refine the number of linear extensions of the poset.
\end{itemize}
We interpret this refinement for all posets as counting linear extensions according to a
statistic that generalizes the number of left-to-right maxima of a permutation.
When the poset is a disjoint union of chains, we interpret this refinement differently,
 using Foata's theory of cycle decomposition for multiset permutations, leading to a simple generating
function compiling these Whitney numbers.
\end{abstract}

\maketitle

%\tableofcontents

%%%%%%%%%%%%%%%%%%%%%%%%%%%%%%%%%%%%%%%%%%%%%%%%%%%%%%%%%%%%%%%%%%%%%%%%%%%%%%%

\section{Introduction}

%%%%%%%%%%%%%%%%%%%%%%%%%%%%%%%%%%%%%%%%%%%%%%%%%%%%%%%%%%%%%%%%%%%%%%%%%%%%%%%

This paper concerns \emph{arrangements} $\A=\{H_1,\ldots,H_m\}$ of
hyperplanes $H_i$, which are affine-linear codimension one subspaces of $V=\R^n$.  Each such arrangement dissects $V$ into the connected components of its complement $V \setminus \bigcup_{i=1}^m H_i$, called \emph{chambers}.  Let $\cham{\A}$ denote the collection of all such chambers.

The theory of hyperplane arrangements is rich and well-explored, with connections to reflection groups, braid groups, random walks, card-shuffling, and discrete geometry of polytopes and oriented matroids; see \cite{OrlikTerao, stanley}.  In particular, the
number of chambers $\#\cham{\A}$ has a famous formula
due to Zaslavsky \cite{zaslavsky-memoir}, expressed in terms of the \emph{intersection poset} $\hypinter{\A}$, consisting of all {\it intersection subspaces} $X=H_{i_1} \cap H_{i_2} \cap \cdots \cap H_{i_k}$, ordered via reverse inclusion.
This poset is known to have the property that every \emph{lower interval}
\[
[V,X]:=\{Y \in \hypinter{\A}: V \leq Y \leq X\}
\]
from its unique bottom element $V$ to any intersection space $X$ forms a \emph{geometric lattice}.
In particular, each such $[V,X]$ is a \emph{ranked} poset,
where the rank of $X$ is $\codim{X}:=n-\dim(X)$.  Zaslavsky's result says that
\begin{equation}
    \label{Zaslavsky-famous-theorem}
\#\cham{\A} = \sum_{ X \in \hypinter{\A}} |\mu(V,X)|
=\sum_{k=0}^n \whit{k}{\A}
=\left[ \Poin(\A,t)\right]_{t=1},
\end{equation}
where $\mu(-,-)$ denotes the \emph{M\"obius function} of $\hypinter{\A}$,
while the nonnegative integers
$$
\whit{k}{\A}:=\sum_{ \substack{X \in \hypinter{\A}:\\ \codim{X}=k} } |\mu(V,X)|,
$$
are often called the \emph{(signless) Whitney numbers of the first kind} for $\A$,
and their generating function
$$
\Poin(\A,t):=\sum_{k\ge0} \whit{k}{\A} t^k
$$
is called the \emph{Poincar\'e polynomial}\footnote{This is because it is the generating function
for the Betti numbers of the complexified complement $\C^n \setminus \A$; see \cite[Chap. 5]{OrlikTerao}}.

Our starting point is a less widely-known generalization of equation (\ref{Zaslavsky-famous-theorem}),
also due to Zaslavsky \cite{zaslavsky}.   It applies more generally to count the chambers of $\A$ that lie
within a \emph{cone} $\K$, defined to be the intersection of any collection of open halfspaces for hyperplanes of $\A$; said differently, a cone $\K$ of $\A$ is a chamber in $\cham{\A'}$ for some subarrangement $\A' \subseteq \A$.  Results on the set $\cham{\K}$ of all chambers of $\A$ inside a cone $\K$ have appeared more recently in work of Brown on random walks \cite{brown}, and in work of Gente on Varchenko determinants \cite[Section 2.4]{gente}, and independently Aguiar and Mahajan \cite[Theorem 8.22]{aguiar}.
Define the \emph{poset of interior intersections for $\K$} to be
the following order ideal in $\hypinter{\A}$:
\[
\inter{\K} = \{X\in\hypinter{\A} \mid X \cap
\K\not=\emptyset\}.
\]
Zaslavsky observed in \cite[Example A, p.  275]{zaslavsky} that \eqref{Zaslavsky-famous-theorem} generalizes to cones $\K$, asserting
\begin{equation}
    \label{Zaslavsky-less-famous-theorem}
\#\cham{\K} = \sum_{ X \in \inter{\K}} |\mu(V,X)|
=\sum_{k=0}^n \whit{k}{\K}
= \left[\Poin(\K,t) \right]_{t=1}.
\end{equation}
Here we again define nonnegative integers, the \emph{(signless) Whitney numbers of the first kind} for the cone $\K$
$$
\whit{k}{\K}:=\sum_{ \substack{X \in \inter{\K}:\\ \codim{X}=k} } |\mu(V,X)|,
$$
with generating function
$\Poin(\K,t):=\sum_k \whit{k}{\K} t^k$, which we call the \emph{Poincar\'e polynomial} of $\K$.

For example, inside $\A = \{H_1, H_2,H_3, H_4, H_5\}$ in $\R^2$ shown
  below on the left, we have shaded one of four possible cones $\K$
  defined by the subarrangement $\A'=\{H_4,H_5\}$, containing $\#\cham{\K}=5$ chambers of $\A$:

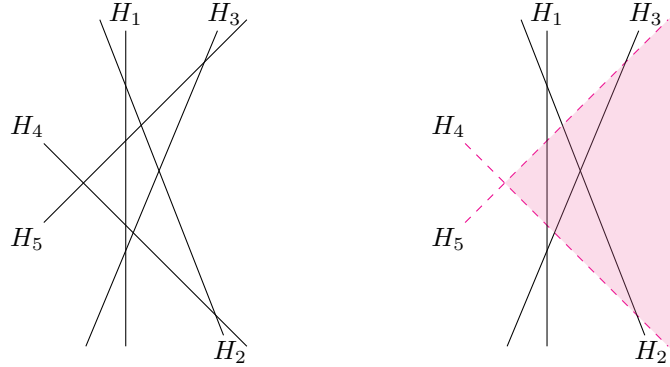
    \begin{center}
      \begin {tikzpicture}[scale=.75]

      \node[invisivertex] (A1) at (-1,1){$H_4$};
      \node[invisivertex] (A2) at (3,-3){};
      \node[invisivertex] (B1) at (.75,-3){};
      \node[invisivertex] (B2) at (.75,3){$H_1$};
      \node[invisivertex] (C1) at (2.6,-3){$H_2$};
      \node[invisivertex] (C2) at (.25,3){};
      \node[invisivertex] (D1) at (0,-3){};
      \node[invisivertex] (D2) at (2.5,3){$H_3$};
      \node[invisivertex] (E1) at (-1,-1){$H_5$};
      \node[invisivertex] (E2) at (3,3){};

      % \node[invisivertex] (arrangement) at (1,4){Arrangement};

      % \node[invisivertex] (chamber1) at (.5,0){\textcolor{magenta}{$1$}};
      % \node[invisivertex] (chamber2) at (1.15,0){\textcolor{magenta}{$2$}};
      % \node[invisivertex] (chamber3) at (1.5,1){\textcolor{magenta}{$3$}};
      % \node[invisivertex] (chamber4) at (1.5,-1){\textcolor{magenta}{$4$}};
      % \node[invisivertex] (chamber5) at (2.5,0){\textcolor{magenta}{$5$}};

      \path[-] (A1) edge [bend left =0] node[above] {} (A2);
      \path[-] (E1) edge [bend left =0] node[above] {} (E2);
      \path[-] (B1) edge [bend left =0] node[above] {} (B2);
      \path[-] (C1) edge [bend left =0] node[above] {} (C2);
      \path[-] (D1) edge [bend left =0] node[above] {} (D2);

      %\fill[magenta!40,nearly transparent] (0,0) -- (3,3) -- (3,-3) -- cycle;
      \end{tikzpicture}
\hspace{2cm}
  \begin {tikzpicture}[scale=.75]

  \node[invisivertex] (A1) at (-1,1){$H_4$};
  \node[invisivertex] (A2) at (3,-3){};
  \node[invisivertex] (B1) at (.75,-3){};
  \node[invisivertex] (B2) at (.75,3){$H_1$};
  \node[invisivertex] (C1) at (2.6,-3){$H_2$};
  \node[invisivertex] (C2) at (.25,3){};
  \node[invisivertex] (D1) at (0,-3){};
  \node[invisivertex] (D2) at (2.5,3){$H_3$};
  \node[invisivertex] (E1) at (-1,-1){$H_5$};
  \node[invisivertex] (E2) at (3,3){};

  % \node[invisivertex] (chamber1) at (.5,0){\textcolor{magenta}{$1$}};
  % \node[invisivertex] (chamber2) at (1.15,0){\textcolor{magenta}{$2$}};
  % \node[invisivertex] (chamber3) at (1.5,1){\textcolor{magenta}{$3$}};
  % \node[invisivertex] (chamber4) at (1.5,-1){\textcolor{magenta}{$4$}};
  % \node[invisivertex] (chamber5) at (2.5,0){\textcolor{magenta}{$5$}};

  \path[-] (A1) edge [bend left =0, color=magenta,dashed] node[above] {} (A2);
  \path[-] (E1) edge [bend left =0, color=magenta,dashed] node[above] {} (E2);
  \path[-] (B1) edge [bend left =0] node[above] {} (B2);
  \path[-] (C1) edge [bend left =0] node[above] {} (C2);
  \path[-] (D1) edge [bend left =0] node[above] {} (D2);

  \fill[magenta!70,nearly transparent] (0,0) -- (3,3) -- (3,-3) -- cycle;
  \end{tikzpicture}
\end{center}

\noindent
Zaslavsky's formula \eqref{Zaslavsky-less-famous-theorem} computes this as follows.  The poset of interior intersections $\inter{\K}$ has Hasse diagram:

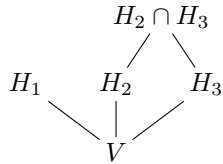
\begin{center}
\begin {tikzpicture}[scale=.6]

\node[invisivertex] (0) at (3,0){$V$};
\node[invisivertex] (H1) at (1,1.5){$H_{1}$};
\node[invisivertex] (H2) at (3,1.5){$H_{2}$};
\node[invisivertex] (H3) at (5,1.5){$H_{3}$};
\node[invisivertex] (X23) at (4,3){$H_{2}\cap H_3$};

\path[-] (0) edge [bend left =0] node[above] {} (H2);
\path[-] (0) edge [bend left =0] node[above] {} (H3);
\path[-] (0) edge [bend left =0] node[above] {} (H1);

\path[-] (H2) edge [bend left =0] node[above] {} (X23);
\path[-] (H3) edge [bend left =0] node[above] {} (X23);
\end{tikzpicture}
\end{center}

\noindent
Here $\mu(V,X)=(-1)^{\codim{X}}$ for all $X$, so that $(c_0(\K),c_1(\K),c_2(\K))=(1,3,1)$, and
$$
\begin{array}{ccccccc}
            & &\Poin(\K,t)&=&c_0(\K)+c_1(\K)t+c_2(\K)t^2&=&1+3t+t^2,\\
\#\cham{\K}&=&\left[\Poin(\K,t)\right]_{t=1}&=&c_0(\K)+c_1(\K)+c_2(\K)&=&5.
\end{array}
$$

\medskip
The object of this paper is to understand the distribution of the signless Whitney numbers as a refinement of $\#\cham{\K}$ as in equation \eqref{Zaslavsky-less-famous-theorem}, for cones $\K$ in the \emph{braid arrangement}. The braid arrangement
$A_{n-1}=\{ H_{ij}\}_{1 \leq i < j \leq n}$ is the set of $\binom{n}{2}$ reflecting hyperplanes
  \[
  H_{ij}=\{(x_1,\dots,x_n)\in V=\R^n\mid x_i-x_j=0\}
  \]
for the symmetric group $\SN$ on $n$ letters,
thought of as the reflection group of type $A_{n-1}$.
There is a well-known and easy bijection between the chambers $\cham{A_{n-1}}$ and the permutations $\sigma=[\sigma_1, \sigma_2, \ldots, \sigma_n]$ in $\SN$,
sending $\sigma$ to the chamber:
\begin{equation}
    \label{braid-arrangement-chambers}
\K_\sigma:=\{ x \in V=\R^n: x_{\sigma_1} < x_{\sigma_2} < \cdots < x_{\sigma_n}\}.
\end{equation}
More generally, one has an easy bijection, reviewed in Section \ref{sect:type A}, between posets on the set $[n]:=\{1,2,\ldots,n\}$ and cones  in the braid arrangement $A_{n-1}$,
sending a poset $P$ to the cone
$$
\K_P:=\{ x \in V=\R^n: x_i < x_j \text{ for }i <_P j\}.
$$
It is readily checked that the chamber $\K_\sigma$ lies in the cone $\cham{\K_P}$ if and only \edit{if} $\sigma$ is a \emph{linear extension} of $P$, meaning that
$i <_P j$ implies $i <_\sigma j$, regarding $\sigma$ as a total order $\sigma_1 < \sigma_2 < \cdots < \sigma_n$ on
$[n]$.  Letting $\linext{P}$ denote the set of all linear extensions of  $P$, this shows that
$\#\cham{\K_P}=\#\linext{P}$, and hence equation \eqref{Zaslavsky-less-famous-theorem}
becomes
\begin{equation}
\label{linear-extensions-as-sums-of-Whitneys}
\#\linext{P} = \sum_{k \geq 0} \whit{k}{P} = \left[ \Poin(P,t) \right]_{t=1}
\end{equation}
abbreviating $\whit{k}{P}:=\whit{k}{\K_P}$ and $\Poin(P,t):=\Poin(\K_P,t)$ from now on.

A motivating special case occurs when $P=\antichain_n$ is the \emph{antichain poset} on $[n]$ that has no order relations. In this
case $\linext{\antichain_n}=\SN$ itself, and the signless Whitney number $\whit{n-k}{\antichain_n}$ of the first kind is well-known
\cite[Prop. 1.3.7]{ec1} to be the \emph{signless Stirling number of the first kind} $c(n,k)$ that counts the permutations in
$\mathfrak{S}_n$ having $k$ cycles.  Consequently, \eqref{linear-extensions-as-sums-of-Whitneys} becomes the easy
summation formula
\begin{equation}
\label{easy-n!-summation}
n! =\#\mathfrak{S}_n= \sum_{k\ge0} c(n,k).
\end{equation}
Our {\bf main results} will generalize
three well-known expressions for $\Poin(\antichain_n,t)$, explained below:
\begin{align}
  \label{antichain-cycles-interpretation}
  \Poin(\antichain_n,t) = \sum_{k\ge0} c(n,k) t^{n-k}
                         &= \sum_{\sigma \in \SN} t^{n-\cyc(\sigma)}\\
  \label{antichain-LRmax-interpretation}
                        &= \sum_{\sigma \in \SN} t^{n-\LRmax(\sigma)}\\
  \label{stirling-number-gf}
                        &= 1(1+t)(1+2t) \cdots (1+(n-1)\thinspace t).
\end{align}

Equation \eqref{antichain-cycles-interpretation} comes from the definition of $c(n,k)$ mentioned above, where
$\cyc(\sigma)$ is the number of cycles of $\sigma$.
This interpretation of the Poincar\'e polynomial
will be generalized to all posets in Theorem~\ref{tranverse-perm-cycles-thm}, and generalized in
a different way to posets which are {\it disjoint unions of chains} in Theorem~\ref{chains-thm}.

Equation \eqref{antichain-LRmax-interpretation} arises from
a well-known bijection $\tau \mapsto \sigma$ from $\Sym_n$ to itself, called the \emph{fundamental bijection}\footnote{\edit{This bijection, often called {\it Foata's first fundamental transformation}, seems to have appeared also in work of R\'enyi \cite[\S 4]{Renyi}}.}, such that
$\cyc(\tau)=\LRmax(\sigma)$, see \cite[Proposition 1.3.1]{ec1}; here $\LRmax(\sigma)$ is the number of \emph{left-to-right maxima} of $\sigma=\sigma_1\sigma_2\dots\sigma_n$, i.e., the number of integers $1\le j\le n$ such that $\sigma_j =
\max\{\sigma_1,\sigma_2,\dots,\sigma_j\}$.  This bijection and interpretation
of the Poincar\'e polynomial
will be generalized to all posets in Theorem~\ref{linear-extensions-LRmax-thm}.

Equation~\eqref{stirling-number-gf} is well-known \cite[Proposition 1.3.7]{ec1}, and will be generalized
to a generating function compiling the Poincar\'e polynomials of all
posets which are disjoint unions of chains in Theorem~\ref{chains-gf-thm} below.

The remainder of this paper is organized as follows.
\vskip.1in

Section~\ref{sect:type A} gives preliminaries on the intersection lattice and cones
in braid arrangements.  In particular, for a poset $P$ on $[n]$,
it gives an explicit combinatorial description of the interior intersections
$\inter{\K_P}$ for a poset cone $\K_P$, as the subset $\PTPar(P)$
of \emph{$P$-transverse partitions} inside the lattice $\Pi_n$ of set partitions of
$[n]$;  roughly speaking, these are the set partitions $\pi$ for which the
quotient preposet $P/\pi$ does not collapse any strict order relations $i <_P j$ into equalities.
With this in hand, it defines the subset $\PTP(P)$ of \emph{$P$-transverse permutations}
inside $\SN$ as those permutations $\sigma$ whose cycles form a $P$-transverse partition, in order to prove
the following generalization of equation \eqref{antichain-cycles-interpretation}.

\begin{thm}
\label{tranverse-perm-cycles-thm}
For any poset $P$ on $[n]$, one has
\edit{\[
  \Poin(P,t) = \sum_{\sigma\in\PTP(P)} t^{n - \cyc(\sigma)}.
\]}
In particular, $\#\linext P=\#\PTP(P)$.
\end{thm}
\noindent
This result even generalizes to a statement (Theorem~\ref{W-transverse-elements-thm})
about cones in any {\it real reflection arrangement}.

In light of Theorem~\ref{tranverse-perm-cycles-thm}, one expects bijections
between $\PTP(P)$ and $\linext P$.
One such bijection is the goal of Section~\ref{sect:bijection}, which defines a notion of \emph{$P$-left-to-right maximum} for a linear
extension of $P$, and
then generalizes the fundamental bijection $\tau \mapsto \sigma$
and equation \eqref{antichain-LRmax-interpretation} as follows.

\begin{thm}
  \label{linear-extensions-LRmax-thm}
  For any poset $P$ on $[n]$, one has a bijection
  $$
  \begin{array}{rcl}
    \PTP(P) &\overset{\Phi}{\longrightarrow}& \linext P \\
    \tau &\longmapsto& \sigma
  \end{array}
  $$
  such that $\cyc(\tau)=\LRmax_{P}(\sigma)$,
  where  $\LRmax_P(\sigma)$ is the number of $P$-left-to-right maxima of $\sigma$.
Therefore,
\[
\Poin(P,t) =  \sum_{\sigma \in \linext{P}} t^{n-\LRmax_{P}(\sigma)}.
\]
\end{thm}

Section \ref{multisets} examines posets $P$ which are disjoint unions of chains,
and produces a second interesting bijection $\PTP(P)\to \linext P$.
Given any \emph{composition} $\overline{a}=(a_1,\ldots,a_\ell)$ of $n$, meaning that
$\overline{a} \in \{0,1,2,\ldots\}^\ell$ and
$\total{a}:=\sum_{i=1}^\ell a_i=n$, let
$\chain{a_i}$ denote a chain (totally ordered) poset on $a_i$ elements, and then
$$
\poset{a}:=\chain{a_1} \sqcup \chain{a_2} \sqcup \cdots \sqcup \chain{a_\ell}
$$
is a disjoint union of incomparable chains of sizes $a_1,a_2,\ldots,a_\ell$.
Here one can identify elements in $\linext{\poset{a}}$ with \emph{multiset permutations} of $\{1^{a_1}, 2^{a_2}, \dots,
\ell^{a_\ell}\}$.  Section~\ref{multisets} reviews the beautiful theory of \emph{prime cycle decompositions}
for such multiset permutations due to
Foata \cite{foatathesis}, and then reinterprets it as giving a bijection
$\linext{\poset{a}} \rightarrow \PTP(\poset{a})$, and the following (second)
generalization of equation \eqref{antichain-cycles-interpretation}.

\begin{thm}
  \label{chains-thm}
  For any composition $\overline{a}$ of $n$, the disjoint union $\poset{a}$ of chains has a bijection
  $$
  \begin{array}{rcl}
    \linext{\poset{a}} &\longrightarrow& \PTP(\poset{a}) \\
                    \sigma &\longmapsto& \tau
  \end{array}
  $$
  with $\cyc(\tau)=\fc$, the number of prime cycles
in Foata's unique decomposition for $\sigma$.  Thus
$$
\Poin(\poset{a},t)=
\sum_{\sigma \in \linext{\poset{a}}} t^{n-\fc}.
$$
\end{thm}

\noindent
%We also discuss the connection between our map $\Phi$ in Section~\ref{sect:bijection} and Foata's prime cycle decomposition.
Foata's theory is then used to prove the following generating function, generalizing
equation \eqref{stirling-number-gf}.

\begin{thm}
\label{chains-gf-thm}

For $\ell=1,2,\ldots$, one has
$$
\sum_{\overline{a} \in \{1,2,\ldots\}^\ell}
\Poin(\poset{a},t) \cdot \mathbf{x}^{\overline{a}}
=\frac{1}{
1-\sum_{j=1}^\ell e_j(\mathbf{x}) \cdot (t-1)(2t-1)\cdots ((j-1)t-1)},
$$
where
$
\mathbf{x}^{\overline{a}}:=x_1^{a_1} \cdots x_\ell^{a_\ell}
$
and
$
e_j(\mathbf{x}):=\sum_{1 \leq i_1 <\cdots<i_j\leq \ell} x_{i_1} \cdots x_{i_j}
 $
is the $j^{th}$ elementary symmetric function.
\end{thm}

Section \ref{sect:width 2} then examines posets of \emph{width two}, that is, posets $P$ decomposable as
$P=P_1 \cup P_2$ where the subposets $P_1, P_2$ are \emph{chains} (i.e. totally ordered subsets) inside $P$.  Here the Whitney numbers $\whit{k}{P}$ are interpreted by a descent-like
statistic on $\sigma$ in $\linext{P}$:
$$
\des_{P_1,P_2}(\sigma)
:=\#\{ i \in [n-1]: \sigma_i \in P_2, \,\, \sigma_{i+1}\in P_1, \text{ with } \sigma_i, \sigma_{i+1}\text{ incomparable in }P\}.
$$
\begin{thm}
\label{width-two-thm}
For a width two poset decomposed into two chains as $P=P_1 \cup P_2$, one has
$$
\Poin(P,t) =
\sum_{\sigma \in \linext{P}} t^{\des_{P_1,P_2}(\sigma)}.
$$
\end{thm}

\noindent
Theorem~\ref{width-two-thm} will be used to show (Corollary~\ref{width-2:coro:descents})
that, for certain very special
width two posets $P$, the Poincar\'e polynomial $\Poin(P,t)$ coincides with \edit{the}
\emph{$P$-Eulerian polynomial},
 \begin{equation}
 \label{P-Eulerian-polynomial}
  \sum_{\sigma \in \linext{P}} t^{\des(\sigma)}
\end{equation}
  which counts linear extensions $\sigma$ of $P$ according to their number of (usual) {\it descents}
  $$
  \des(\sigma):=\#\{i \in [n-1]: \sigma_i > \sigma_{i+1}\},
  $$
  assuming that $P$ has been \emph{naturally labeled} in the sense that the identity permutation $\sigma=[1, 2,\dots, n]$ lies in $\linext{P}$.
  In particular, Example~\ref{two-row-skew-example} uses this to deduce that
  for $P=\chain{2} \times \chain{n}$,
the {\it Cartesian product} of chains having sizes $2$ and $n$,
the Whitney numbers $\whit{k}{\chain{2} \times \chain{n}}$
are \emph{Narayana numbers}, counting $2 \times n$ standard tableaux according to their number of descents.    This $P$-Eulerian polynomial has many interpretations, e.g.,
as the \emph{$h$-polynomial} of the order complex
  for the \emph{distributive lattice} $J(P)$ of order ideals in $P$, or of the {\it $P$-partition
  triangulation} of the {\it order polytope} for $P$; see \cite[Proposition 2.1, Proposition 2.2]{reiner-welker} and \cite[Sections 3.4, 3.8, 3.13]{ec1} for more on this.    Unfortunately, in general
  the $P$-Eulerian polynomial differs from the Poincar\'e polynomial $\Poin(P,t)$
  considered here. For example, when $P$ is an antichain with three elements,
  the $P$-Eulerian polynomial is $1+4t+t^2$, while $\Poin(P,t)=1+3t+2t^2$.

  Lastly, we note that the cone $\K_P$ associated to a poset $P$ has appeared in
  many places in the literature, for example, \edit{implicitly
in Stanley's theory of {\it $P$-partitions} \cite[\S 3.15]{ec1} and related work of Gessel
on {\it quasisymmetric functions} \cite{Gessel}, more explicitly in work of Bj\"orner and Wachs \cite{BjornerWachs1, BjornerWachs2, BjornerWachs3}},
and as the {\it ranking COMs} of Bandelt, Chepoi and Knauer \cite{bandelt}. A pointed version of this cone, called an \emph{order cone} is discussed in Beck and Sanyal \cite[Chap. 6]{becksanyal}.  When the underlying Hasse diagram of $P$ is a tree, results on the number of chambers $\#\cham{\K_P}=\#\linext{P}$ in these cones appear in work of K. Saito \cite{saito}.

%%%%%%%%%%%%%%%%%%%%%%%%%%%%%%%%%%%%%%%%%%%%%%%%%%%%%%%%%%%%%%%%%%%%%%%%%%%%%%%

\section{The braid arrangement, its intersection lattice, and its cones}
\label{sect:type A}

%%%%%%%%%%%%%%%%%%%%%%%%%%%%%%%%%%%%%%%%%%%%%%%%%%%%%%%%%%%%%%%%%%%%%%%%%%%%%%%

\subsection{Preliminaries on arrangements}
We begin with some preliminaries on hyperplane arrangements, focusing on braid arrangements.  Good references include \cite{OrlikTerao}, \cite{stanley}, \cite[\S 3.11]{ec1}, \cite[\S 3.3]{postnikov}.

\begin{defn}
A \emph{hyperplane} in $V=\R^n$ is an affine linear subspace of codimension one. An \emph{arrangement} of hyperplanes in $\R^n$ is a finite collection $\A=\{H_1,\dots, H_m\}$ of distinct hyperplanes.
A \emph{chamber} of $\A$ is an open, connected component of $\R^n\backslash \bigcup_{H\in\A} H$. The set of all chambers of $\A$ is denoted by $\cham{\A}$.
\end{defn}

\begin{ex}\label{A_n basics}
  The \emph{type A reflection arrangement}, $A_{n-1}$, also called the braid arrangement,
  consists of the $\binom{n}{2}$ hyperplanes of the form
  \[
  H_{ij}=\{(x_1,\dots,x_n)\in\R^n\mid x_i-x_j=0\}
  \]
  for integers \edit{$1\leq i< j\leq n$}. There are $n!$ chambers $\K_\sigma=\{x \in \R^n: x_{\sigma_1} < \cdots < x_{\sigma_n}\}$ of $A_{n-1}$, naturally indexed by the permutations $\sigma=[\sigma_1, \sigma_2, \ldots, \sigma_n]$ of $[n]$, that give the strict inequalities ordering the coordinates within the chamber, as in \eqref{braid-arrangement-chambers}.  For example, when $n=4$,
  $$
  \begin{aligned}
  \K_{1243}&=\{x \in \R^4: x_1<x_2<x_4<x_3\},\\
  \K_{4213}&=\{x \in \R^4: x_4<x_2<x_1<x_3\}
  \end{aligned}
  $$
are two out of the $4!=24$ chambers of $\cham{A_{4-1}}$.
\end{ex}

\begin{defn}
  Let $\A$ be a hyperplane arrangement in $\R^n$. An \emph{intersection} of $\A$ is a nonempty subspace of the form
  $
  X = H_{i_1} \cap H_{i_2} \cap \cdots \cap H_{i_k}
$
where $\{H_{i_1},H_{i_2},\dots, H_{i_k}\}\subseteq \mathcal{A}$.
Here the ambient vector space $V=\R^n$ is considered to be the intersection $\bigcap_{H \in \varnothing} H$ of
the empty set of hyperplanes.
We denote the set of intersections of $\A$ by $\hypinter{\A}$.
\end{defn}

\begin{ex}\label{A_n basics3}
  The intersections of $A_{n-1}$ are described by equalities between the variables.
  \begin{itemize}
  \item For all $n\geq 1$, the line $x_1=x_2=\dots=x_n$
  is the intersection of all the hyperplanes of $A_{n-1}$.
  \item When $n=4$ the intersection of $H_{12}$ and $H_{34}$ in the subspace of $\R^4$ in which $x_1=x_2$ and $x_3=x_4$. On the other hand, the intersection of $H_{12}$ and $H_{13}$ is the subspace of $\R^4$ in which
  $x_1=x_2=x_3$.
  \end{itemize}
  More generally, there is a bijection $\pi \mapsto X_\pi$ between the collection $\Pi_n$ of all \emph{set partitions} $\pi=\{B_1,\dots,B_k\}$ of $[n]=\{1,2,\dots,n\}$ and the set of all intersections of $A_{n-1}$.  The bijection sends
  the set partition $\pi$ to
  the subspace $X_\pi$ where one has equal coordinates
  $x_i=x_j$ whenever $i,j$ lie in a common block $B_k$ of $\pi$.
  We sometimes denote the set partition $\pi=\{B_1,\dots,B_k\}$ with the notation $\pi=B_1 | B_2 | \cdots | B_k$, and may or may not include commas and set braces around the elements of each block $B_i$. E.g., $1\mid 23 \mid 456$ and
  $\{ \{1\},\{2,3\},\{4,5,6\}\}$ represent the same set partition of $[6]$.

  \begin{itemize}
  \item For example, the set partition $1|2|\cdots |n$ in which all elements appear as singletons corresponds to $X_{1|2|\cdots|n}=V=\R^n$, the empty intersection, which is the ambient space.
  \item For all $n\geq 1$,  the set partition $123\cdots n$ having all the elements in the same block corresponds to
  the line $X_{123\cdots n}$ defined by $x_1=x_2=\dots=x_n$.
  \item When $n=4$, one has $ X_{12|34}=H_{12} \cap H_{34}$ and $ X_{123|4}=H_{12} \cap H_{13}$.
  \end{itemize}
\end{ex}

The collection $\hypinter{\A}$ of all intersections of an arrangement $\A$
will be partially ordered by reverse inclusion, and called the \emph{intersection poset} of $\A$.
It has a unique minimal element, namely the intersection
\[
\bigcap_{H \in \varnothing}H=V=\R^n.
\]

Let $\Pi_n$ denote the lattice of set partitions on $[n]$, ordered via refinement:  $\pi_1 \leq \pi_2$
  if $\pi_1$ refines $\pi_2.$  For the braid arrangement $A_{n-1}$, the intersection poset $\hypinter{A_{n-1}}$ is
  isomorphic to $\Pi_n$.

\begin{prop}[{{\cite[pp. 26-27]{stanley}}}] \label{arrangements:thm:isomorphism}
  The map $\pi \mapsto X_\pi$  in Example \ref{A_n basics3} is a poset isomorphism
 $
   \Pi_n \cong \hypinter{A_{n-1}}.
$
\end{prop}

For any hyperplane arrangement $\A$, each of the \emph{lower intervals}
$
[V,X]:=\{Y \in \hypinter{\A}: V \leq Y \leq X\}
$
forms a \emph{geometric lattice} \cite[Prop. 3.11.2]{ec1}.
In particular, this implies that each such lower interval is a \emph{ranked poset}, with rank function given by
the \emph{codimension} $\codim{X}=\dim(V)-\dim(X)$.  Furthermore,
this implies that its \emph{M\"obius function} values $\mu(V,X)$,
defined recursively by
$\mu(V,V):=1$ and $\mu(V,X):=-\sum_{Y: V \leq Y < X} \mu(V,Y)$,
will \emph{alternate in sign} in the sense that $(-1)^{\codim{X}} \mu(V,X) \geq 0$.

For the braid arrangement, these M\"obius function values have a simple expression.

\begin{prop}[{{\cite[Example 3.10.4]{ec1}}}]
\label{cones:prop:mobius isom}
For any set partition $\pi=\{B_1,\dots,B_k\}$ in $\Pi_n$, one has
\[
\mu(V,X_\pi)
= (-1)^{n-k}\prod_{i=1}^k (\#B_i - 1)!
\quad
\left(\,\, =\mu(\,\, 1|2|\cdots|n \, , \, \pi \,\,) \,\, \right)
\]
with the convention $0!:=1$.
Here $\mu(V,X_\pi),\mu(1|2|\cdots|n,\pi)$ are $\mu(-,-)$ values in $\hypinter{A_{n-1}}, \Pi_n$, respectively.
\end{prop}

For the sake of our later discussion, we point out a re-interpretation of this formula involving permutations.
Given a permutation $\sigma$ in $\SN$, when considered as acting on $V=\R^n$, its \emph{fixed space}
$
V^\sigma:=\{\mathbf{x} \in \R^n:\sigma(\mathbf{x})=\mathbf{x}\}
$
will be the intersection subspace $V^\sigma=X_\pi$,
where $\pi=\{B_1,B_2,\ldots,B_k\}$ is the set partition in $\Pi_n$ given by the {\it cycles} of $\sigma$.
Consequently, one can re-interpret

\begin{equation}
\label{type-A-mobius-interpretation}
\begin{aligned}
|\mu(V,X_\pi)|
&= \prod_{i=1}^k (\#B_i - 1)! \\
&=\#\{\sigma \in \SN \text{ with cycle partition }\pi \}\\
&=\#\{\sigma \in \SN: V^\sigma=X_\pi\},
\end{aligned}
\end{equation}
since each block $B_i$ of $\pi$ has $(\#B_i-1)!$ choices of a cyclic orientation.

\begin{defn}
Let $\A$ be an arrangement of hyperplanes in $\R^n$. For $0\leq k\leq n$, the
\emph{$k$th signless Whitney number} of $\hypinter{\A}$ of the first kind is
\[
  \whit{k}{\A} = \sum_{\substack{X\in \hypinter{\A}:\\{\codim{X}=k}}} |\mu(V,X)|
  =(-1)^k \sum_{\substack{X\in \hypinter{\A}:\\{\codim{X}=k}}} \mu(V,X).
\]
\end{defn}

\noindent
Henceforth, we call $c_k(\A)$, $0\le k\le n$, the \emph{Whitney numbers of $\A$}.
One of the standard ways to compile them into a generating function is their \emph{Poincar\'e polynomial}
$\Poin(\A,t):=\sum_{k=0}^n\whit{k}{\A} t^k$; see \cite[\S 2.3]{OrlikTerao}.

As mentioned in the Introduction, we aim to understand the chambers, intersections, and Whitney numbers for cones in $\A$;   the chambers, intersections, and Whitney numbers for $\A$ are a special case.

\begin{defn}
Let $\A$ be an arrangement of hyperplanes in $V=\R^n$.
A \emph{cone}\footnote{Aguiar and Mahajan \cite{aguiar} call these objects \emph{top-cones}.} $\K$ of $\A$ is any nonempty
intersection $\varnothing \neq \K \subseteq V=\R^n$ of (open) halfspaces defined by a subset $\A'$ of the hyperplanes from $\A$.  That is, a cone $\K$ is any one of the (open) chambers
from the set of all chambers $\cham{\A'}$ for some subarrangement $\A' \subseteq \A$. For example, in the following arrangement in $\R^2$ there are four cones defined by the dashed hyperplanes. One such cone $\K$ is shaded below.

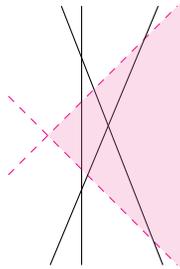
\begin{center}
 \begin {tikzpicture}[scale=.6]

  \node[invisivertex] (A1) at (-1,1){};
  \node[invisivertex] (A2) at (3,-3){};
  \node[invisivertex] (B1) at (.75,-3){};
  \node[invisivertex] (B2) at (.75,3){};
  \node[invisivertex] (C1) at (2.6,-3){};
  \node[invisivertex] (C2) at (.25,3){};
  \node[invisivertex] (D1) at (0,-3){};
  \node[invisivertex] (D2) at (2.5,3){};
  \node[invisivertex] (E1) at (-1,-1){};
  \node[invisivertex] (E2) at (3,3){};

  % \node[invisivertex] (chamber1) at (.5,0){\textcolor{magenta}{$1$}};
  % \node[invisivertex] (chamber2) at (1.15,0){\textcolor{magenta}{$2$}};
  % \node[invisivertex] (chamber3) at (1.5,1){\textcolor{magenta}{$3$}};
  % \node[invisivertex] (chamber4) at (1.5,-1){\textcolor{magenta}{$4$}};
  % \node[invisivertex] (chamber5) at (2.5,0){\textcolor{magenta}{$5$}};

  \path[-] (A1) edge [bend left =0, color=magenta,dashed] node[above] {} (A2);
  \path[-] (E1) edge [bend left =0, color=magenta,dashed] node[above] {} (E2);
  \path[-] (B1) edge [bend left =0] node[above] {} (B2);
  \path[-] (C1) edge [bend left =0] node[above] {} (C2);
  \path[-] (D1) edge [bend left =0] node[above] {} (D2);

  \fill[magenta!70,nearly transparent] (0,0) -- (3,3) -- (3,-3) -- cycle;
  \end{tikzpicture}
\end{center}

\noindent Each cone $\K$ of $\A$ has its collection of \emph{chambers}, namely those chambers in
$\cham{\A}$ that lie inside $\K$:
\[
\cham{\K} = \{C\in\cham{\A}:C \subset \K\}.
\]
The poset of \emph{interior intersections} of the cone $\K$  is the following order ideal within the poset $\hypinter{\A}$:
\[
\mathcal{L}^{\text{int}}(\K) =\mathcal{L}_{\A}^{\text{int}}(\K) = \{X\in\hypinter{\A} \mid X \cap
\K\not=\emptyset\}.
\]
For each $X$ in $\mathcal{L}_{\A}^{\text{int}}(\K)$, its lower interval $[V,X]$ is still a geometric lattice, with same rank function
$\codim{X}$, so that one can define the $k^{th}$ \emph{(signless) Whitney number of $\K$} by
  \[
    \whit{k}{\K} = \sum_{\substack{X\in\inter{\K}:\\\codim{X}=k}}|\mu(V,X)|
     =(-1)^k \sum_{\substack{X\in \inter{\K}:\\{\codim{X}=k}}} \mu(V,X),
  \]
along with their generating function $\Poin(\K,t):= \sum_{k=0}^n \whit{k}{\K} t^k$,
the \emph{Poincar\'e polynomial} for $\K$.
 \end{defn}

The starting point for our study is the following result of Zaslavsky~\cite{zaslavsky} counting
the number $\#\cham{\K}$ of chambers of an arrangement $\A$ lying inside one of its cones $\K$.

\begin{thm}[{\cite[Example A, p. 275]{zaslavsky}}]
\label{zaslavsky's theorem}
   Let $\K$ be a cone of an arrangement $\A$ in $V=\R^n$. Then
   \[
     \#
     \cham{\K}
     =\sum_{X \in \inter{\K}} |\mu(V,K)|
     = \whit{0}{\K} + \whit{1}{\K} + \cdots +\whit{n}{\K}
     =\left[\Poin(\K,t)\right]_{t=1}.
   \]
\end{thm}

\noindent
Zaslavsky proved in his doctoral thesis \cite{zaslavsky-memoir} the better-known special case of Theorem \ref{zaslavsky's theorem} for the full arrangement, that is, where $\K=V=\R^n$.

The following two examples illustrate Theorem~\ref{zaslavsky's theorem} for two cones in $A_3$.

\begin{ex}\label{arrangements:ex:3<4}
Consider the braid arrangement $\A=A_4=\{H_{12},H_{13},H_{14},H_{23},H_{24},H_{34}\}$ inside $V=\R^4$.
On the left below we have drawn a linearly equivalent picture of
its intersection with the hyperplane where $x_1+x_2+x_3+x_4=0$,
isomorphic to $\R^3$, and depicted the intersection of the hyperplanes with the unit $2$-sphere in this $3$-dimensional
space.  Here we pick the cone $\K$ to be the one defined by the halfspace $x_3<x_4$ for the hyperplane $H_{34}$,
and draw the intersection of $H_{34}$ with the unit sphere as the equatorial circle,
with the other five hyperplanes $H_{ij}$ depicted as great circles intersecting the hemisphere where
$x_3<x_4$.  On the right below the non-hyperplane interior intersection subspaces $X_\pi$ are labeled.

  
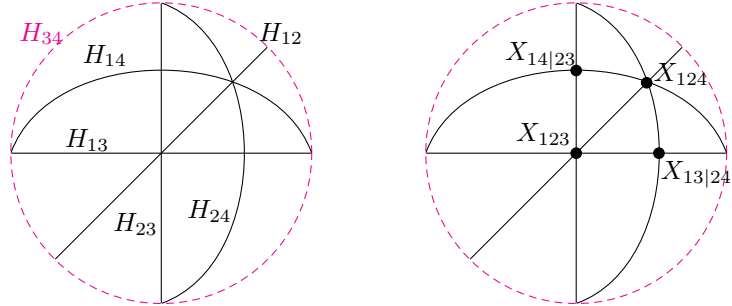
\begin{center}
    \begin{tikzpicture}[scale=.5]
    \draw[boundaryEdge] (0,0) circle (4cm);
    %Boundary Vertices
    \node[origVertex]  (14A) at (180:4) {};
    % \node[origVertex]  (34A) at (65:4) {};
    \node[origVertex] (13A) at (180:4) {};
    \node[origVertex]  (23A) at (90:4) {};
    \node[origVertex]  (12A) at (45:4) {};
    \node[origVertex]  (24A) at (90:4) {};
    \node[origVertex]  (14B) at (0:4) {};
    % \node[origVertex]  (34B) at (240:4) {};
    \node[origVertex] (13B) at (0:4) {};
    \node[origVertex]  (23B) at (270:4) {};
    \node[origVertex]  (12B) at (225:4) {};
    \node[origVertex]  (24B) at (270:4) {};

    %Hyperplane Labels
    \node[invisivertex] (H14) at (120:3) {$H_{14}$};
    \node[invisivertex] (H13) at (170:2) {$H_{13}$};
    \node[invisivertex] (H23) at (250:2) {$H_{23}$};
    \node[invisivertex] (H12) at (45:4.5) {$H_{12}$};
    \node[invisivertex] (H24) at (310:2) {$H_{24}$};
    \node[invisivertex] (H34) at (135:4.5) {\textcolor{magenta}{$H_{34}$}};

    % % Intersection Markers
    % \node[point] (13|24) at (0:2.2) {};
    % \node[point] (124) at (45:2.65) {};
    % \node[point] (14,23) at (90:2.2) {};
    % \node[point] (123) at (0:0) {};
    % \node[point] (234+) at (90:4) {};
    % \node[point] (234-) at (270:4) {};
    % \node[point] (12|34+) at (45:4) {};
    % \node[point] (12|34-) at (225:4) {};
    % \node[point] (134+) at (180:4) {};
    % \node[point] (134-) at (0:4) {};
    %
    % % Intersection Labels
    % \node[invisivertex]  (X 13|24) at (5:3) {$X_{13|24}$};
    % \node[invisivertex]  (X 124) at (37:3.25) {$X_{124}$};
    % \node[invisivertex]  (X 14|23) at (107:2.75) {$X_{14|23}$};
    % \node[invisivertex]  (X 123) at (150:.5) {$X_{123}$};
    % \node[invisivertex] (X 234+) at (90:4.5) {$X_{234}$};
    % \node[invisivertex] (X 234-) at (270:4.5) {$X_{234}$};
    % \node[invisivertex] (X 12|34+) at (45:4.5) {$X_{12|34}$};
    % \node[invisivertex] (X 12|34-) at (225:4.5) {$X_{12|34}$};
    % \node[invisivertex] (X 134+) at (180:4.5) {$X_{134}$};
    % \node[invisivertex] (X 134-) at (0:4.5) {$X_{134}$};

    % Hyperplanes (edges)

    \path[-] (13B) edge[-] (13A);
    \path[-] (14B) edge[out=100, in=200, bend right=70] (14A);
    \path[-] (23B) edge[] (23A);
    \path[-] (12B) edge[] (12A);
    \path[-] (24B) edge[out=220, in=340, bend right=70] (24A);

    % Shading
    %\fill[magenta!40,nearly transparent] (132:1.75) -- (180:1.8) -- (200:2.2)  -- (220:2.85)-- (240:4) -- (255:4) -- (270:4) -- (280:4)-- (295:4)-- (310:4) --  (320:4) -- cycle;
    \end{tikzpicture}
\hspace{.5in}
    \begin{tikzpicture}[scale=.5]
    \draw[boundaryEdge] (0,0) circle (4cm);
    %Boundary Vertices
    \node[origVertex]  (14A) at (180:4) {};
    % \node[origVertex]  (34A) at (65:4) {};
    \node[origVertex] (13A) at (180:4) {};
    \node[origVertex]  (23A) at (90:4) {};
    \node[origVertex]  (12A) at (45:4) {};
    \node[origVertex]  (24A) at (90:4) {};
    \node[origVertex]  (14B) at (0:4) {};
    % \node[origVertex]  (34B) at (240:4) {};
    \node[origVertex] (13B) at (0:4) {};
    \node[origVertex]  (23B) at (270:4) {};
    \node[origVertex]  (12B) at (225:4) {};
    \node[origVertex]  (24B) at (270:4) {};

    %Hyperplane Labels
    % \node[invisivertex] (H14) at (8:3.5) {$H_{14}$};
    % \node[invisivertex] (H13) at (10:1.5) {$H_{13}$};
    % \node[invisivertex] (H23) at (255:2) {$H_{23}$};
    % \node[invisivertex] (H12) at (45:4.5) {$H_{12}$};
    % \node[invisivertex] (H24) at (325:2) {$H_{24}$};
    % \node[invisivertex] (H24) at (135:4.5) {$H_{34}$};

    % % Intersection Markers
    \node[point] (13|24) at (0:2.2) {};
    \node[point] (124) at (45:2.65) {};
    \node[point] (14,23) at (90:2.2) {};
    \node[point] (123) at (0:0) {};
    % \node[point] (234+) at (90:4) {};
    % \node[point] (234-) at (270:4) {};
    % \node[point] (12|34+) at (45:4) {};
    % \node[point] (12|34-) at (225:4) {};
    % \node[point] (134+) at (180:4) {};
    % \node[point] (134-) at (0:4) {};
    %
    % % Intersection Labels
    \node[invisivertex]  (X 13|24) at (350:3.25) {$X_{13|24}$};
    \node[invisivertex]  (X 124) at (37:3.5) {$X_{124}$};
    \node[invisivertex]  (X 14|23) at (110:2.75) {$X_{14|23}$};
    \node[invisivertex]  (X 123) at (150:1) {$X_{123}$};
    % \node[invisivertex] (X 234+) at (90:4.5) {$X_{234}$};
    % \node[invisivertex] (X 234-) at (270:4.5) {$X_{234}$};
    % \node[invisivertex] (X 12|34+) at (45:4.5) {$X_{12|34}$};
    % \node[invisivertex] (X 12|34-) at (225:4.5) {$X_{12|34}$};
    % \node[invisivertex] (X 134+) at (180:4.5) {$X_{134}$};
    % \node[invisivertex] (X 134-) at (0:4.5) {$X_{134}$};

    % Hyperplanes (edges)
    \path[-] (13B) edge[-] (13A);
    \path[-] (14B) edge[out=100, in=200, bend right=70] (14A);
    \path[-] (23B) edge[] (23A);
    \path[-] (12B) edge[] (12A);
    \path[-] (24B) edge[out=220, in=340, bend right=70] (24A);

    % Shading
    %\fill[magenta!40,nearly transparent] (132:1.75) -- (180:1.8) -- (200:2.2)  -- (220:2.85)-- (240:4) -- (255:4) -- (270:4) -- (280:4)-- (295:4)-- (310:4) --  (320:4) -- cycle;
    \end{tikzpicture}
\end{center}
Therefore the intersection poset $\inter{\K}$ of this cone is
    
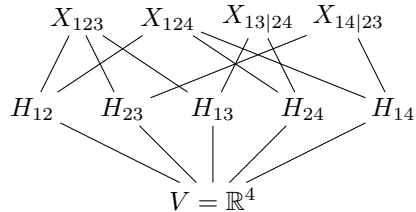
\begin{center}
    \begin{tikzpicture}[scale=.6]

    % Lattice Elements, Grouped by Rank

    \node[invisivertex] (0) at (5,0){$V=\R^4$};

    \node[invisivertex] (H12) at (1,2){$H_{12}$};
    \node[invisivertex] (H23) at (3,2){$H_{23}$};
    \node[invisivertex] (H13) at (5,2){$H_{13}$};
    \node[invisivertex] (H24) at (7,2){$H_{24}$};
    \node[invisivertex] (H14) at (9,2){$H_{14}$};
    % \node[invisivertex] (H34) at (10,2){$H_{34}$};

    \node[invisivertex] (X123) at (2,4){$X_{123}$};
    \node[invisivertex] (X124) at (4,4){$X_{124}$};
    % \node[invisivertex] (X12 34) at (3,4){$X_{12|34}$};
    \node[invisivertex] (X13 24) at (6,4){$X_{13|24}$};
    \node[invisivertex] (X14 23) at (8,4){$X_{14|23}$};
    % \node[invisivertex] (X134) at (9,4){$X_{134}$};
    % \node[invisivertex] (X234) at (11,4){$X_{234}$};

    % \node[invisivertex] (1) at (5,6){$\hat{1}=\text{span}(\vec{1})$};

    % Order Relations, Grouped by Rank

    \path[-] (H12) edge [bend left =0] node[above] {} (0);
    \path[-] (H23) edge [bend left =0] node[above] {} (0);
    \path[-] (H13) edge [bend left =0] node[above] {} (0);
    \path[-] (H24) edge [bend left =0] node[above] {} (0);
    \path[-] (H14) edge [bend left =0] node[above] {} (0);
    % \path[-] (H34) edge [bend left =0] node[above] {} (0);

    \path[-] (H12) edge [bend left =0] node[above] {} (X123);
    \path[-] (H12) edge [bend left =0] node[above] {} (X124);
    % \path[-] (H12) edge [bend left =0] node[above] {} (X12 34);
    \path[-] (H23) edge [bend left =0] node[above] {} (X123);
    \path[-] (H23) edge [bend left =0] node[above] {} (X14 23);
    % \path[-] (H23) edge [bend left =0] node[above] {} (X234);
    \path[-] (H13) edge [bend left =0] node[above] {} (X123);
    \path[-] (H13) edge [bend left =0] node[above] {} (X13 24);
    % \path[-] (H13) edge [bend left =0] node[above] {} (X134);
    \path[-] (H24) edge [bend left =0] node[above] {} (X124);
    \path[-] (H24) edge [bend left =0] node[above] {} (X13 24);
    % \path[-] (H24) edge [bend left =0] node[above] {} (X234);
    \path[-] (H14) edge [bend left =0] node[above] {} (X124);
    \path[-] (H14) edge [bend left =0] node[above] {} (X14 23);
    % \path[-] (H14) edge [bend left =0] node[above] {} (X134);
    % \path[-] (H34) edge [bend left =0] node[above] {} (X134);
    % \path[-] (H34) edge [bend left =0] node[above] {} (X234);
    % \path[-] (H34) edge [bend left =0] node[above] {} (X12 34);

    % \path[-] (X123) edge [bend left =0] node[above] {} (1);
    % \path[-] (X124) edge [bend left =0] node[above] {} (1);
    % \path[-] (X12 34) edge [bend left =0] node[above] {} (1);
    % \path[-] (X13 24) edge [bend left =0] node[above] {} (1);
    % \path[-] (X14 23) edge [bend left =0] node[above] {} (1);
    % \path[-] (X134) edge [bend left =0] node[above] {} (1);
    % \path[-] (X234) edge [bend left =0] node[above] {} (1);

    %\fill[magenta!40,nearly transparent] (0,0) -- (3,3) -- (3,-3) -- cycle;
    \end{tikzpicture}
  \end{center}

  We have $(\whit{0}{\K},\whit{1}{\K},\whit{2}{\K}) = (1,5,6)$.
  Summing these gives
  \(
  1+5+6=12,
  \)
  and a quick visual verification
  assures that there are $12$ chambers in this cone.
\end{ex}

\begin{ex}\label{1<2 and 3<4}
Consider the cone $\K$ of $A_3$ in which $x_3<x_4$ and $x_1<x_2$.
On the left below we have drawn the same picture as Example \ref{arrangements:ex:3<4} with
the cone corresponding to $\K$ shaded. We depict $\inter{\K}$ on the right.

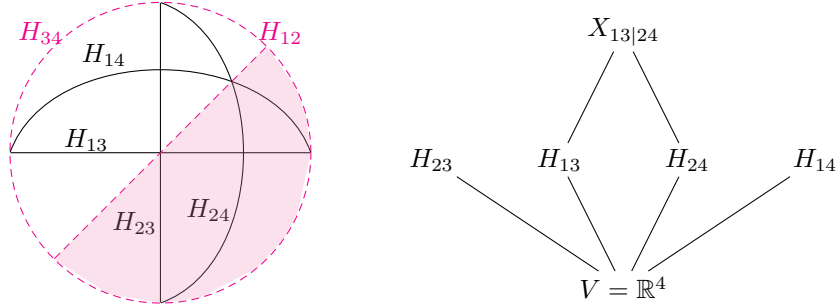
\begin{center}
  \begin{tikzpicture}[scale=.5]
  \draw[boundaryEdge] (0,0) circle (4cm);
  %\draw[dashed] (zero) circle (4cm) {};

  %Boundary Vertices
  \node[origVertex]  (14A) at (180:4) {};
  % \node[origVertex]  (34A) at (65:4) {};
  \node[origVertex] (13A) at (180:4) {};
  \node[origVertex]  (23A) at (90:4) {};
  \node[origVertex]  (12A) at (45:4) {};
  \node[origVertex]  (24A) at (90:4) {};
  \node[origVertex]  (14B) at (0:4) {};
  % \node[origVertex]  (34B) at (240:4) {};
  \node[origVertex] (13B) at (0:4) {};
  \node[origVertex]  (23B) at (270:4) {};
  \node[origVertex]  (12B) at (225:4) {};
  \node[origVertex]  (24B) at (270:4) {};

  %Hyperplane Labels
  \node[invisivertex] (H14) at (120:3) {$H_{14}$};
  \node[invisivertex] (H13) at (170:2) {$H_{13}$};
  \node[invisivertex] (H23) at (250:2) {$H_{23}$};
  \node[invisivertex] (H24) at (310:2) {$H_{24}$};
  \node[invisivertex] (H12) at (45:4.5) {\textcolor{magenta}{$H_{12}$}};
  \node[invisivertex] (H24) at (135:4.5) {\textcolor{magenta}{$H_{34}$}};

  % % Intersection Markers
  % \node[point] (124) at (45:2.65) {};
  % \node[point] (14,23) at (90:2.2) {};
  % \node[point] (123) at (0:0) {};
  % \node[point] (234+) at (90:4) {};
  % \node[point] (234-) at (270:4) {};
  % \node[point] (12|34+) at (45:4) {};
  % \node[point] (12|34-) at (225:4) {};
  % \node[point] (134+) at (180:4) {};
  % \node[point] (134-) at (0:4) {};
  %
  % % Intersection Labels
  % \node[invisivertex]  (X 124) at (37:3.25) {$X_{124}$};
  % \node[invisivertex]  (X 14|23) at (107:2.75) {$X_{14|23}$};
  % \node[invisivertex]  (X 123) at (150:.5) {$X_{123}$};
  % \node[invisivertex] (X 234+) at (90:4.5) {$X_{234}$};
  % \node[invisivertex] (X 234-) at (270:4.5) {$X_{234}$};
  % \node[invisivertex] (X 12|34+) at (45:4.5) {$X_{12|34}$};
  % \node[invisivertex] (X 12|34-) at (225:4.5) {$X_{12|34}$};
  % \node[invisivertex] (X 134+) at (180:4.5) {$X_{134}$};
  % \node[invisivertex] (X 134-) at (0:4.5) {$X_{134}$};

  % Hyperplanes (edges)
  \path[-] (13B) edge[-] (13A);
  \path[-] (14B) edge[out=100, in=200, bend right=70] (14A);
  \path[-] (23B) edge[] (23A);
  \path[-] (12B) edge[boundaryEdge] (12A);
  \path[-] (24B) edge[out=220, in=340, bend right=70] (24A);

  % Shading
  \fill[magenta!60,nearly transparent] (225:4) -- (45:4) -- (25:4) -- (15:4) -- (0:4) -- (335:4)
  -- (325:4) -- (310:4) -- (285:4) -- (270:4) -- (250:4) -- cycle;
  \end{tikzpicture}
  \hspace{1cm}
% \end{center}
% The intersection poset of this cone is
%     \begin{center}
    \begin{tikzpicture}[scale=.85]

    % Lattice Elements, Grouped by Rank

    \node[invisivertex] (0) at (5,0){$V=\R^4$};

    % \node[invisivertex] (H12) at (1,2){$H_{12}$};
    \node[invisivertex] (H23) at (2,2){$H_{23}$};
    \node[invisivertex] (H13) at (4,2){$H_{13}$};
    \node[invisivertex] (H24) at (6,2){$H_{24}$};
    \node[invisivertex] (H14) at (8,2){$H_{14}$};
    % \node[invisivertex] (H34) at (10,2){$H_{34}$};

    % \node[invisivertex] (X123) at (2,4){$X_{123}$};
    % \node[invisivertex] (X124) at (4,4){$X_{124}$};
    % \node[invisivertex] (X12 34) at (3,4){$X_{12|34}$};
    \node[invisivertex] (X13 24) at (5,4){$X_{13|24}$};
    % \node[invisivertex] (X14 23) at (8,4){$X_{14|23}$};
    % \node[invisivertex] (X134) at (9,4){$X_{134}$};
    % \node[invisivertex] (X234) at (11,4){$X_{234}$};

    % \node[invisivertex] (1) at (5,6){$\hat{1}=\text{span}(\vec{1})$};

    % Order Relations, Grouped by Rank

    % \path[-] (H12) edge [bend left =0] node[above] {} (0);
    \path[-] (H23) edge [bend left =0] node[above] {} (0);
    \path[-] (H13) edge [bend left =0] node[above] {} (0);
    \path[-] (H24) edge [bend left =0] node[above] {} (0);
    \path[-] (H14) edge [bend left =0] node[above] {} (0);
    % \path[-] (H34) edge [bend left =0] node[above] {} (0);

    % \path[-] (H12) edge [bend left =0] node[above] {} (X123);
    % \path[-] (H12) edge [bend left =0] node[above] {} (X124);
    % \path[-] (H12) edge [bend left =0] node[above] {} (X12 34);
    % \path[-] (H23) edge [bend left =0] node[above] {} (X123);
    % \path[-] (H23) edge [bend left =0] node[above] {} (X14 23);
    % \path[-] (H23) edge [bend left =0] node[above] {} (X234);
    % \path[-] (H13) edge [bend left =0] node[above] {} (X123);
    \path[-] (H13) edge [bend left =0] node[above] {} (X13 24);
    % \path[-] (H13) edge [bend left =0] node[above] {} (X134);
    % \path[-] (H24) edge [bend left =0] node[above] {} (X124);
    \path[-] (H24) edge [bend left =0] node[above] {} (X13 24);
    % \path[-] (H24) edge [bend left =0] node[above] {} (X234);
    % \path[-] (H14) edge [bend left =0] node[above] {} (X124);
    % \path[-] (H14) edge [bend left =0] node[above] {} (X14 23);
    % \path[-] (H14) edge [bend left =0] node[above] {} (X134);
    % \path[-] (H34) edge [bend left =0] node[above] {} (X134);
    % \path[-] (H34) edge [bend left =0] node[above] {} (X234);
    % \path[-] (H34) edge [bend left =0] node[above] {} (X12 34);

    % \path[-] (X123) edge [bend left =0] node[above] {} (1);
    % \path[-] (X124) edge [bend left =0] node[above] {} (1);
    % \path[-] (X12 34) edge [bend left =0] node[above] {} (1);
    % \path[-] (X13 24) edge [bend left =0] node[above] {} (1);
    % \path[-] (X14 23) edge [bend left =0] node[above] {} (1);
    % \path[-] (X134) edge [bend left =0] node[above] {} (1);
    % \path[-] (X234) edge [bend left =0] node[above] {} (1);

    %\fill[magenta!40,nearly transparent] (0,0) -- (3,3) -- (3,-3) -- cycle;
    \end{tikzpicture}
  \end{center}

  We have $\whit{0}{\K} = 1$, $\whit{1}{\K} = 4$, and $\whit{2}{\K} = 1$. Summing
  these gives
  \(
  1+4+1=6 = \#\cham{\K}.
  \)
\end{ex}

For the remainder of this paper, we focus on cones $\K$ inside braid arrangements $A_{n-1}$.
It is well-known (see \cite[\S 3.3]{postnikov}, for example) and easy to see
that such cones correspond bijectively with posets $P$ on $[n]$ via this rule:
one has $x_i<x_j$ for all points in the cone $\K$ if and only if $i<_P j$.
We will denote the cone associated to $P$ by $\K_P$, and abbreviate
$\whit{k}{P}:=\whit{k}{\K_P}$ and $\Poin(P,t):=\Poin(\K_P,t)$.

\begin{ex}
\label{two-union-of-chain-examples}
The cone inside $A_3$ in Example~\ref{arrangements:ex:3<4} given by the inequality $x_3 < x_4$ on $V=\R^4$
has defining poset $P_1$ with order relation $3<_{P_1} 4$ on $[4]=\{1,2,3,4\}$, while the cone in Example~\ref{1<2 and 3<4} given by the inequalities $x_1<x_2$ and $x_3<x_4$ has defining poset $P_2$ with
  order relations $1<_{P_2} 2$ and $3<_{P_2} 4$.  These posets $P_1, P_2$ are shown here:

\begin{center}
\begin{tikzpicture}[scale=.8]
  % \node[invisivertex] (A) at (-2,.5){$\leftrightarrow$};
  \node[invisivertex] (P) at (-1,.5){$P_1=$};
  \node[circledvertex] (1) at (0,0){$1$};
  \node[circledvertex] (2) at (.75,0){$2$};
  \node[circledvertex] (3) at (1.5,0){$3$};
  \node[circledvertex] (4) at (1.5,.75){$4$};

%   \path[-] (1) edge [bend left =0] node[above] {} (2);
  \path[-] (3) edge [bend left =0] node[above] {} (4);
\end{tikzpicture}
\hspace{2cm}
\begin{tikzpicture}[scale=.8]
  % \node[invisivertex] (A) at (-2,.5){$\leftrightarrow$};
  \node[invisivertex] (P) at (-1,.5){$P_2=$};
  \node[circledvertex] (1) at (0,0){$1$};
  \node[circledvertex] (2) at (0,.75){$2$};
  \node[circledvertex] (3) at (.75,0){$3$};
  \node[circledvertex] (4) at (.75,.75){$4$};

  \path[-] (1) edge [bend left =0] node[above] {} (2);
  \path[-] (3) edge [bend left =0] node[above] {} (4);
\end{tikzpicture}
\end{center}

\end{ex}

\subsection{Preposets, posets, cones and a characterization of $\inter{\K_P}$}

By Theorem \ref{arrangements:thm:isomorphism}, the intersection poset $\hypinter{A_{n-1}}$ is
isomorphic to the set partition lattice $\Pi_n$, and hence for each cone $\K_P$ in $\A_{n-1}$,
one should be able to identify the interior intersection poset $\inter{\K_P}$ as some order ideal
inside $\Pi_n$.  This is our next goal, which will be aided
by recalling some facts about preposets, posets, binary relations, and cones.

\begin{defn}
Recall that a \emph{preposet} $Q$ on $[n]$ is a binary relation
$Q \subseteq [n] \times [n]$
which is both
\begin{itemize}
\item \emph{reflexive}, meaning $(i,i) \in Q$ for all $i$, and
\item \emph{transitive}, meaning $(i,j),(j,k) \in Q$ implies $(i,k) \in Q$.
\end{itemize}
If in addition $Q$ is \emph{antisymmetric}, meaning $(i,j),(j,i) \in Q$ implies $i=j$,
then $Q$ is called a \emph{poset} on $[n]$;  in this case, we sometimes write $i \leq_Q j$ when $(i,j) \in Q$.

A set partition $\pi\in\Pi_n$ is identified with an \emph{equivalence relation}
$\pi \subseteq [n] \times [n]$ having  $(i,j)\in\pi$ when $i,j$ appear in the same block of $\pi$.
That is, $\pi$ is reflexive, transitive, and \emph{symmetric}, meaning $(i,j) \in \pi$ implies $(j,i) \in \pi$.
We will sometimes write this binary relation as $i \equiv_\pi j$ when $(i,j) \in \pi$.

The \emph{union} $Q_1 \cup Q_2 \subseteq [n] \times [n]$ of two preposets will be
reflexive, but possibly not transitive, so not always a preposet.
However, the transitive closure operation $Q \mapsto \overline{Q}$ lets one
complete it to a preposet $\overline{Q_1\cup Q_2}$.
\end{defn}

We will use a slight rephrasing of the folklore \emph{cone-preposet dictionary}, as discussed by Postnikov, Reiner, and Williams in \cite[Section 3.3]{postnikov}.
This dictionary is a bijection between preposets $Q$ on $[n]$ and \emph{closed} cones of any dimension that are intersections in $V=\R^n$ of closed halfspaces of the form
$\{x_i \leq x_j\}$.  Under this bijection, any such closed cone $C$ corresponds to a preposet $Q_C$ via
\[
C \mapsto
Q_C:=\{(i,j)\mid x_i\leq x_j \text{ for all }\mathbf{x}\in C \}.
\]
Conversely, any preposet $Q$ on $[n]$
corresponds to a closed cone $C_Q$ via
\[
Q\mapsto C_Q:=\bigcap_{(i,j) \in Q} \{x_i \leq x_j\}
=\{\mathbf{x}\in\R^n \mid x_i\leq x_j \text{ for all }(i,j)\in Q\}.
\]

\noindent
For a subset $A\subseteq\R^n$, denote its \emph{interior} and \emph{relative interior} by
$\cint{A}, \relint{A}$. Then for a preposet $Q$,
\begin{equation}
    \label{relative-interior-description}
\relint{C_Q}
      =
      \left\{\mathbf{x}\in\R^n :
      \begin{matrix}
       x_i<x_j\text{ if }(i,j)\in Q \text{ but }(j,i)\not\in Q,\\
        x_i=x_j \text{ if both } (i,j),(j,i)\in Q
    \end{matrix}
      \right\}.
\end{equation}

Also, one has the following assertions, using the notation of this dictionary:
\begin{itemize}
\item for $\pi$ in $\Pi_n$, the subspace denoted $X_\pi$ is the (non-pointed) cone $C_\pi$,
regarding $\pi$ as a preposet, and
    \item
for any poset $P$ on $[n]$, the open $n$-dimensional cone denoted $\K_P$ earlier is $\relint{C_P}$ $( ={\sf{int}}(C_P) )$.
\end{itemize}

We will need one further dictionary fact.
\begin{prop}[{{\cite[Proposition 3.5]{postnikov}}}]\label{cones:prop:cone-preposet dictionary}
For preposets $Q, Q'$, one has
$
C_Q \cap C_{Q'}
=  C_{\overline{Q\cup Q'}}.
$
\end{prop}

The following definition will help to characterize the set partitions $\pi$ having $X_\pi$ in $\inter{\K_P}$.

\begin{defn}
Given a poset $P$ on $[n]$ and a set partition
$\pi=\{B_1,\dots,B_k\}$ in $\Pi_n$, define a preposet $P/\pi$ on the set $\{B_1,\dots, B_k\}$ as the transitive closure of the (reflexive)
binary relation having $(B_i, B_j)\in P/\pi$ whenever there exist $p\in B_i$ and $q\in B_j$ with $p \leq_P q$.
\end{defn}

\begin{prop}\label{cones:prop:transverse}
  For $P$ a poset on $[n]$ and $\pi=\{B_1,\dots,B_k\}$ a set partition in $\Pi_n$, the following are equivalent.
  \begin{enumerate}
    \item[(i)] $X_\pi \in\inter{\K_P}$, that is, one has a nonempty intersection $X_\pi \cap \K_P \neq \varnothing$.
    \item[(ii)] If $i <_P j$, meaning that $i \leq_P j$ and $i \neq j$, then $(j,i)\not\in\overline{P\cup \pi}$.
    \item[(iii)] Every block $B_i\in\pi$ is an antichain of $P$, and the preposet $P/\pi$ is actually a poset.
  \end{enumerate}
\end{prop}

We give a  proof of Proposition \ref{cones:prop:transverse} toward the end of
this subsection, after some discussion and examples.

\begin{defn}
  Let $P$ be a poset on $[n]$. A set partition $\pi$ of $P$ is called a \emph{$P$-transverse partition} if it satisfies one of the equivalent
  conditions in Proposition~\ref{cones:prop:transverse}. We denote by $\PTPar(P)$
  the induced subset of $\Pi_n$ consisting of $P$-transverse partitions.
\end{defn}

\begin{rem}
Aguiar and Mahajan \cite[p.230]{aguiar} have a similar concept, which they call a \emph{prelinear extension} of $P$. A prelinear extension of $P$ is  equivalent to a $P$-transverse partition $\pi$ together with a linear ordering on the blocks of $\pi$ that extends the partial order $P/\pi$ from Proposition~\ref{cones:prop:transverse}(iii).
\end{rem}

Proposition~\ref{cones:prop:transverse} and
Corollary~\ref{cones:prop:mobius isom} immediately imply the following corollary.

\begin{coro}\label{cor:Lint}
  Let $P$ be a poset on $[n]$. Then $\PTPar(P)$ and $\inter{\K_P}$ are isomorphic as posets.
  Consequently,
$$
%\begin{aligned}
    \Poin(P,t)
    =\sum_{\pi \in \PTPar(P)} |\mu(V,X_\pi)| \cdot t^{n-\#\blocks(\pi)}\\
    =\sum_{\substack{\pi=\{B_1,\ldots,B_k\} \\ \text{ in }\PTPar(P)}}
    \quad \prod_{i=1}^k (\#B_i-1)! \cdot t^{n-k}.
%\end{aligned}
$$
\end{coro}

\begin{ex}
Let $P:=P_2$ be the second poset from Example~\ref{two-union-of-chain-examples}, with
$x_1 <_P x_2$ and $x_3 <_P x_4$.   Then
\begin{itemize}
    \item $\pi=13|24$ is $P$-transverse.
    \item $\pi= 12|3|4$ is \emph{not} $P$-transverse as it fails condition (ii): $1<_P 2$, but $(2,1) \in \pi \subset \overline{P\cup \pi}.$
    \item  $\pi=14|23$ is \emph{not} $P$-transverse, failing condition (ii):  $1<_P 2$,
but $(2,1) \in \overline{P\cup \pi}$, though $(2,1) \not\in P\cup \pi$.
\end{itemize}
The six $P$-transverse partitions give this subposet $\PTPar(P)$ of $\Pi_4$ isomorphic to
$\inter{\K_P}$, as in Example~\ref{1<2 and 3<4}:

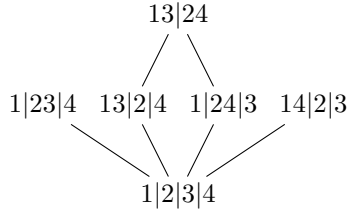
\begin{center}
\begin{tikzpicture}[scale=.6]

% Lattice Elements, Grouped by Rank

\node[invisivertex] (0) at (5,0){$1|2|3|4$};

% \node[invisivertex] (H12) at (1,2){$H_{12}$};
\node[invisivertex] (H23) at (2,2){$1|23|4$};
\node[invisivertex] (H13) at (4,2){$13|2|4$};
\node[invisivertex] (H24) at (6,2){$1|24|3$};
\node[invisivertex] (H14) at (8,2){$14|2|3$};
% \node[invisivertex] (H34) at (10,2){$H_{34}$};

% \node[invisivertex] (X123) at (2,4){$X_{123}$};
% \node[invisivertex] (X124) at (4,4){$X_{124}$};
% \node[invisivertex] (X12 34) at (3,4){$X_{12|34}$};
\node[invisivertex] (X13 24) at (5,4){$13|24$};
% \node[invisivertex] (X14 23) at (8,4){$X_{14|23}$};
% \node[invisivertex] (X134) at (9,4){$X_{134}$};
% \node[invisivertex] (X234) at (11,4){$X_{234}$};

% \node[invisivertex] (1) at (5,6){$\hat{1}=\text{span}(\vec{1})$};

% Order Relations, Grouped by Rank

% \path[-] (H12) edge [bend left =0] node[above] {} (0);
\path[-] (H23) edge [bend left =0] node[above] {} (0);
\path[-] (H13) edge [bend left =0] node[above] {} (0);
\path[-] (H24) edge [bend left =0] node[above] {} (0);
\path[-] (H14) edge [bend left =0] node[above] {} (0);
% \path[-] (H34) edge [bend left =0] node[above] {} (0);

% \path[-] (H12) edge [bend left =0] node[above] {} (X123);
% \path[-] (H12) edge [bend left =0] node[above] {} (X124);
% \path[-] (H12) edge [bend left =0] node[above] {} (X12 34);
% \path[-] (H23) edge [bend left =0] node[above] {} (X123);
% \path[-] (H23) edge [bend left =0] node[above] {} (X14 23);
% \path[-] (H23) edge [bend left =0] node[above] {} (X234);
% \path[-] (H13) edge [bend left =0] node[above] {} (X123);
\path[-] (H13) edge [bend left =0] node[above] {} (X13 24);
% \path[-] (H13) edge [bend left =0] node[above] {} (X134);
% \path[-] (H24) edge [bend left =0] node[above] {} (X124);
\path[-] (H24) edge [bend left =0] node[above] {} (X13 24);
% \path[-] (H24) edge [bend left =0] node[above] {} (X234);
% \path[-] (H14) edge [bend left =0] node[above] {} (X124);
% \path[-] (H14) edge [bend left =0] node[above] {} (X14 23);
% \path[-] (H14) edge [bend left =0] node[above] {} (X134);
% \path[-] (H34) edge [bend left =0] node[above] {} (X134);
% \path[-] (H34) edge [bend left =0] node[above] {} (X234);
% \path[-] (H34) edge [bend left =0] node[above] {} (X12 34);

% \path[-] (X123) edge [bend left =0] node[above] {} (1);
% \path[-] (X124) edge [bend left =0] node[above] {} (1);
% \path[-] (X12 34) edge [bend left =0] node[above] {} (1);
% \path[-] (X13 24) edge [bend left =0] node[above] {} (1);
% \path[-] (X14 23) edge [bend left =0] node[above] {} (1);
% \path[-] (X134) edge [bend left =0] node[above] {} (1);
% \path[-] (X234) edge [bend left =0] node[above] {} (1);

%\fill[magenta!40,nearly transparent] (0,0) -- (3,3) -- (3,-3) -- cycle;
\end{tikzpicture}
\end{center}
% The chambers are labelled with linear orders on $[n]$ in which $a<b<c<d$ is
% denoted by $abcd$. Note that the chambers of $\K_P$ are precisely the
% linear extensions of $P$.

It happens that here $|\mu(V,X_\pi)|=1$ for all $\pi$ in $\PTPar(P)$, so that $\Poin(P,t)=1+4t+t^2$.
\end{ex}

\begin{ex}\label{cones:ex:colored-blocks}
  Let $P$ be the following poset:

  
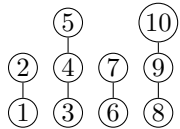
\begin{center}
  \begin{tikzpicture}[scale=.6]
    % \node[invisivertex] (P) at (-1.5,1){$\poset{a}=$};\begin{ex}\label{cones:ex:transverse}
    \node[circledvertex] (1) at (0,0){$1$};
    \node[circledvertex] (2) at (0,1){$2$};
    \node[circledvertex] (3) at (1,0){$3$};
    \node[circledvertex] (4) at (1,1){$4$};
    \node[circledvertex] (5) at (1,2){$5$};
    \node[circledvertex] (6) at (2,0){$6$};
    \node[circledvertex] (7) at (2,1){$7$};
    \node[circledvertex] (8) at (3,0){$8$};
    \node[circledvertex] (9) at (3,1){$9$};
    \node[circledvertex] (10) at (3,2){$10$};

    \path[-] (1) edge [bend left =0] node[above] {} (2);
    \path[-] (3) edge [bend left =0] node[above] {} (4);
    \path[-] (4) edge [bend left =0] node[above] {} (5);
    \path[-] (6) edge [bend left =0] node[above] {} (7);
    \path[-] (8) edge [bend left =0] node[above] {} (9);
    \path[-] (9) edge [bend left =0] node[above] {} (10);

    % Shading
    % \fill[orange!50,nearly transparent] (-.25,1.25) -- (1.25,1.25) -- (1.25,.75)  -- (-.25,.75)-- cycle;
    % \fill[green!50,nearly transparent] (-.25,1.75)--(-.25,2.25) -- (2.25,2.25) -- (2.25,.75)  -- (1.75,.75) -- (1.75,1.75)  -- cycle;
    % \fill[blue!50,nearly transparent] (3.25,-.5) -- (-.25,-.5)--(-.25,.25) -- (1.25,.25) -- (1.25,-.25)  -- (2.75,-.25) -- (2.75,.25) -- (3.25,.25) -- cycle;
    % \fill[purple!50,nearly transparent] (2,-.35) -- (1.5,.15) -- (3,1.5)  -- (3.5,.75)-- cycle;
  \end{tikzpicture}
\end{center}
\noindent
Then $\pi = \{ \{1,4,7\},\{2,5,10\}, \{3,6,8\},\{9\}\}$ in $\Pi_{10}$ is $P$-transverse, represented here by shading the blocks:

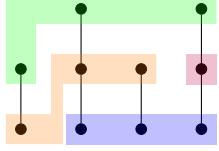
\begin{center}
  \begin{tikzpicture}[scale=.8]
    % \node[invisivertex] (P) at (-1.5,1){$\poset{a}=$};
    \node[point] (1) at (0,0){};
    \node[point] (2) at (0,1){};
    \node[point] (3) at (1,0){};
    \node[point] (4) at (1,1){};
    \node[point] (5) at (1,2){};
    \node[point] (6) at (2,0){};
    \node[point] (7) at (2,1){};
    \node[point] (8) at (3,0){};
    \node[point] (9) at (3,1){};
    \node[point] (10) at (3,2){};

    \path[-] (1) edge [bend left =0] node[above] {} (2);
    \path[-] (3) edge [bend left =0] node[above] {} (4);
    \path[-] (4) edge [bend left =0] node[above] {} (5);
    \path[-] (6) edge [bend left =0] node[above] {} (7);
    \path[-] (8) edge [bend left =0] node[above] {} (9);
    \path[-] (9) edge [bend left =0] node[above] {} (10);

    % Shading
    \fill[blue,nearly transparent] (.75,-.25) -- (3.25,-.25) -- (3.25,.25) -- (.75,.25) -- cycle;
    \fill[orange,nearly transparent] (.7,.75) -- (2.25,.75) -- (2.25,1.25)  -- (.5,1.25) -- (.5,.25)-- (-.25,.25) -- (-.25,-.25) -- (.7, .-.25)  -- cycle;
    \fill[green,nearly transparent] (-.25,.75)--(-.25,2.25) -- (3.25,2.25) -- (3.25,1.75)  -- (.25,1.75) -- (.25,.75)  -- cycle;
    \fill[purple,nearly transparent] (2.75,.75) -- (2.75,1.25) -- (3.25,1.25)  -- (3.25,.75)-- cycle;
  \end{tikzpicture}
\end{center}

Viewed in this way, Proposition \ref{cones:prop:transverse}(iii), roughly speaking, states that $\pi$ is $P$-transverse if and only if its blocks are antichains that can be ``stacked without crossings'' with respect to the Hasse diagram for $P$.
\end{ex}

\begin{proof}[Proof of Proposition~\ref{cones:prop:transverse}]
We will show a cycle of implications: (i) $\Rightarrow$ (ii) $\Rightarrow$ (iii) $\Rightarrow$ (i).
\vskip.1in
\noindent
{\sf (i) implies (ii):}\\
Assume (i), so that there exists some $\mathbf{x}$ in $\R^n$ lying
in the nonempty set
$$
\begin{aligned}
X_\pi \cap \K_P
=X_\pi \cap {\sf int}(C_P)
=\relint{X_\pi \cap C_P}
&=\relint{C_{\overline{P \cup \pi}}}\\
&=
\left\{\mathbf{x}\in\R^n :
      \begin{matrix}
       x_i<x_j\text{ if }(i,j)\in  \overline{P \cup \pi} \text{ but }(j,i)\not\in \overline{P \cup \pi},\\
        x_i=x_j \text{ if both } (i,j),(j,i)\in \overline{P \cup \pi}
    \end{matrix}
      \right\},
\end{aligned}
$$
where the first equality comes from the definition of $\K_P$ and
$C_P$, the second from the fact that $\K_P,C_P$ are full $n$-dimensional, \edit{the third from Proposition~\ref{cones:prop:cone-preposet dictionary}, and the fourth from
equation \eqref{relative-interior-description} above.}
Now to see that (ii) holds, given any pair $i,j$ with $i <_P j$, then $x_i < x_j$ since
$\mathbf{x} \in \K_P$, but then since $(i,j) \in P \subseteq \overline{P \cup \pi}$, the conditions above imply $(j,i) \not\in \overline{P \cup \pi}$, as desired for (ii).

\vskip.1in
\noindent
{\sf (ii) implies (iii):}\\
Assume (ii) holds.  Then every block $B$ of $\pi$ must be an antichain in $P$, else there exists $i \neq j$ in $B$ with $i<_P j$, and then
$(j,i) \in \pi \subseteq \overline{P \cup \pi}$, contradicting (ii).

Now suppose for the sake of contradiction that $P/\pi$ is not a poset. Since $P/\pi$ is a preposet, it can only fail to be antisymmetric, that is, there are blocks $B \neq B'$ of $\pi$ having
both $(B,B'), (B',B)$ in $P/\pi$.  Since both $P, \pi$ are transitive binary relations, this means there
must exist a (periodic) sequence of elements of the form
$$
\cdots \equiv_\pi p_1 <_P p_2 \equiv_\pi p_3  <_P p_4  \equiv_\pi \cdots
<_P p_{m-2}
\equiv_\pi p_{m-1} <_P p_m \equiv_\pi p_1 <_P p_2 \equiv \cdots
$$
alternating relations $(p_i,p_{i+1})$ lying in $P$ and in $\pi$.
Then $p_1 <_P p_2$ and $(p_2,p_1) \in \overline{P \cup \pi}$,
contradicting (ii).

\vskip.1in
\noindent
{\sf (iii) implies (i):}\\
Assume (iii), that is, the blocks of $\pi$ are antichains of $P$, and $P/\pi$ is a poset.
One can then reindex the blocks of $\pi$ such that
$(B_1,B_2,\ldots, B_k)$ is a linear extension of $P/\pi$.
Use this indexing to define a point $\mathbf{x} \in \R^n$ whose $p^{th}$ coordinate $x_p=i$ if $p$ lies in block $B_i$ of $\pi$.

We claim $\mathbf{x}$ lies in $X_\pi \cap \K_P$, verifying (i).
By construction $\mathbf{x}$ lies in $X_\pi$, since its coordinates are constant within
the blocks of $\pi$.  To verify $\mathbf{x} \in \K_P$, given $p <_P q$, one
must check that $x_p < x_q$.  Assume that $p, q$ lie in blocks $B_i, B_j$ of $\pi$,
so that $x_p=i$ and $x_q=j$.  Since the blocks of $\pi$ are antichains in $P$ and $p <_P q$,
one has $i \neq j$, and since $(B_1,B_2,\ldots, B_k)$ is a linear extension of $P/\pi$,
one must have $i < j$, that is, $x_p < x_q$.
\end{proof}

It will help later in identifying
$P$-transverse partitions
to also have the following {\it recursive} characterization.

\begin{prop}
  \label{P-transverse-partitions-recursively}
  Let $P$ be a poset on $[n]$, and $\pi$ a set partition of $[n]$.
  Then $\pi$ lies in $\PTPar(P)$ if and only if it contains a block $B$ with
  these two properties:
  \begin{itemize}
  \item[(a)] $B$ is a subset of the minimal elements of $P$, and
  \item[(b)] if $\hat{\pi}:=\pi\backslash\{B\}$ and $\hat{P}$ is the poset on $[n]\backslash B$ obtained from $P$ by
    removing the elements in $B$, then $\hat{\pi}\in \PTPar(\hat P)$.
  \end{itemize}
\end{prop}

\begin{proof}
  For the forward implication, assume $\pi \in \PTPar(P)$ and $B\in \pi$ is a minimal block. We use Proposition
  \ref{cones:prop:transverse}(iii) to show that (a) and (b) hold. For (a), assume there exists some $x\in B$ which is
  {\it not} minimal in $P$, i.e., there is some $y\in P$ with $x>_P y$. The block $B$ is an antichain, so necessarily $y \not\in B$ and so $y$ lies in a block $B' \neq B$ of $\pi$.
  But then $B >_{P/\pi} B'$ contradicts the minimality of $B$ in $P/\pi$.

  For (b), note that the blocks of $\hat{\pi}$ are a subset of the blocks of
  $\pi$, so they are antichains in $P$ and necessarily also antichains in $\hat{P}$.  Furthermore,
  the preposet $\hat{P}/\hat{\pi}$ on the blocks of $\hat{\pi}$ must be
  a poset, else if \edit{$B', B''$} were two blocks of $\hat{\pi}$ having
  $B' \geq_{\hat{P}/\hat{\pi}} B''$ and $B'' \geq_{\hat{P}/\hat{\pi}} B'$,
  then these same two blocks $B',B''$ in $\pi$ and would have
  have $B' \geq_{P/\pi} B''$ and $B'' \geq_{P/\pi} B'$, a contradiction.

  For the backward implication, assume $\pi$ in $\Pi_n$ has
  a block $B$ satisfying properties (a), (b).
  We use Proposition~\ref{cones:prop:transverse}(i)
  to show that $\pi \in  \PTPar(P)$.
  Since $\hat{\pi} \in \PTPar(\hat{P})$, there is a point
  $\hat{\mathbf{x}} \in \R^{[n] \setminus B}$ in the (nonempty) set
  $X_{\hat{\pi}} \cap \K_{\hat{P}}$, i.e., the coordinates of $\mathbf{x}$
  are constant within each block of $\hat{\pi}$, and
  $\hat{\mathbf{x}}_p < \hat{\mathbf{x}}_q$ whenever $p <_{\hat{P}} q$. Let
  $
  \mathfrak{m}:=\min_\R\{\hat{\mathbf{x}}_p: p \in [n] \setminus B\}
  $
  be the smallest coordinate of $\hat{\mathbf{x}}$, and then
  extend $\hat{\mathbf{x}} \in \R^{[n] \setminus B}$ to a point
  $\mathbf{x} \in \R^n$ by assigning all of the new coordinates
  $x_p$ for $p \in B$ to have the same value, but strictly smaller
  than $ \mathfrak{m}$, e.g., $x_p := \mathfrak{m} - 1$ for all $p \in B$.  One then checks that
  this $\mathbf{x}$ lies in $X_\pi \cap \K_P$:  it lies in $\K_p$ due to the fact that $B$ is a subset of the minimal elements of $P$, and it lies in $X_\pi$ because
  it is constant on the new block $B$ of $\pi$ not already
  in $\hat{\pi}$, as well as constant on the blocks of $\hat{\pi}$.
\end{proof}

\subsection{More examples of $\PTPar(P)$}

\begin{ex}
Given a poset $P$ on $[n]$, its \emph{dual} or \emph{opposite} poset $P^{\opp}$ has the same underlying set $[n]$, but with opposite order relation:  $i \leq_P j$ if and only if $j \leq_{P^\opp} i$.  One can readily check that conditions (ii) and (iii) in Proposition~\ref{cones:prop:transverse} are self-dual in the sense that
$\pi$ in $\Pi_n$ is $P$-transverse if and only if it is $P^{\opp}$-transverse.  Consequently,
$\PTPar(P^{\opp})=\PTPar(P)$, and hence
$
\Poin(P^{\opp},t)=\Poin(P,t).
$
\end{ex}

\begin{ex}
  Given posets $P_1, P_2$, respectively, their \emph{ordinal sum} $P_1 \oplus P_2$ is the poset whose underlying set is the disjoint union $P_1 \sqcup P_2$, and having order relations $x \leq_{P_1 \oplus P_2} y$ if either
  \begin{itemize} \item $x,y\in P_i$ and $x_{\leq_{P_i}} y$ for some $i=1,2$, or
  \item $x \in P_1$ and $y \in P_2$.
  \end{itemize}
  If the underlying sets for $P_1, P_2$ are $[n_1], [n_2]$, one can
  readily check from either of Proposition~\ref{cones:prop:transverse} (ii) or (iii) that a partition $\pi$ of $[n_1] \sqcup [n_2]$ is $P_1 \oplus P_2$-transverse if and only if it is of form $\pi=\{A_1,\ldots,A_k,B_1,\ldots,B_\ell\}$ where $\pi_1=\{A_i\}_{i=1}^k$ and $\pi_2=\{B_j\}_{j=1}^\ell$ are $P_1$-transverse and $P_2$-transverse partitions of $[n_1]$ and $[n_2]$, respectively.
  Bearing in mind that $V=\R^{n_1+n_2}=V_1\oplus V_2$
  where $V_i=\R^{n_i}$ for $i=1,2$, this gives isomorphisms
  $$
  \begin{aligned}{}
  \PTPar(P\oplus P_2) & \cong  \PTPar(P_1) \times  \PTPar(P_2)\\
  [V,X_\pi] &\cong [V_1,X_{\pi_1}] \times  [V_2,X_{\pi_2}]\\
%\mu(V,X_\pi) &= \mu(V_1,X_{\pi_1}) \cdot \mu(V_2,X_{\pi_2})
  \end{aligned}
  $$
  and therefore also
$$
\Poin(P_1 \oplus P_2,t)=\Poin(P_1,t) \cdot \Poin(P_2,t).
$$
\end{ex}

%%%%%%%%%%%%%%%%%%%%%%%
\subsection{Linear extensions, $P$-transverse permutations,
  and proof of Theorem~\ref{tranverse-perm-cycles-thm}}
We recall here the bijection between the chambers of braid arrangement $A_{n-1}$
inside a cone $\K_P$ and the linear extensions of $P$. We then define $P$-transverse permutations
and use them to combinatorially re-interpret $\Poin(P,t)$.

\begin{defn}
  Given two posets $P,Q$ on $[n]$, say that $Q$ \emph{extends} $P$ if $i \leq_P j$ implies $i \leq_Q j$,
  that is, $P \subseteq Q$ as binary relations on $[n]$, or equivalently, their cones satisfy $\K_Q \subseteq \K_P$.
  When $Q$ is a total or linear order $\sigma_1 < \cdots <\sigma_n$ on $[n]$, we identify it with a
  permutation $\sigma = \sigma_1\dots \sigma_n$, and call $\sigma$ a \emph{linear extension} of $P$.  Let $\linext{P}$ denote the set of all linear extensions of $P$.
\end{defn}

Example \ref{A_n basics} noted that chambers of the braid arrangement $A_{n-1}$ are
of the form $\K_\sigma$ for permutations $\sigma$.  Then $\K_\sigma$ is a
chamber lying in the cone $\cham{\K_P}$ if and only if $\sigma$ lies in $\linext{P}$,
giving a bijection
$$
\begin{array}{rcl}
\linext{P} & \longrightarrow &\cham{\K_P}\\
\sigma &\longmapsto & \K_\sigma.
\end{array}
$$
See also \cite[Example 1.3]{stanley}.  Consequently, as noted in \eqref{linear-extensions-as-sums-of-Whitneys}
one has
$$
\#\linext{P} = \sum_{k \geq 0} \whit{k}{P} = \left[ \Poin(P,t) \right]_{t=1}.
$$

\begin{ex} \label{cones:ex:1<2,3<4}
  The poset $P$ defined by $1<_P 2$ and $3 <_P 4$ from Example~\ref{1<2 and 3<4}
  has six linear extensions, shown here labeling the chambers in $\cham{\K_P}$:

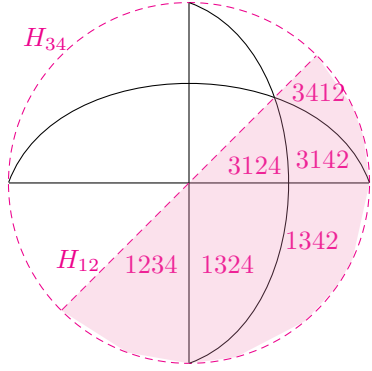
\begin{center}
  \begin{tikzpicture}[scale=.6]
  \draw[boundaryEdge] (0,0) circle (4cm);
  %\draw[dashed] (zero) circle (4cm) {};

  %Boundary Vertices
  \node[origVertex]  (14A) at (180:4) {};
  % \node[origVertex]  (34A) at (65:4) {};
  \node[origVertex] (13A) at (180:4) {};
  \node[origVertex]  (23A) at (90:4) {};
  \node[origVertex]  (12A) at (45:4) {};
  \node[origVertex]  (24A) at (90:4) {};
  \node[origVertex]  (14B) at (0:4) {};
  % \node[origVertex]  (34B) at (240:4) {};
  \node[origVertex] (13B) at (0:4) {};
  \node[origVertex]  (23B) at (270:4) {};
  \node[origVertex]  (12B) at (225:4) {};
  \node[origVertex]  (24B) at (270:4) {};

  % Chamber Labels
  \node[invisivertex] (H14) at (35:3.5) {\textcolor{magenta}{$3412$}};
  \node[invisivertex] (3124) at (15:1.5) {\textcolor{magenta}{$3124$}};
  \node[invisivertex] (3142) at (10:3) {\textcolor{magenta}{$3142$}};
  \node[invisivertex] (1234) at (245:2) {\textcolor{magenta}{$1234$}};
  \node[invisivertex] (H12) at (215:3) {\textcolor{magenta}{$H_{12}$}};
  \node[invisivertex] (1324) at (295:2) {\textcolor{magenta}{$1324$}};
  \node[invisivertex] (1342) at (335:3) {\textcolor{magenta}{$1342$}};
  \node[invisivertex] (H24) at (135:4.5) {\textcolor{magenta}{$H_{34}$}};

  % % Intersection Markers
  % \node[point] (13|24) at (0:2.2) {};
  % \node[point] (124) at (45:2.65) {};
  % \node[point] (14,23) at (90:2.2) {};
  % \node[point] (123) at (0:0) {};
  % \node[point] (234+) at (90:4) {};
  % \node[point] (234-) at (270:4) {};
  % \node[point] (12|34+) at (45:4) {};
  % \node[point] (12|34-) at (225:4) {};
  % \node[point] (134+) at (180:4) {};
  % \node[point] (134-) at (0:4) {};
  %
  % % Intersection Labels
  % \node[invisivertex]  (X 13|24) at (5:3) {$X_{13|24}$};
  % \node[invisivertex]  (X 124) at (37:3.25) {$X_{124}$};
  % \node[invisivertex]  (X 14|23) at (107:2.75) {$X_{14|23}$};
  % \node[invisivertex]  (X 123) at (150:.5) {$X_{123}$};
  % \node[invisivertex] (X 234+) at (90:4.5) {$X_{234}$};
  % \node[invisivertex] (X 234-) at (270:4.5) {$X_{234}$};
  % \node[invisivertex] (X 12|34+) at (45:4.5) {$X_{12|34}$};
  % \node[invisivertex] (X 12|34-) at (225:4.5) {$X_{12|34}$};
  % \node[invisivertex] (X 134+) at (180:4.5) {$X_{134}$};
  % \node[invisivertex] (X 134-) at (0:4.5) {$X_{134}$};

  % Hyperplanes (edges)
  \path[-] (13B) edge[-] (13A);
  \path[-] (14B) edge[out=100, in=200, bend right=70] (14A);
  \path[-] (23B) edge[] (23A);
  \path[-] (12B) edge[boundaryEdge] (12A);
  \path[-] (24B) edge[out=220, in=340, bend right=70] (24A);

  % Shading
  \fill[magenta!60,nearly transparent] (225:4) -- (45:4) -- (25:4) -- (15:4) -- (0:4) -- (335:4)
  -- (325:4) -- (310:4) -- (285:4) -- (270:4) -- (250:4) -- cycle;
  \end{tikzpicture}
\end{center}

\end{ex}

Recall our goal of finding combinatorial interpretations
for $\Poin(P,t)$.
Comparing the expression for $\Poin(P,t)$ given by Corollary~\ref{cor:Lint} in terms of $P$-transverse partitions $\pi$,
and the interpretation of $\mu(V,X_\pi)$ in terms of permutations given in   \eqref{type-A-mobius-interpretation}
motivates the following definition.

\begin{defn}
Given a poset $P$ on $[n]$,  a \emph{$P$-transverse permutation} is a permutation $\sigma$ in $\SN$
for which the set partition $\pi$ given by the cycles of $\sigma$ is a $P$-transverse partition.
Denote by $\PTP(P)$ the set of all $P$-transverse permutations.
\end{defn}

%Note that just as $\PTPar(P)\subseteq\Pi_n$ is a collection of special set partitions, $\PTP(P)\subseteq\Sym_n$ is also a
%collection of special permutations. For $\pi\in\PTP(P)$ we also define $P/\pi$ to be the poset on the set of cycles of
%$\pi$ with the inherited relations $C\le C'$ if there are elements $x\in C$ and $x'\in C'$ with $x\le_P x'$.
%By Proposition~\ref{cones:prop:transverse} we have $\pi\in\PTP(P)$ if and only if the elements in every cycle of $\pi$ form %an antichain in $P$ and $P/\pi$ is a well-defined poset.

Corollary~\ref{cor:Lint} and equation \eqref{type-A-mobius-interpretation}
then immediately imply  this interpretation for $\Poin(P,t)$.

\vskip.1in
\noindent
{\bf Thereom~\ref{tranverse-perm-cycles-thm}.}
{\it
For any poset $P$ on $[n]$, one has
\edit{\[
  \Poin(P,t) = \sum_{\sigma\in\PTP(P)} t^{n - \cyc(\sigma)}.
\]}
In particular, setting $t=1$, one has $\#\linext P=\#\PTP(P)$.
}
\vskip.1in

%%%%%%%%
\subsection{An aside on cones in reflection arrangements}

We digress here to generalize Theorem~\ref{tranverse-perm-cycles-thm}
from posets $P$ parametrizing cones in the type $A_{n-1}$ reflection arrangement,
to {\it any real reflection arrangement}.

  We start first at the level of generality of a {\it complex reflection group} $W$ acting on
  $V=\C^n$.  This means that $W$ is a finite subgroup of $GL(V) \cong GL_n(\C)$ generated by
  {\it (complex, unitary, pseudo-) reflections}, which are elements $w$ in $W$ whose fixed space $V^w=H$ is a
  hyperplane, that is, a codimension one $(\C-)$linear subspace. Let $\A_W$ denote
  the arrangement of all such reflecting hyperplanes, and $\hypinter{\A_W}$ its poset
  of intersections, as before.  This generalizes the type $A_{n-1}$
  setting, where $W=\SN$ and $\hypinter{\A_W} \cong \Pi_n$ is the poset of set partitions of $[n]$.   There is also a well-known generalization of the map  $\SN \rightarrow \Pi_n$
  that sends a permutation $\sigma$ to the set partition $\pi=\{B_1,B_2,\ldots,B_k\}$ whose blocks $B_i$ are the cycles of $\sigma$, given by
  \begin{equation}
    \label{O-S-map}
  \begin{array}{rcl}
    W &\longrightarrow& \hypinter{\A_W}\\
    w &\longmapsto& V^w.
  \end{array}
  \end{equation}
  Orlik and Solomon proved \cite[\S4]{orlik-solomon} the following facts about this map.

  \begin{prop}
    \label{orlik-solomon-lemmas}
    For any finite complex reflection group $W$, the map defined in \eqref{O-S-map}
    has these properties:
    \begin{itemize}
    \item[(a)] The map is well-defined: $V^w$ is an intersection of reflecting hyperplanes, so it
      lies in $\hypinter{\A_W}$.
    \item[(b)] The map {\bf surjects} $W \twoheadrightarrow \hypinter{\A_W}$.
    \item[(c)] The M\"obius function values for lower intervals $[V,X]$ in $\hypinter{\A_W}$ can be expressed via this map as
              $$\mu(V,X)  =  \sum_{w \in W: V^w=X} \det(w).$$
    \end{itemize}
  \end{prop}
  \begin{proof}
    Parts (a) and (b) are rephrasings of \cite[Lemma (4.4)]{orlik-solomon} and part (c) is \cite[Lemma (4.7)]{orlik-solomon}.
    However, we include here a shorter proof, due to C. Athanasiadis\footnote{A version of this proof for {\it real} reflection groups appears (in Greek)  within the proof of Theorem 5.1 on pages 33-34 in the Masters Thesis of Athanasiadis' student C. Savvidou.}.
    Recall that the values $\mu(V,X)$ are uniquely determined by the identity
    $
    \delta_{X,V}=\sum_{Y:V \supseteq Y \supseteq X} \mu(V,Y)
    $
    where $\delta_{X,V}=1$ if $X=V$ and $0$ otherwise.  It therefore suffices to
    check these equalities:
    $$
    \sum_{Y:V \supseteq Y \supseteq X}
      \left( \sum_{w \in W: V^w=Y} \det(w) \right)
    = \sum_{w \in W:  V^w \supseteq X} \det(w) =\sum_{w \in W_X} \det(w)=\delta_{X,V}
    $$
    where here $W_X$ denotes the subgroup of $W$ that fixes $X$ pointwise.
    The last equality follows from Steinberg's Theorem \cite[Thm. 1.5]{steinberg}: he showed $W_X$ is generated by the reflections whose hyperplane contains $X$,
    so that $W_X=\{1\}$ when $X=V$ (implying $\sum_{w \in W_X} \det(w)=\det(1)=1$),
    and otherwise if $X \neq V$,
    summing the (nontrivial) character $\det(-)$ over $W_X$ yields $\sum_{w \in W_X} \det(w)=0$.
   \end{proof}

\noindent
   For {\it real} reflection groups,  part (c) above
   has the following reformulation, generalizing equation \eqref{type-A-mobius-interpretation} above; see \cite[Lemma 5.17]{aguiar},
   \cite[\S2, pp. 413-414]{lehrer-solomon}, and \cite[Prop. 7.2]{RSW}.

  \begin{coro}
    \label{real-O-S-Mobius-cor}
    Let $W$ be a finite real reflection group acting on $V=\R^n$. For any $X$ in $\hypinter{\A_W}$, one has
    $$
    \mu(V,X) = (-1)^{n-\dim(X)} \#\{ w \in W: V^w = X\}.
    $$
  \end{coro}
  \begin{proof}
    Note $W$ acts orthogonally.  We claim that any $w$ acting {\it orthogonally} on $\R^n$  has $\det(w) = (-1)^{n-\dim(V^w)}$; given this,  Lemma~\ref{orlik-solomon-lemmas}(c) would finish the proof.
        To see this claim, note that the eigenvalues $\lambda$ of $w$
        lie on the unit circle in $\C$ , so $\lambda \overline{\lambda}=1$. If the eigenvalue $\lambda$ has multiplicity $m_\lambda$, then $m_{\overline{\lambda}}=m_\lambda$.  Thus
        \begin{align*}
        \det(w)
        &= \prod_{\lambda \in \C} \lambda^{m_\lambda}
        =(+1)^{m_1} (-1)^{m_{-1}}
        \prod_{\text{pairs} \{\lambda \neq \overline{\lambda}\}} (\lambda \overline{\lambda})^{m_\lambda}
        =(-1)^{m_{-1}}.
        \end{align*}
    Modulo two, however, we have
        \begin{align*}
        m_{-1} & = {n-\sum_{\lambda \neq -1} m_\lambda}
        ={n-m_{+1}-\sum_{\lambda \neq \overline{\lambda}} 2m_\lambda}
        ={n-m_{+1}}
        ={n-\dim(V^w)},
      \end{align*}
    and so $\det(w) = (-1)^{n-\dim(V^w)}$.
  \end{proof}

  We specialize now to real reflection groups $W$.
  Here it is known (see, e.g., \cite[Chap. 1]{humphreys}) that $W$ permutes the chambers $\cham{\A_W}$ simply transitively.
  Thus by fixing a choice of base chamber $C_0$, every other chamber $wC_0$ has a
  unique label by some $w$ in $W$, giving a bijection $\cham{\A_W} \leftrightarrow W$.

  Cones $\K$ inside the reflection arrangement $\A_W$ correspond
  to what were called {\it parsets} by the third author \cite[Chap. 3]{reiner},
  or {\it Coxeter cones} by Stembridge \cite{stembridge}, where they were studied as well-behaved
  generalizations of posets $P$ on $[n]$. In particular, the set of chambers
  $\cham{\K}$ inside $\K$ generalizes the set $\linext{P}$ of linear extensions of $P$.
  For a cone $\K$ in $\A_W$, we consider as before the subposet $\inter{\K}$ of
  intersection subspaces interior to $\K$, playing the role of the $P$-transverse
  set partitions
  $\Pi^{\pitchfork}(P)$ in type $A_{n-1}$.  Generalizing the $P$-transverse permutations
  $\PTP(P)$, define the subset
  $$
  W^\pitchfork(\K):=\{ w \in W: V^w \in \inter{\K}\}.
  $$
The real reflection group generalization of Theorem~\ref{tranverse-perm-cycles-thm}  is the following.

\begin{thm}
  \label{W-transverse-elements-thm}
  Any cone $\K$ in the reflection arrangement $\A_W$ for a finite real reflection group $W$ has
\[
  \Poin(\K,t) = \sum_{w \in W^\pitchfork(\K)} t^{n - \dim(V^w)}.
\]
In particular, setting $t=1$, one has $\#\cham{\K} = \# W^\pitchfork(\K)$.
\end{thm}
\begin{proof}
$$
  \begin{aligned}
    \Poin(\K,t)
    &= \sum_{X \in \inter{\K}} | \mu(V,X) | \cdot t^{n-\dim(X)}\\
    &= \sum_{X \in \inter{\K}} \#\{ w \in W: V^w=X \}\cdot t^{n-\dim(X)}
    = \sum_{w \in W^\pitchfork(\K)} t^{n - \dim(V^w)}
 \end{aligned}
$$
 where the second equality used Corollary~\ref{real-O-S-Mobius-cor},
 and the last equality used the definition of $W^\pitchfork(\K)$.
\end{proof}

%%%%%%%%%%%%%%%
\begin{comment}

\begin{rem}
Note that by \eqref{eq:1}, one way to give a combinatorial meaning to $\whit kP$ is by providing a map
$$
\linext{P} \overset{f}{\longrightarrow} \PTPar(P)
$$
such that $\pi=\{B_1,B_2,\ldots\}$ has $\#f^{-1}(\pi) =  \prod_{i}(\#B_i-1)!$,
and then define $\stat(\sigma)=n-\blocks(f(\sigma))$.
There is a motivation for finding such a map $f$ above,
  and more generally, for any cone $\K$ in an arrangement $\A$, to seek a map
  $f:\cham{\K} \longrightarrow \inter{\K}$ having $\#f^{-1}(X)=|\mu(V,X)|$
  for all $X \in \inter{K}$.  Brown \cite[Section 4.2]{brown} considered random walks on $\cham{K}$ that generalize the Bidigare-Hanlon-Rockmore random walks on $\cham{\A}$. He completely
  analyzed the spectrum of their transition matrices in \cite[Theorem 2]{brown}, showing that for
  each $X$ in $\inter{K}$ one has an easily computable eigenvalue $\lambda(X)$ whose  multiplicity is $|\mu(V,X)|$.
\end{rem}

\end{comment}

%%%%%%%%%%%%%%%%%%%%%%%%%%%%%%%%%%%%%%%%%%%%%%%%%%%%%%%%%%%%%%%%%%%%%%%%%%%%%%%%

\section{Bijecting $P$-transverse permutations and linear extensions:
  proof of Theorem~\ref{linear-extensions-LRmax-thm}}
\label{sect:bijection}

%%%%%%%%%%%%%%%%%%%%%%%%%%%%%%%%%%%%%%%%%%%%%%%%%%%%%%%%%%%%%%%%%%%%%%%%%%%%%%%%

The goal in the next few subsections is to define mutually inverse bijections
$$
\begin{array}{rcl}
\PTP(P) & \overset{\Phi}{\longrightarrow} &\linext P\\
\linext P & \overset{\Psi}{\longrightarrow} &\PTP(P)\\
\end{array}
$$
along with the notion of {\it $P$-left-to-right maxima} for $\sigma$ in $\linext P$,
to prove this result from the Introduction.

\vskip.1in
\noindent
{\bf Theorem~\ref{linear-extensions-LRmax-thm}.}
{\it
  For any poset $P$ on $[n]$, one has a bijection $\Phi:\PTP(P)\to \linext P$
  sending the number of cycles to the number of $P$-left-to-right-maxima.
Therefore,
\[
\Poin(P,t) =  \sum_{\sigma \in \linext{P}} t^{n-\LRmax_{P}(\sigma)}
\]
where $\LRmax_P(\sigma)$ denotes the number of $P$-left-to-right maxima of $\sigma$.
}
\vskip.1in

To this end, we first recall the special case the antichain poset $P=\antichain_n$ on $[n]$,
where $\Phi$ is known as the \emph{fundamental bijection} \cite[Proposition 1.3.1]{ec1}.  A permutation
$\tau$ in $\SN$ sending $i \mapsto \tau(i)$ may be written
\begin{itemize}
\item
 in a one-line notation as $\tau=[\tau(1),\tau(2),\ldots,\tau(n)]$, or
 \item
 in a two-line notation
 $\tau=\begin{pmatrix}
    1 & 2 & \cdots& n\\
    \tau(1)& \tau(2) & \cdots & \tau(n) \end{pmatrix}
 $, or
 \item
in various cycle notations that list the $\tau$-orbits on $[n]$, called its {\it cycles},
 in some arbitrarily chosen order, with each cycle listed as $(j,\tau(j), \tau^2(j),\ldots )$ for some arbitrary choice of the first element $j$.
 \end{itemize}
 One way to make the choices non-arbitrary and put the cycle notation in \emph{standard form}
 insists that the first element $j$ listed within each cycle $\tau^{(i)}$
 is the maximum element of the cycle,
 and then insists that the cycles $\tau^{(1)}, \tau^{(2)}, \dots$ are listed with their
 maximum elements in increasing order as integers, that is, $j_1 <_\Z j_2 <_\Z \cdots$.
 The fundamental bijection $\SN \rightarrow \SN$ sends $\tau$ to $\Phi(\tau):=\sigma=[\sigma_1,\sigma_2,\ldots,\sigma_n]$
 by erasing the parentheses around the standard form cycle notation for $\tau$.

 \begin{ex}\label{ex:X}
  The permutation $\tau=[7,5,9,4,2,8,3,6,1]$ in $\Sym_9$ in one-line notation can also be written in
  two-line notation and factored according to its $\tau$-orbits or cycles
  $$
  \tau=\begin{pmatrix}
    1 & 2 &3& 4&5&6&7&8&9\\
    7&5&9&4&2&8&3&6&1
 \end{pmatrix}
 =\begin{pmatrix}
    1  &3&7&9\\
    7&9&3&1
 \end{pmatrix}
 \begin{pmatrix}
    2 &5\\
    5&2
 \end{pmatrix}
 \begin{pmatrix}
    4\\
    4
 \end{pmatrix}
  \begin{pmatrix}
    6&8\\
    8&6
 \end{pmatrix}.
 $$
 Its cycle notation in standard form and image $\sigma=\Phi(\tau)$ are then
$$
\begin{array}{rcccl}
& &\tau &=&(\mathbf{4})(\mathbf{5},2)(\mathbf{8},6)(\mathbf{9},1,7,3)\\
\Phi(\tau)&=&\sigma&=&[4,5,2,8,6,9,1,7,3].
\end{array}
$$
\end{ex}

The inverse map $\Psi$ starts with $\sigma=[\sigma_1,\ldots,\sigma_n]$ in
one-line notation, and determines where to re-insert the parenthesis pairs in the sequence to obtain the
standard form for the cycles of $\tau$.  One only needs to know the locations of the left parentheses,
since then the right parenthesis locations are determined.  There will be one left parenthesis just
to the left of each $\sigma_j$ which is a \emph{left-to-right maximum} (or
\emph{LR-maximum} for short) in $\sigma$, meaning  $\sigma_j = \max_\Z\{\sigma_1,\sigma_2,\dots,\sigma_j\}$.
It is not hard to check that $\Phi, \Psi$ are mutual inverses, and if $\sigma=\Phi(\tau)$, one
has $\cyc(\tau)=\LRmax(\sigma)$, the number of LR-maxima of $\sigma$.

\subsection{The map $\Phi: \PTP(P)\to \linext P$}

To define the map $\Phi: \PTP(P)\to \linext P$ on a $P$-transverse permutation $\tau$, we will first use the $P$-transverse \edit{partition} $\pi$ whose blocks are the cycles of $\tau$ to separate the blocks of $\pi$ and the elements of $P$ into {\it levels}, and then define a notion of \emph{essential} elements. Recall that because $\sigma$ lies in $\PTP(P)$, meaning $\pi$ lies in $\PTPar(P)$,
the quotient preposet $P/\pi$
on the blocks of $\pi$ is actually a poset. This leads to the following definition:

\begin{defn}
Say that a block of $\pi$ which is minimal in $P/\pi$ is of \emph{Level $1$}.
For $k\ge2$, the blocks of  $\pi$ of \emph{Level $k$} are the minimal ones
in the poset obtained from $P/\pi$ by removing all blocks of Level
less than $k$.
\end{defn}

\noindent In other words, a block $B$ of $\pi$ is of Level $k$ if and only if
$$
k=\max\{\ell: \text{ there exists a chain }B=:B_1 >_{P/\pi} B_2 >_{P/\pi}  \cdots >_{P/\pi} B_\ell\},
$$
or even more concretely, $k$ is the maximum among integers $\ell$
with the property that there exist blocks $B=:B_1, B_2,  \dots, B_\ell$ of $\pi$ and elements
$x_i >_P y_{i+1}$ with $x_i \in B_i, y_{i+1} \in B_{i+1}$ for each $i=1,2,\ldots,\ell-1$.

For an element $x$ in $[n]$, define the \emph{Level} of $x$ to
be the Level of the unique block of $\pi$ containing $x$.

\begin{defn}
An element $x$ is {\it essential} if
$x$ has Level $k$ and there exists some $y$ of Level $k-1$ with $x >_P y$; by convention,
all Level $1$ elements $x$ are essential.
\end{defn}

In order to define the map $\Phi: \PTP(P)\to \linext P$, we first introduce a
\emph{standard form} for a $P$-transverse permutation $\tau$.
Let $\tau$ have cycle partition $\pi=\{B_1,\ldots,B_m\}$, lying in $\PTPar(P)$.
List the cycles of $\tau$ in order
$\tau^{(1)},\tau^{(2)},\ldots,\tau^{(m)}$, with the block $B_i$ of $\pi$ corresponding to the cycle $\tau^{(i)}$, and write the cycle $\tau^{(i)}$ as
$\tau^{(i)}=(x_i,\tau(x_i), \tau^2(x_i),\dots)$ for some $x_i$ in $B_i$.
%Note that the blocks of $\pi$ are in bijection with the cycles of $\tau$ and their exists a map sending $c_i$ to $B_i\in \pi$ such that the underlying sets are equal.
Then this listing is the {\it standard form} of $\tau$ if and only if it has the
following properties:
\begin{itemize}
\item If blocks $B_i, B_j$
  have Levels $k, k+1$ in $\pi$, respectively,
  then their indices satisfy $i\leq_\Z j$.
\item For each \edit{$i$}, the first element $x_i$
  listed in the cycle $\tau^{(i)}=(x_i,\tau(x_i), \tau^2(x_i),\dots)$ is the maximal essential element of $B_i$; by Lemma~\ref{lem:Phi}(b), below, every block $B_i$ contains an essential element.
  \item\edit{If $B_i, B_j$ are blocks of Level $k$ with $i<_\Z j$  then $x_i<_\Z x_j$.}
  % The $x_i$ are in increasing order, i.e., $x_1<_\Z x_2<_\Z \cdots <_\Z x_m$.
\end{itemize}
Following the fundamental bijection, the map $\Phi: \PTP(P)\to \linext P$ is defined by putting $\tau\in \PTP(P)$ into standard form and erasing parenthesis.
The following example illustrates this process.

% Given $\tau$ in $\PTP(P)$,
% with cycle partition $\pi$ in $\PTPar(P)$, we will list the cycles of $\tau$ (or the blocks of $\pi$) in a
% \emph{standard form}, and then remove the parentheses to obtain $\Phi(\tau)=\sigma=[\sigma_1,\ldots,\sigma_n]$.
% The standard form is uniquely defined by these properties:  for all $k\ge1$,
% \begin{itemize}
% \item every Level $k$ cycle occurs earlier than every Level $(k+1)$ cycle,
% \item every cycle $(x,\tau(x),\tau^2(x),\ldots)$ starts with the largest integer $x$ {\bf among its
% essential elements} (the cycle contains at least one essential element by Lemma~\ref{lem:Phi}(b) below), and
% \item the starting elements $x_1,x_2,\ldots$ among the cycles of Level $k$ are in increasing order:
% $x_1 <_\Z x_2 <_\Z < \cdots$.
% \end{itemize}

\begin{ex}\label{ex:Phi}
 Let $P$ be the following poset on $[13]$:
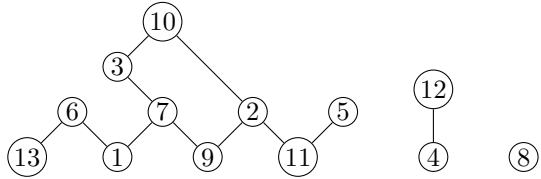
\begin{center}
  \begin{tikzpicture}[scale=.6]
    % \node[invisivertex] (P) at (-1.5,1){$\poset{a}=$};
    \node[circledvertex] (13) at (0,0){$13$};
    \node[circledvertex] (1) at (2,0){$1$};
    \node[circledvertex] (9) at (4,0){$9$};
    \node[circledvertex] (11) at (6,0){$11$};
    \node[circledvertex] (6) at (1,1){$6$};
    \node[circledvertex] (7) at (3,1){$7$};
    \node[circledvertex] (2) at (5,1){$2$};
    \node[circledvertex] (5) at (7,1){$5$};
    \node[circledvertex] (3) at (2,2){$3$};
    \node[circledvertex] (10) at (3,3){$10$};
    \node[circledvertex] (4) at (9,0){$4$};
    \node[circledvertex] (12) at (9,1.5){$12$};
    \node[circledvertex] (8) at (11,0){$8$};

    \path[-] (13) edge [bend left =0] node[above] {} (6);
    \path[-] (6) edge [bend left =0] node[above] {} (1);
    \path[-] (1) edge [bend left =0] node[above] {} (7);
    \path[-] (7) edge [bend left =0] node[above] {} (9);
    \path[-] (9) edge [bend left =0] node[above] {} (2);
    \path[-] (2) edge [bend left =0] node[above] {} (11);
    \path[-] (4) edge [bend left =0] node[above] {} (12);
    \path[-] (11) edge [bend left =0] node[above] {} (5);
    \path[-] (3) edge [bend left =0] node[above] {} (7);
    \path[-] (3) edge [bend left =0] node[above] {} (10);
    \path[-] (2) edge [bend left =0] node[above] {} (10);
  \end{tikzpicture}
  \end{center}
Let
  \[
  \begin{array}{rcll}
   \tau&=& (4) (6,3) (9) (10) (11,7) (12,5,8,2) (13,1) &\in\PTP(P), \text{ so  that}\\
   \pi&=&\{ \{4\}, \{3,6\},\{9\},\{7,10\},  \{2,5,8,12\},\{1,13\} \} &\in \PTPar(P),
   \end{array}
  \]
and  the poset $P/\pi$, drawn as a poset on the cycles of $\tau$, looks as follows:
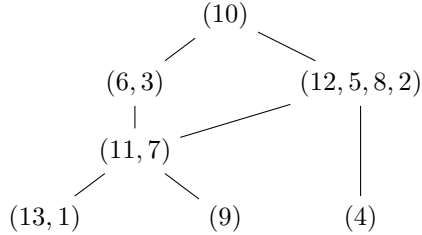
\begin{center}
  \begin{tikzpicture}[scale=.6]
    % \node[invisivertex] (P) at (-1.5,1){$\poset{a}=$};
    \node[] (13) at (0,0){$(13,1)$};
    \node[] (9) at (4,0){$(9)$};
    \node[] (7) at (2,1.5){$(11,7)$};
    \node[] (5) at (7,3){$(12,5,8,2)$};
    \node[] (3) at (2,3){$(6,3)$};
    \node[] (10) at (4,4.5){$(10)$};
    \node[] (4) at (7,0){$(4)$};

    \path[-] (13) edge [bend left =0] node[above] {} (7);
    \path[-] (7) edge [bend left =0] node[above] {} (9);
    \path[-] (7) edge [bend left =0] node[above] {} (5);
    \path[-] (4) edge [bend left =0] node[above] {} (5);
    \path[-] (3) edge [bend left =0] node[above] {} (7);
    \path[-] (3) edge [bend left =0] node[above] {} (10);
    \path[-] (5) edge [bend left =0] node[above] {} (10);
  \end{tikzpicture}
  \end{center}
  One can check that
  \begin{itemize}
  \item
  the Level $1$ cycles are $(13,1),(9),(4)$, with essential elements $1,4,9,13$,
  \item
  there is one Level $2$ cycle $(11,7)$, with one essential element $7$,
  \item
  the Level $3$ cycles are $(6,3),(12,5,8,2)$, with essential elements $2,3,5$,
   \item
  there is one Level $4$ cycle $(10)$, with one essential element  $10$.
  \end{itemize}
Here is $\tau$ in standard form, with essential elements overlined, and Levels separated by bars:
  \[
   \tau=(\overline{4}) (\overline9) (\overline{13},\overline 1) \:\big|\:
    (\overline{7},11) \:\big|\: (\overline 3,6) (\overline 5,8,\overline 2,12) \:\big|\: (\overline{10}),
  \]
  Removing the parentheses (and bars), one obtains its image under $\Phi$:
  \[
    \Phi(\tau) = \sigma =[4,9,13,1,7,11,3,6,5,8,2,12,10] \in \linext P.
  \]
\end{ex}

The next lemma is used to prove the image of $\Phi$ lies in
$\linext P$,
and $\PTP(P) \overset{\Phi}{\rightarrow} \linext P$ is bijective.

\begin{lem}\label{lem:Phi}
  Let $P$ be a  poset on $[n]$ and $\tau \in \PTP(P)$. Then the following properties hold:
  \begin{enumerate}
  \item[(a)] For each $k\ge1$, the Level $k$ elements of $[n]$ form an antichain in $P$.
  \item[(b)] Every cycle contains at least one essential element.
  \item [(c)]The image $\Phi(\tau)$ of $\tau$ is a linear extension of $P$.
  \end{enumerate}
\end{lem}
\begin{proof}
For (a), assume that there were two comparable elements $x <_P y$ with $x,y$ both of Level $k$.  Either $x, y$  lie in
the {\it same} block of $\pi$, contradicting $P$-transversality, or they lie in {\it different} blocks of $\pi$, which would be comparable in $P/\pi$, contradicting both blocks being of Level $k$.

 For (b), note that our most concrete description of a block $B$ having Level $k$ shows that there exist
 blocks $B=:B_1, B_2,  \dots, B_k$ of $\pi$ and elements
$x_i >_P y_{i+1}$ with $x_i \in B_i, y_{i+1} \in B_{i+1}$ for each $i=1,2,\ldots,k-1$.  But then this shows that
the block $B_2$ must be of Level $k-1$, and $x_1 >_P y_2$ shows $x_1$ is essential in $B_1(=B)$.

  For (c), one can show via induction on $k$ that the restriction of $\Phi(\tau)$ to the order ideal of elements of $P$ having Level at most $k$ is a linear extension.  In both the base case $k=1$, and in the inductive step, one notes that one can add in the elements of Level $k$ in any order, because they form an antichain by part (a).
 \end{proof}

%%%
\subsection{The inverse map $\Psi$}

To define the inverse map $\Psi=\Phi^{-1}$ on a linear extension $\sigma$,
we proceed similarly to the previous subsection.  We first cut $\sigma$ into consecutive strings, suggestively
called \emph{Levels}, define a notion of \emph{essential} elements of $\sigma$,
and a notion of  \emph{$P$-left-to-right-maximum}.

\begin{defn}
Given $\sigma=[\sigma_1, \ldots, \sigma_n]$ in $\linext P$, we
recursively break $\sigma$ into disjoint contiguous sequences
$[\sigma_i, \sigma_{i+1}, \ldots, \sigma_{i+j}]$, each forming
an antichain of $P$, and each maximal in the sense that
the slightly longer sequence
$[\sigma_i, \sigma_{i+1}, \ldots, \sigma_{i+j}, \sigma_{i+j+1}]$ is {\it not} an antichain of $P$:
  \begin{itemize}
  \item Let $[\sigma_1,\sigma_2, \dots,\sigma_r]$ be the longest initial
    segment of $\sigma$ whose underlying set
    $\{\sigma_1,\sigma_2,\dots, \sigma_r\}$ is
    an antichain of $P$;  call $\{\sigma_1,\sigma_2, \dots,\sigma_r\}$
    the \emph{Level 1} elements of $\sigma$.
  \item The \emph{Level 2} elements of $\sigma$ are
    $\{\sigma_{r+1}, \sigma_{r+2}, \ldots, \sigma_s\}$,
    where $[\sigma_{r+1}, \sigma_{r+2}, \ldots, \sigma_s]$
    is the longest initial segment of $[\sigma_{r+1}, \sigma_{r+2},\ldots,\sigma_n]$   that forms an antichain in $P$.
  \item Similarly, for $k \geq 3$, the Level $k$ elements of $\sigma$ are defined
    as follows: if the union of all elements of Levels $1,2,\ldots,k-1$
    are $\{\sigma_1,\sigma_2,\dots,\sigma_t\}$, then the set of
    Level $k$ elements is the underlying set of the longest
    initial segment in $[\sigma_{t+1},\sigma_{t+2},\ldots,\sigma_n]$
    that forms an antichain in $P$.
 \end{itemize}
%   the elements of  $\sigma$ of \emph{Level $1$} are those
% contained in the longest initial segment $\{\sigma_1,\sigma_2,\ldots,\sigma_\ell\}$ that forms an antichain in $P$;
% in particular, $\{\sigma_1,\sigma_2,\ldots,\sigma_\ell,\sigma_{\ell+1}\}$ does \emph{not} form an antichain in $P$.
% For each $k \geq 2$, proceed recursively, having already identified the elements of Levels $1,2,\ldots,k-1$, as an
% initial segment $\{\sigma_1,\sigma_2,\ldots,\sigma_\ell\}$ of $\sigma$, and define the elements of \emph{Level $k$}
% be those contained in the longest initial segment $\{\sigma_{\ell+1},\sigma_{\ell+2},\ldots,\sigma_m\}$ of what remains that forms an antichain in $P$.
\end{defn}

\begin{defn} As in the previous section, say $x$ in $\sigma$
  is \emph{essential} if $x$
has Level $k$  in $\sigma$  and there exists an element $y$ of Level $k-1$ in $\sigma$ with $x >_P y$;
again by convention, Level $1$ elements of $\sigma$ are all essential.
\end{defn}

\begin{defn}
Say that an element $x$ is a \emph{$P$-left-to-right-maximum} of $\sigma$, or \emph{$P$-LR-maximum} for short, if
$x$ is essential in $\sigma$, and $x$ appears as a LR-maximum in the usual sense among the subsequence of essential
elements of $\sigma$ having the same Level as $x$.  In other words, if $x$ has Level $k$, it is a $P$-LR-maximum if
the subsequence of essential Level $k$ elements is $[\sigma_{i_1},\sigma_{i_2},\ldots,\sigma_{i_r}]$ for some indices $i_1 < i_2 < \cdots < i_r$, and there is some $j$ with $1 \leq j \leq r$ for which
$\sigma_{i_j}=x=\max\{\sigma_{i_1},\sigma_{i_2},\ldots,\sigma_{i_j}\}$.
We denote the number of $P$-LR maxima of $\sigma$ by $\LRmax_{P}(\sigma)$.
\end{defn}

A map  $\Psi: \linext P \to \SN$ can now be defined in much the same
way as for the fundamental bijection:  starting with $\sigma=[\sigma_1,\ldots,\sigma_n]$ in $\linext P$,
one must determine where to re-insert the parenthesis pairs in the sequence to recover the cycles of $\tau$.
In fact, one only needs to know the locations of the left parentheses, since then the right parenthesis locations are determined,
and there will be one left parenthesis just to the left of each $x=\sigma_j$ which is a $P$-LR maximum.

\begin{ex}\label{ex:poset}
  Let $P$ be the poset in Example~\ref{ex:Phi} and let
  \[
    \sigma = [4,9,13,1,7,11,3,6,5,8,2,12,10] \in \linext P,
  \]
which is $\Phi(\tau)$ from Example~\ref{ex:Phi}.
 The Level decomposition of $\sigma$, with essential elements overlined, looks like
  \[
    \overline{4}, \overline9,\overline{13},\overline 1 \:\big|\:
    \overline{7},11 \:\big|\: \overline 3,6,\overline 5,8,\overline 2,12 \:\big|\: \overline{10}.
  \]
Now create cycles by placing left parentheses just before each $P$-LR maximum:
  \[
   (\overline{4}) (\overline9) (\overline{13},\overline 1) \:\big|\:
    (\overline{7},11) \:\big|\: (\overline 3,6) (\overline 5,8,\overline 2,12) \:\big|\: (\overline{10}).
  \]
The resulting cycle structure gives the $P$-transverse permutation $\Psi(\sigma)$:
  \[
   \Psi(\sigma)= (4) (9) (13,1) (7,11) (3,6) (5,8,2,12)(10),
  \]
which is the $P$-transverse permutation \edit{$\tau$} in Example~\ref{ex:Phi}.
\end{ex}

%The following lemma will be used to prove that the map $\Psi: \linext P \to \PTP(P)$ is a bijection.

%\begin{lem}\label{lem:Psi}
%  Let $P$ be a poset on $[n]$ and $\sigma=\sigma_1\cdots\sigma_n\in \linext P$. Then for each $k \geq 1$,
%\begin{enumerate}
%  \item[(a)] the Level $k$ elements of $\sigma$ form an antichain in $P$, and
%  \item[(b)] the leftmost Level $k$ element of $\sigma$ is always essential.
% \item[(c)] The image $\Psi(\sigma)$ of $\sigma$ is a $P$-transverse permutation.
%  \end{enumerate}
%\end{lem}
%\begin{proof}
%  Property (a) follows from the definition of Level $k$ elements of $\sigma$.
%
%  For (b), the case $k=1$ is vacuous, since all Level $1$ elements are essential.  So assume $k \geq 2$,
%  and that the Level $k-1$ elements of $\sigma$ are the antichain $\{\sigma_{\ell+1},\sigma_{\ell+2},\ldots,\sigma_m\}$.
%  This means that the leftmost Level $k$ element is $x:=\sigma_{m+1}$.  Also, by definition, since $\{\sigma_{\ell+1},
% \sigma_{\ell+2},\ldots,\sigma_m,\sigma_{m+1}\}$ is \emph{not} an antichain, there must be some $y$ in
%   $\{\sigma_{\ell+1},\sigma_{\ell+2},\ldots,\sigma_m\}$ with $x >_P y$.  Thus $x$ is essential.
%b\end{proof}

Note that it is not yet clear that the image of $\Psi$ lies in $\PTP(P)$, but it will follow from the proof of Theorem~\ref{linear-extensions-LRmax-thm}. First we need a technical lemma.

\begin{lem}\label{lem:level k coincides}
Let \edit{$\sigma\in \linext{P}$ and $\tau \in \PTP(P)$} such that $\sigma=\Phi(\tau)$. Then
the set of Level $k$ elements of $\tau$ is precisely the set of Level $k$ elements of $\sigma$.
\end{lem}

\begin{proof}
We prove this by induction on $k$. For the base case ($k=1$), note that the Level $1$ elements of $\tau$ will form an initial
segment $\{\sigma_1,\sigma_2,\ldots,\sigma_\ell\}$ of $\sigma=[\sigma_1,\ldots,\sigma_n]$,
by the definition of $\Phi$.  These Level $1$ elements of $\tau$
will also form an antichain of $P$ by Lemma~\ref{lem:Phi}(a).  On the other hand, we claim that the longer initial
segment
$\{\sigma_1,\sigma_2,\ldots,\sigma_\ell,\sigma_{\ell+1}\}$
\emph{cannot} form an antichain in $P$, because $\sigma_{\ell+1}$ is of Level $2$ in $\tau$ by definition of
$\Phi$, and it is also essential in $\tau$ because it is leftmost in its cycle in the standard form for $\tau$,
and all such elements are essential.  Thus $\sigma_{\ell+1} >_P \sigma_i$ for some $i=1,2,\ldots,\ell$,
showing that the longer segment is not an antichain of $P$.  By definition of $\Psi$, this means
that the Level $1$ elements of
$\sigma$ will be those in the shorter segment $\{\sigma_1,\sigma_2,\ldots,\sigma_\ell\}$.

For the inductive step ($k\geq 2$), we perform the same argument as in the base case, but
replace $P, \tau,\sigma$ with their counterparts $\hat{P}, \hat{\tau},\hat{\sigma}$ in which the elements of Level $1,2,\ldots,k-1$ have
been removed.  One must check that $\hat{\tau}$ lies in $\PTP(\hat{P})$, but this is straightforward,
because if $\sigma, \hat{\sigma}$ have cycle partitions $\pi, \hat{\pi}$, then
$\hat{P}/\hat{\pi}$ is obtained from $P/\pi$ by removing its minimal blocks.
\end{proof}

\begin{proof}[Proof of Theorem~\ref{linear-extensions-LRmax-thm}]
Note $\#\PTP(P)=\#\linext P$ by Theorem~\ref{tranverse-perm-cycles-thm},
and $\Phi$ maps $\PTP(P) \to \linext P$ by Lemma~\ref{lem:Phi}(c).
We therefore claim that it suffices to check $\Psi(\Phi(\tau))=\tau$ for all $\tau$ in $\PTP(P)$.
This will imply that $\Phi$ is injective,
hence bijective, with $\Psi$ its inverse bijection, and thus the image of $\Psi$ is $\PTP(P)$.
Note that by construction, if $\Psi(\sigma)=\tau$, then $\LRmax_P(\sigma)=\cyc(\tau)$, so
Theorem~\ref{linear-extensions-LRmax-thm} would follow.

By Lemma~\ref{lem:level k coincides}, we have if $\sigma = \Phi(\tau)$, then the sets of Level $k$ elements of $\sigma$ and $\tau$ coincide. It follows immediately from Lemma~\ref{lem:level k coincides}
that the \emph{essential} elements in $\sigma$ and $\tau$ coincide, since in each case, their
definition uses only the order $P$ and the partition by Levels. Therefore, we can focus our attention on each Level $k$ separately, where the definition of
$\Phi$ and $\Psi$ coincides almost exactly with their definition in the fundamental bijection,
ignoring the non-essential elements carried along in each Level.  It then
follows that $\Psi(\Phi(\tau))=\tau$ via the same argument as for the fundamental bijection.
\end{proof}

%  If $P$ is a disjoint union of chains, one can similarly relate $\PTP(P)$ and $\linext P$ with multiset permutations. In this case the map $\Phi :\PTP(P)\to \linext P$ induces a bijection on multiset permutations sending a collection of cycles to a word of integers. There is a similar but different bijection due to Foata.  Roughly speaking, our map $\Phi$ is a generalization of the fundamental bijection $\pi\mapsto\widehat{\pi}$ on $\Sym_n$ and Foata's map is a generalization of the process on $\Sym_n$ converting the cycle notation of a permutation into its one-line notation. In the next section we will investigate this situation in more details.

\section{Disjoint Unions of Chains and proofs of Theorems~\ref{chains-thm} and \ref{chains-gf-thm}}
\label{multisets}

%%%%%%%%%%%%%%%%%%%%%%%%%%%%%%%%%%%%%%%%%%%%%%%%%%%%%%%%%%%%%%%%%%%%%%%%%%%%%%%%

The goal of this section is to understand $\Poin(P,t)$ for a poset $P$ that is a disjoint union of chains.  Although Theorem~\ref{linear-extensions-LRmax-thm}
applies to {\it any} poset, when $P$
is a disjoint union of chains, there turns out to be another
elegant expression for $\Poin(P,t)$ stemming from
Foata's theory of multiset permutations,
generalizing equation \eqref{stirling-number-gf} for the antichain poset $P$.

In Subsection \ref{chains:subsect:multiset perms} we review Foata's theory of multiset permutations, in particular his
work with the \emph{intercalation product} and \emph{prime cycle decompositions}. Subsection \ref{subsec:
  partial-commutation-monoids} reviews its relation to partial commutation monoids.  Subsection
\ref{chains:subsect:whitney} shows how the results in Section~\ref{sect:bijection} can be rephrased in terms of multiset
permutations when $P$ is a disjoint union of chains.  Theorem~\ref{chains-thm} is also proved in this subsection.
Finally Subsection \ref{chains:subsect:mmt} employs Foata's theory to give a generalization of MacMahon's Master Theorem
which specializes to Theorem~\ref{chains-gf-thm}, a generating function compiling the Poincar\'e polynomials for
disjoint unions of chains.

\subsection{Multiset Permutations}\label{chains:subsect:multiset perms}

This subsection gives background on the theory
of multiset permutations as introduced by Foata in his PhD
thesis \cite[Section 3.2]{foatathesis}, and extended in later publications \cite[Chapters 3-5]{foata}.
It also appears in Knuth \cite[Section 5.1.2]{knuth}.

\begin{defn}
 Recall that a \emph{(weak) composition $\overline{a}=(a_1,\ldots,a_\ell)$ of $n$}
 is a sequence of nonnegative integers having sum $|a|:=\sum_i a_i=n$.  We will regard $\overline{a}$
 as specifying the multiplicities in a \emph{multiset} $\Ma:=\{1^{a_1}, 2^{a_2}, \dots, \ell^{a_\ell}\}$, that is, a set with repetitions
 $$
 \Ma=\{
 \underbrace{1,1,\ldots,1}_{a_1 \text{ times}},
 \underbrace{2,2,\ldots,2}_{a_2 \text{ times}},\ldots,
 \underbrace{\ell,\ell,\ldots,\ell}_{a_\ell \text{ times}}\}.
 $$
  A \emph{multiset permutation} $\sigma=[\sigma_1, \ldots, \sigma_n]$ is a rearrangement of the elements of $\Ma$, which we will often write in a two-line notation that generalizes that of permutations:
  $$
  \sigma=
  \begin{pmatrix}
  1& \cdots&1 &2& \cdots &2&\cdots&\ell& \cdots & \ell \\
  \sigma_1& \cdots &\sigma_{a_1}&
  \sigma_{a_1+1}& \cdots &\sigma_{a_1+a_2}&
  \cdots&
  \sigma_{a_1+\cdots+a_{\ell-1}+1}& \cdots &\sigma_{n}
  \end{pmatrix}.
  $$

  We denote the set of all multiset permutations of $\Ma$ by $\Sym_\Ma$.  For any $\sigma \in \Sym_\Ma$, we call $\Ma$
  the \emph{support} of $\sigma$, and write $\Ma=\supp{\sigma}$.
\end{defn}

\begin{ex}
  \label{bigchain}
The composition $\overline{a}=(2,3,2,3)$
gives the multiplicities of the multiset
$$
\Ma=\{1^2, 2^3, 3^2, 4^3\}=\{1,1,2,2,2,3,3,4,4,4\}.
$$
Then the following multiset permutation $\sigma$ is an element of $\Sym_\Ma$:
  \[
    \sigma=
    \begin{pmatrix}
      1&1&2&2&2&3&3&4&4&4\\
      2&4&4&3&1&2&1&3&4&2
    \end{pmatrix}.
  \]
\end{ex}

Foata \cite[\S 3.2]{foatathesis} defined an associative \emph{intercalation product} operation on multiset permutations $(\sigma,\rho) \mapsto \sigma \intercal \rho$.
Knuth \cite[\S 5.1.2]{knuth} describes it algorithmically: think of $\sigma, \rho$ in two-line notation as sequences of columns $\left(\begin{smallmatrix}i\\j\end{smallmatrix}\right)$, and juxtapose these sequences of columns. Then perform swaps to \emph{sort} the columns according to their top entries, never swapping two with the \emph{same} top entry.
For example,
\begin{align*}
  \begin{pmatrix}
    2 & 3 & 4\\
    4 & 2 & 3
  \end{pmatrix}
  \intercal
  \begin{pmatrix}
    1&1&2&2&3&4&4 \\
    2&4&3&1&1&4&2\end{pmatrix}
                 = & \left(
                      \begin{array}{ccc|ccccccc}
    2 & 3 & 4 & 1 & 1 & 2 & 2 & 3 & 4 & 4\\
    4 & 2 & 3 & 2 & 4 & 3 & 1 & 1 & 4 & 2
                      \end{array}
                                        \right)\\
    = & \begin{pmatrix}
      1&1&2&2&2&3&3&4&4&4\\
      2&4&4&3&1&2&1&3&4&2
    \end{pmatrix}.
  \end{align*}

\begin{defn}
For each $\ell$, the \emph{intercalation monoid} $\Q$ is the submonoid of all multiset permutations $\sigma$ whose
support $M=\{1^{a_1},  2^{a_2}, \dots, \ell^{a_\ell}\}$
involves only the letters in $\{1,2,\ldots,\ell\}$.
The empty permutation $()$ is
the identity element for $\intercal$,
since
$()\intercal \sigma = \sigma = \sigma\intercal()$.
% so we will denote it by $1:=()$ in $\Q$ when we want to emphasize the monoid structure.
\end{defn}

Note that, just as permutations in the
symmetric group $\mathfrak{S}_n$
do not  commute in general, the monoid $\Q$
is not commutative.  For example
\[
\begin{pmatrix}
  1&2\\
  2&1
\end{pmatrix}
\intercal
\begin{pmatrix}
  1&3\\
  3&1
\end{pmatrix}
=
\begin{pmatrix}
  1&1&2&3\\
  2&3&1&1
\end{pmatrix}
\not=
\begin{pmatrix}
  1&1&2&3\\
  3&2&1&1
\end{pmatrix}
=
\begin{pmatrix}
  1&3\\
  3&1
\end{pmatrix}
\intercal
\begin{pmatrix}
  1&2\\
  2&1
\end{pmatrix}.
\]

However, one can check that $\sigma \intercal\rho
= \rho \intercal \sigma$ when $\sigma,\rho$ are \emph{disjoint}, that is,
  $\supp{\sigma} \cap \supp{\rho}=\emptyset$.

\begin{defn}
  Say $\sigma$ in $\Q$ is \emph{prime} if
  the only factorizations $\sigma=\rho\intercal\tau$
  have either $\rho=()$ or $\tau=()$.
\end{defn}

\begin{ex}
  The permutation $\begin{pmatrix} 2 & 4 & 5 & 7\\ 5 & 7 & 4 & 2 \end{pmatrix}$
  is prime. However, $\begin{pmatrix} 1 & 1 & 2 & 3\\ 2 & 3 & 1 & 1 \end{pmatrix}$ is not prime, since
  \[
    \begin{pmatrix} 1 & 1 & 2 & 3\\ 2 & 3 & 1 & 1 \end{pmatrix}
      = \begin{pmatrix} 1 & 2 \\ 2 & 1 \end{pmatrix} \intercal
      \begin{pmatrix} 1 & 3\\ 3 & 1 \end{pmatrix}.
  \]
  On the other hand $\begin{pmatrix} 2 & 4 & 5 & 7\\ 5 & 7 & 2 & 4 \end{pmatrix}$ is not
  prime, even though its support is multiplicity free, since
  \[
    \begin{pmatrix} 2 & 4 & 5 & 7\\ 5 & 7 & 2 & 4 \end{pmatrix}
      = \begin{pmatrix} 2 & 5 \\ 5 & 2 \end{pmatrix} \intercal
      \begin{pmatrix} 4 & 7 \\ 7 & 4 \end{pmatrix}
      = \begin{pmatrix} 4 & 7 \\ 7 & 4 \end{pmatrix} \intercal
      \begin{pmatrix} 2 & 5 \\ 5 & 2 \end{pmatrix}.
  \]
\end{ex}

It is not obvious, but turns out to be true that $\sigma$ is prime if and only if both
\begin{itemize}
    \item $\supp{\sigma}=M$ is multiplicity free, that is, $M$ is a set not a multiset, and
    \item $\sigma$ has only one cycle when considered as an ordinary
      permutation of the set $M$.
\end{itemize}
We therefore call \emph{prime} elements \emph{prime cycles}.
More generally, one has the following.

\begin{thm}[Foata, 1969 \cite{foata,knuth}]\label{unique}
  Let $\sigma$ be a multiset permutation. Then $\sigma$ has a decomposition into
  a product of prime cycles.
  That is, there exist $k\geq 0$ prime cycles $\sigma^{(1)},\dots,\sigma^{(k)}$ such that
  \begin{align}
  \sigma = \sigma^{(1)}\intercal \sigma^{(2)} \intercal \cdots \intercal\sigma^{(k)}.
  \label{foata decomp}
  \end{align}
   Further, this cycle decomposition of $\sigma$ is unique up to successively
  interchanging pairs of adjacent prime cycles with disjoint support. In particular
  $k$ is unique.
\end{thm}

\begin{defn}
Call $\foata{\sigma}:=k$ the number of prime cycles in the decomposition of $\sigma$
from Theorem~\ref{unique}.
\end{defn}

\begin{ex}\label{chains:ex:fcyc}
  The element $\sigma$ from Example \ref{bigchain} has
  $\foata{\sigma}=4$ and two prime cycle decompositions
  \begin{align*}
  \begin{pmatrix}
      1&1&2&2&2&3&3&4&4&4\\
      2&4&4&3&1&2&1&3&4&2
  \end{pmatrix}
  & =
  \begin{pmatrix}
    2 & 3 & 4\\
    4 & 2 & 3
  \end{pmatrix}
  \intercal
  \begin{pmatrix}
    1 & 2 & 3\\
    2 & 3 & 1
  \end{pmatrix}
  \intercal
  \begin{pmatrix}
    4\\
    4
  \end{pmatrix}
  \intercal
  \begin{pmatrix}
    1 & 2 & 4\\
    4 & 1 & 2
  \end{pmatrix}\\
  & =
  \begin{pmatrix}
    2 & 3 & 4\\
    4 & 2 & 3
  \end{pmatrix}
  \intercal
  \begin{pmatrix}
    4\\
    4
  \end{pmatrix}
  \intercal
  \begin{pmatrix}
    1 & 2 & 3\\
    2 & 3 & 1
  \end{pmatrix}
  \intercal
  \begin{pmatrix}
    1 & 2 & 4\\
    4 & 1 & 2
  \end{pmatrix}.
\end{align*}
\end{ex}

We describe here an algorithm to find a prime cycle decomposition of a multiset permutation, which can be deduced from  \cite[\S 5.1.2]{knuth},
and illustrate how it produces the first of the two
decompositions in Example~\ref{chains:ex:fcyc}.
Encode a multiset permutation $\sigma$ in $\Sym_\Ma$
as two pieces of data:
\begin{itemize}
  \item
  a directed graph $D_\sigma$ on vertex set $\{1,2,\ldots,\ell\}$
having one copy of the directed arc $i \to j$ for each occurrence
of the column $\binom{i}{j}$
in its two-line notation, along with
\item
  specification for each vertex $x$ in $\{1,2,\ldots,\ell\}$ the linear ordering
of the arcs $x \to y$ emanating from $x$, indicating the
left-to-right ordering of the corresponding columns in the two-line
notation.
\end{itemize}
The resulting digraphs are those with
the outdegree equal to the indegree equal to $a_i$ for each $i$.
E.g., the $\sigma$ from Example~\ref{chains:ex:fcyc} has this directed graph
$D_\sigma$, with linear orderings indicated on the arcs out of each vertex:
$$
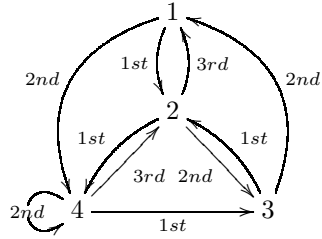
\xymatrix@R=25pt@C=25pt{
  &1\ar@/_/_{1st}[d] \ar@/_2pc/_{2nd}[ddl]               & \\
  & 2\ar_{2nd}[dr] \ar@/_/_{1st}[dl] \ar@/_/_{3rd}[u] & \\
4 \ar@(ul,dl)[]|{2nd} \ar_{1st}[rr] \ar_{3rd}[ur] &  & 3\ar@/_/_{1st}[ul] \ar@/_2pc/_{2nd}[uul]}
$$

With this identification, one factors $\sigma$ recursively.
First produce a prime cycle $\sigma^{(1)}$
for which
\begin{equation}
  \label{recursive-Foata-decomp-step}
\sigma = \sigma^{(1)} \intercal \hat{\sigma}
\end{equation}
via the following algorithm that takes a directed walk in
$D_\sigma$.

\begin{itemize}
  \item Start at the smallest
vertex $i_0$ in $\{1,2,\ldots,\ell\}$ with outdegree $a_{i_0} \geq 1$,
and follow its first outward arc $i_0 \to i_1$.  Then
follow  $i_1$'s first outward arc $i_1 \to i_2$,
follow $i_2$'s first outward arc $i_2 \to i_3$, etc.
\item Repeat until first arriving at a
  previously-visited\footnote{During the process, when one enters a new vertex
  $i_j$ along an arc $i_{j_1} \to i_j$, there will always be
  at least one outward arc $i_j \to i_{j+1}$ leaving $i_j$,
  because each vertex started with its indegree matching its outdegree.} vertex $i_s$,
say $i_s=i_r$ with $r < s$; possibly $r=s-1$.
\item
The directed {\it circuit} $C$ of arcs
$
i_r \to i_{r+1} \to i_{r+2} \to \cdots \to i_{s-1} \to i_s(=i_r)
$
corresponds to a prime cycle $\sigma^{(1)}$ that one can factor
out to the left as in \eqref{recursive-Foata-decomp-step}:
by construction, each of its corresponding columns $\binom{i_t}{j_t}$ occurs
as the {\it leftmost} column of $\sigma$
having $i_t$ as its top element.
\item
Complete the factorization recursively,
replacing $\sigma$ by $\hat{\sigma}$,
removing the arcs $C$ from $D_\sigma$
to give $D_{\hat{\sigma}}$.
\end{itemize}

\begin{ex}
  Here is the algorithm for $\sigma$ above,
  with dotted arrows showing the directed walks in $D_\sigma$:
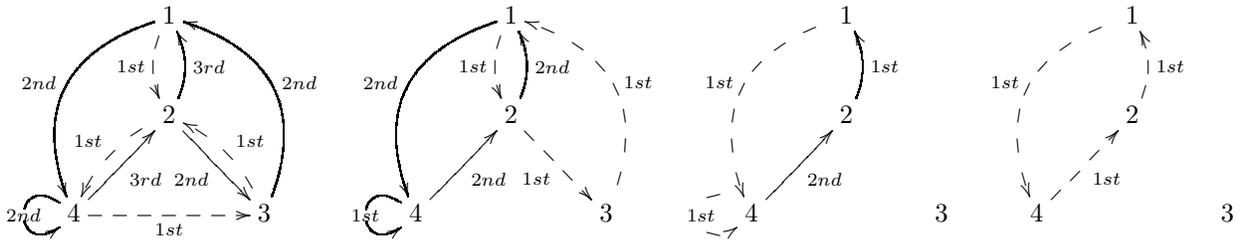
$$
\xymatrix@R=25pt@C=25pt{
  &1\ar@{-->}@/_/_{1st}[d] \ar@/_2pc/_{2nd}[ddl]               & \\
  & 2\ar_{2nd}[dr] \ar@{-->}@/_/_{1st}[dl] \ar@/_/_{3rd}[u] & \\
4 \ar@(ul,dl)[]|{2nd} \ar@{-->}_{1st}[rr] \ar_{3rd}[ur] &  & 3\ar@{-->}@/_/_{1st}[ul] \ar@/_2pc/_{2nd}[uul]}
\quad
\xymatrix@R=25pt@C=25pt{
  &1\ar@{-->}@/_/_{1st}[d] \ar@/_2pc/_{2nd}[ddl]               & \\
  & 2\ar@{-->}_{1st}[dr] \ar@/_/_{2nd}[u] & \\
  4 \ar@(ul,dl)[]|{1st}  \ar_{2nd}[ur] &  & 3 \ar@{-->}@/_2pc/_{1st}[uul]}
\quad
\xymatrix@R=25pt@C=25pt{
  &1 \ar@{-->}@/_2pc/_{1st}[ddl]               & \\
  & 2 \ar@/_/_{1st}[u] & \\
4 \ar@{-->}@(ul,dl)[]|{1st}  \ar_{2nd}[ur] &  & 3}
\quad
\xymatrix@R=25pt@C=25pt{
  &1 \ar@{-->}@/_2pc/_{1st}[ddl]               & \\
  & 2 \ar@{-->}@/_/_{1st}[u] & \\
4  \ar@{-->}_{1st}[ur] &  & 3}
$$
\vskip.2in
$$
\begin{aligned}
\begin{pmatrix}
      1&1&2&2&2&3&3&4&4&4\\
      2&4&4&3&1&2&1&3&4&2
    \end{pmatrix}
&    =\begin{pmatrix}
      2&3&4\\
      4&2&3
    \end{pmatrix}
    \intercal
  \begin{pmatrix}
      1&1&2&2&3&4&4\\
      2&4&3&1&1&4&2
    \end{pmatrix} \\
&=\begin{pmatrix}
      2&3&4\\
      4&2&3
    \end{pmatrix}
    \intercal
  \begin{pmatrix}
      1&2&3\\
      2&3&1
    \end{pmatrix}
        \intercal
  \begin{pmatrix}
      1&2&4&4\\
      4&1&4&2
    \end{pmatrix} \\
&=\begin{pmatrix}
      2&3&4\\
      4&2&3
    \end{pmatrix}
    \intercal
  \begin{pmatrix}
      1&2&3\\
      2&3&1
    \end{pmatrix}
        \intercal
  \begin{pmatrix}
      4\\
      4
    \end{pmatrix}
    \intercal
  \begin{pmatrix}
      1&2&4\\
      4&1&2
    \end{pmatrix}.
\end{aligned}
$$
\end{ex}

\subsection{Partial Commutation Monoids}\label{subsec: partial-commutation-monoids}
It will be helpful to view the intercalation monoid $\Q$ as a partial commutation monoid.
We briefly review some relevant facts about partial commutation monoids.

\begin{defn}
Given a set $\mathbb{A}$, which we call an \emph{alphabet} and a subset of its pairs $C\subseteq\binom{\mathbb{A}}{2}$,
the associated \emph{partial commutation monoid} $\mathcal{M}$ is defined to be the set of equivalence classes on
words $\alpha_1\alpha_2\dots\alpha_k$ in the alphabet $\mathbb{A}$ under the equivalence
relation
\begin{equation}
  \label{eq:equiv_rel}
\alpha_1\alpha_2\dots \alpha_i\alpha_{i+1}\dots\alpha_k
\equiv
\alpha_1\alpha_2\dots \alpha_{i+1}\alpha_{i}\dots\alpha_k
\end{equation}
if $\{\alpha_i,\alpha_{i+1}\}\in C$.
\end{defn}

From this perspective, Foata's Theorem \ref{unique} asserts that $\Q$ is a partial
commutation monoid, whose associated alphabet $\mathbb{A}$ is the set of all prime cycles,
and $C$ the pairs of prime cycles with disjoint supports.

For later use, we point out the following (nontrivial) proposition, see \cite[\S 5.1.2, Exercise 11]{knuth}
and \cite[Exercise 3.123]{ec1}. Given a factorization
of an element
$\alpha=\alpha_1\alpha_2\dots\alpha_k$ in $\mathcal{M}$
a partial commutation monoid, define a poset $\mathcal{P}_\alpha$ on $[k]$ as
the transitive closure of the binary relation containing $(i,j)\in\mathcal{P}_{\alpha}$
when $i<_{\Z}j$ and either $\alpha_i=\alpha_j$ or $\alpha_i\alpha_j\not\equiv \alpha_j\alpha_i$.

\begin{prop}\label{chain:prop:partial-comm}
Given a factorization of
$\alpha=\alpha_1\alpha_2\dots\alpha_k$ in $\mathcal{M}$
a partial commutation monoid,
\begin{enumerate}
\item $\mathcal{P}_\alpha$ does not depend on the choice of factorization of $\alpha$, and
\item there is a bijection between $\linext{\mathcal{P}_\alpha}$ and the factorizations of
$\alpha$ given by
\[
(i_1,\dots,i_k)\mapsto \alpha_{i_1}\dots\alpha_{i_k}.
\]
\end{enumerate}
\end{prop}

\begin{ex}\label{chains:ex:partial-comm}
  The multiset permutation $\sigma$ from Example \ref{chains:ex:fcyc} had
  two prime cycle factorizations
  $$
  \begin{aligned}
  \sigma&=\sigma^{(1)} \intercal \sigma^{(2)} \intercal \sigma^{(3)} \intercal \sigma^{(4)} \\
  &=\sigma^{(1)} \intercal \sigma^{(3)} \intercal \sigma^{(2)} \intercal \sigma^{(4)}
  \end{aligned}
  $$
corresponding to the two linear extensions of the poset $\mathcal{P}_\sigma$ on $[4]$ with
this Hasse diagram:
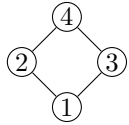
\begin{center}
  \begin{tikzpicture}[scale=.6]
    % \node[invisivertex] (P) at (-1.5,1){$\poset{a}=$};
    \node[circledvertex] (1) at (1,0){$1$};
    \node[circledvertex] (2) at (0,1){$2$};
    \node[circledvertex] (3) at (2,1){$3$};
    \node[circledvertex] (4) at (1,2){$4$};

    \path[-] (1) edge [bend left =0] node[above] {} (2);
    \path[-] (1) edge [bend left =0] node[above] {} (3);
    \path[-] (2) edge [bend left =0] node[above] {} (4);
    \path[-] (3) edge [bend left =0] node[above] {} (4);
  \end{tikzpicture}
  \end{center}
\end{ex}

%%%%%%%
\subsection{Connection with linear extensions and $P$-transverse permutations}
\label{chains:subsect:whitney}
We wish to use Foata's prime cycle decomposition to define
a bijection $\linext{\poset{a}} \rightarrow \PTP(\poset{a})$,
and use this to prove Theorem~\ref{chains-thm}.

We begin with an easy identification of $\linext{\poset{a}}$ with
$\Sym_\Ma$.  For this purpose, given any weak composition
$\overline{a}=(a_1,\ldots,a_\ell)$ of $n$, consider two labelings
of $\poset{a}$, one by the elements $\{1,2,\ldots,n\}$ which
we will call the {\it standardized labeling},
and the second by the elements of the multiset $\Ma$, which we will call the
{\it multiset labeling}.  The standardized labeling labels the
first chain $\chain{a}_1$ by $1,2,\ldots,a_1$ from bottom-to-top,
then the second chain $\chain{a}_2$ by $a_1+1,\ldots,a_1+a_2$
from bottom-to-top, and so on.  The multiset labeling labels
the elements in the first chain $\chain{a}_1$ all by $1$, the second chain  $\chain{a}_2$ all by $2$, etc.

\begin{ex}
  For $n=10$ and $\overline{a}=(2,3,2,3)$, the standardized and multiset labelings of $\poset{a}$ are
  
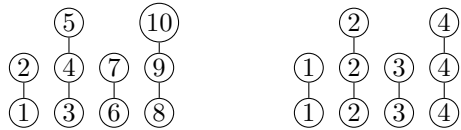
\begin{center}
  \begin{tikzpicture}[scale=.6]
    % \node[invisivertex] (P) at (-1.5,1){$\poset{a}=$};\begin{ex}\label{cones:ex:transverse}
    \node[circledvertex] (1) at (0,0){$1$};
    \node[circledvertex] (2) at (0,1){$2$};
    \node[circledvertex] (3) at (1,0){$3$};
    \node[circledvertex] (4) at (1,1){$4$};
    \node[circledvertex] (5) at (1,2){$5$};
    \node[circledvertex] (6) at (2,0){$6$};
    \node[circledvertex] (7) at (2,1){$7$};
    \node[circledvertex] (8) at (3,0){$8$};
    \node[circledvertex] (9) at (3,1){$9$};
    \node[circledvertex] (10) at (3,2){$10$};

    \path[-] (1) edge [bend left =0] node[above] {} (2);
    \path[-] (3) edge [bend left =0] node[above] {} (4);
    \path[-] (4) edge [bend left =0] node[above] {} (5);
    \path[-] (6) edge [bend left =0] node[above] {} (7);
    \path[-] (8) edge [bend left =0] node[above] {} (9);
    \path[-] (9) edge [bend left =0] node[above] {} (10);
  \end{tikzpicture}
  \qquad \qquad
  \begin{tikzpicture}[scale=.6]
    % \node[invisivertex] (P) at (-1.5,1){$\poset{a}=$};\begin{ex}\label{cones:ex:transverse}
    \node[circledvertex] (1) at (0,0){$1$};
    \node[circledvertex] (2) at (0,1){$1$};
    \node[circledvertex] (3) at (1,0){$2$};
    \node[circledvertex] (4) at (1,1){$2$};
    \node[circledvertex] (5) at (1,2){$2$};
    \node[circledvertex] (6) at (2,0){$3$};
    \node[circledvertex] (7) at (2,1){$3$};
    \node[circledvertex] (8) at (3,0){$4$};
    \node[circledvertex] (9) at (3,1){$4$};
    \node[circledvertex] (10) at (3,2){$4$};

    \path[-] (1) edge [bend left =0] node[above] {} (2);
    \path[-] (3) edge [bend left =0] node[above] {} (4);
    \path[-] (4) edge [bend left =0] node[above] {} (5);
    \path[-] (6) edge [bend left =0] node[above] {} (7);
    \path[-] (8) edge [bend left =0] node[above] {} (9);
    \path[-] (9) edge [bend left =0] node[above] {} (10);
  \end{tikzpicture}
\end{center}
\end{ex}

With this in hand, the following proposition is a straightforward observation.
\begin{prop}
  \label{multiset-perms-are-linear-extensions}
  For any weak composition $\overline{a}$ of $n$,
  one has a bijection
  $$
  \begin{array}{rcl}
    \linext{\poset{a}} &\longrightarrow& \Sym_\Ma\\
      \lambda=[\lambda_1,\ldots,\lambda_n] & \longmapsto & \sigma=[\sigma_1,\ldots,\sigma_n]
  \end{array}
  $$
  replacing $\lambda_i$ by its corresponding multiset label $\sigma_i$, that is, if $\lambda_i$ lies on the $j^{th}$ chain $\chain{a}_j$
  in $\poset{a}$, then $\sigma_i:=j$.
\end{prop}
  \begin{proof}
    The inverse map recovers $\lambda$ from $\sigma$ by
    labeling the $a_j$ occurrences of the value $j$ within $\sigma$
    from left-to-right with the integers in the interval
    $[a_1+a_2+\cdots+a_{j-1}+1, a_1+a_2+\cdots+a_{j-1}+a_j]$.
  \end{proof}

  \begin{ex}
    \label{example-of-multiset-labeling}
    For $\overline{a}=(2,3,2,3)$, this
    bijection maps $\lambda=[3,8,9,6,1,4,2,7,10,5]$ in $\linext{\poset{a}}$ to
    $$
    \sigma=[2,4,4,3,1,2,1,3,4,2]=
    \begin{pmatrix}
      1&1&2&2&2&3&3&4&4&4\\
      2&4&4&3&1&2&1&3&4&2
    \end{pmatrix}.
    $$
  \end{ex}

  We can now define a map $\varphi: \linext{\poset{a}} \rightarrow \SN$,
  which will turn out to be a bijection onto $\PTP(\poset{a})$.

  \begin{defn}
    Fix a weak composition $\overline{a}$ of $n$.  Given $\lambda$ in
    $\linext{\poset{a}}$,
    \begin{itemize}
    \item
      let $\sigma \in \Sym_\Ma$ be its corresponding
      multiset permutation from Proposition~\ref{multiset-perms-are-linear-extensions},
    \item
      label the entries in the top line of $\sigma$'s two-line notation
      with subscripts $1,2,\ldots,n$ from left-to-right,
    \item
      use Foata's Theorem~\ref{unique} to decompose
    $\sigma=\sigma^{(1)} \intercal \cdots \intercal \sigma^{(\ell)}$
    into prime cycles $\sigma^{(i)}$, carrying along the subscripts in
    the top line, and finally
  \item
    replace each prime cycle $\sigma^{(i)}$ with the cyclic
    permutation $\tau^{(i)}$ of the subscripts of its top line.
  \end{itemize}
  Then $\varphi(\lambda):=\tau=\tau^{(1)} \cdots \tau^{(\ell)}$ in $\SN$.
\end{defn}

\begin{ex}
  \label{running-example}
  We continue Example~\ref{example-of-multiset-labeling}.
  Let $\overline{a}=(2,3,2,3)$ and
  $
  \lambda=[3,8,9,6,1,4,2,7,10,5]
  $
  in $\linext{\poset{a}}$. Subscript the top line of its
  corresponding $\sigma$ in $\Sym_\Ma$, and factor as
  in Theorem~\ref{unique}, carrying along subscripts:
$$
    \begin{aligned}
      \sigma=&
    \begin{pmatrix}
      1_1&1_2&2_3&2_4&2_5&3_6&3_7&4_8&4_9&4_{10}\\
      2&4&4&3&1&2&1&3&4&2
    \end{pmatrix}\\
&=\begin{pmatrix}
      2_3&3_6&4_8\\
      4&2&3
    \end{pmatrix}
    \intercal
  \begin{pmatrix}
      1_1&2_4&3_7\\
      2&3&1
    \end{pmatrix}
        \intercal
  \begin{pmatrix}
      4_9\\
      4
    \end{pmatrix}
    \intercal
  \begin{pmatrix}
      1_2&2_5&4_{10}\\
      4&1&2
    \end{pmatrix}\\
&=(2_3, 4_8,3_6) \intercal (1_1, 2_4,3_7) \intercal (4_9) \intercal (1_2, 4_{10}, 2_5)\\
&=\sigma^{(1)} \intercal \sigma^{(2)} \intercal \sigma^{(3)} \intercal \sigma^{(4)}
\end{aligned}
$$
Replacing each prime cycle $\sigma^{(i)}$ with the cycle $\tau^{(i)}$
on its subscripts gives $\varphi(\lambda)=\tau \in \Sym_{10}$:
$$
  \varphi(\lambda)=\tau = \tau^{(1)} \tau^{(2)} \tau^{(3)} \tau^{(4)}
  =(3, 8, 6) (1, 4, 7) (9) (2, 10, 5).
$$
\end{ex}

We can now prove Theorem~\ref{chains-thm}, whose statement we recall here.
\vskip.1in
\noindent
{\bf Theorem~\ref{chains-thm}.}
{\it
  For any composition $\overline{a}$ of $n$, the disjoint union $\poset{a}$ of chains has a bijection
  $$
  \begin{array}{rcl}
    \linext{\poset{a}} &\longrightarrow& \PTP(\poset{a}) \\
                    \sigma &\longmapsto& \tau
  \end{array}
  $$
  with $\cyc(\tau)=\fc$, the number of prime cycles
in Foata's unique decomposition for $\sigma$.  Thus
$$
\Poin(\poset{a},t)=
\sum_{\sigma \in \linext{\poset{a}}} t^{n-\fc}.
$$
}
\vskip.1in
\noindent
\begin{proof}
  We claim that the above map $\varphi: \linext{\poset{a}} \rightarrow \SN$
  is the desired bijection.  Since Theorem~\ref{tranverse-perm-cycles-thm} showed
  $\#\linext{\poset{a}}=\#\PTP(\poset{a})$,
  it suffices to show that the image of $\varphi$ lies in
  $\PTP(\poset{a})$, and that $\varphi$ is injective.

  To see that every $\lambda$ in
  $\linext{\poset{a}}$ has $\varphi(\lambda)=\tau$ lying in
  $\PTP(\poset{a})$, we will induct on the number $\ell$ of
  cycles in $\tau=\tau^{(1)} \cdots \tau^{(\ell)}$,
  which is also the number of prime cycles
  in the decomposition $\sigma=\sigma^{(1)} \intercal \cdots
  \intercal \sigma^{(\ell)}$.  By definition of $\PTP(\poset{a})$,
  we must check that the cycle partition $\pi=\{B_1,\ldots,B_\ell\}$ of $\tau$
  lies in $\PTPar(\poset{a})$, where $B_i$ is the set underlying
  the cycle $\tau^{(i)}$.
  To this end, assume that we produce the factorization
  $\sigma=\sigma^{(1)} \intercal \cdots
  \intercal \sigma^{(\ell)}$ according to the algorithm presented in
  Subsection~\ref{chains:subsect:multiset perms}, and let us check
  that the block $B=B_1$ underlying $\tau^{(1)}$ satisfies the two
  properties (a), (b) in the recursive characterization of $\PTPar(\poset{a})$ from
  Proposition~\ref{P-transverse-partitions-recursively}:
  \begin{itemize}
  \item For (a), the elements of $B$ are all
    minimal in $\poset{a}$ because, in the initial
    factorization step
    $
    \sigma=\sigma^{(1)} \intercal \hat{\sigma},
    $
    each column $\binom{i}{j}$ in the two-line notation of $\sigma^{(1)}$
    is the {\it leftmost} column of $\sigma$ having $i$ as its
    top element, so it corresponds to the {\it bottom} element in the
    $i^{th}$ chain $\chain{a_i}$ of $\poset{a}$.
\item For (b), note that after that initial factorization step,
  the poset $\hat{\poset{a}}$ and partition $\hat{\pi}=\{B_2,\ldots,B_\ell\}$
  will correspond to $\hat{\sigma}$ in the above factorization,
  coming from a $\hat{\lambda}$ in $\linext{\hat{\poset{a}}}$ with $\varphi(\hat{\lambda})=\hat{\sigma}$
  for which the result holds by induction on $\ell$.
\end{itemize}

To show that $\varphi$ is injective, we must give an algorithm
to recover $\lambda$ from $\tau=\varphi(\lambda)$.  It would
be equivalent to recover $\tau$'s
multiset labeled image $\sigma$ in $\Sym_\Ma$
from the bijection in Proposition~\ref{multiset-perms-are-linear-extensions}.
Factoring $\tau$ lets us recover its {\it unordered} set of cycles
$\{ \tau^{(1)}, \ldots,\tau^{(\ell)}\}$, and hence also the
{\it unordered} set of prime cycles $\{ \sigma^{(1)}, \ldots,\sigma^{(\ell)}\}$
that will appear in an intercalation factorization of $\sigma$.  We
would like to know how to properly index $\{ \sigma^{(i)} \}_{i=1,\ldots,\ell}$,
up to interchanging commuting elements, so that we could recover $\sigma$
as their intercalation product
$\sigma=\sigma^{(1)} \intercal \cdots \intercal \sigma^{(\ell)}$.
We claim that this (partial)
ordering is already contained in the information of the unordered set of cycles
$\{ \tau^{(i)} \}_{i=1,\ldots,\ell}$ as follows.  When
two prime cycles $\sigma^{(r)}, \sigma^{(s)}$ do not commute,
it is because they share a common element $i$, so there must exist
two elements $x,y$ in $[n]$ that
come from the $i^{th}$ chain $\chain{a}_i$ in the standardized labeling
of $\poset{a}$, with $x \in \tau^{(r)}, y \in \tau^{(s)}$.
If $x <_\Z y$, then $\sigma^{(r)}$ must occur to the left of $\sigma^{(s)}$
in the intercalation product.
\end{proof}

\subsection{Proof of Theorem~\ref{chains-gf-thm}}
\label{chains:subsect:mmt}

Our goal here is to find a generating function compiling
the Poincar\'e polynomials $\Poin(\poset{a},t)$
for all compositions $\overline{a}$ of length $\ell$.
This uses more of Foata's theory for
the intercalation monoid $\Q$, similar to his deduction of MacMahon's \emph{Master Theorem}.

Since each multiset permutation $\sigma$
has only finitely many intercalation factorizations
$\sigma=\rho \intercal \tau$, one can define
a convolution algebra on the set of
functions $\phi: \Q \rightarrow \Z$
with pointwise addition:
$$
(\phi_1 * \phi_2)(\sigma)
:=\sum_{\rho \intercal \tau=\sigma}
     \phi_1(\rho) \cdot \phi_2(\tau).
$$

Let $\zeta: \Q \rightarrow \Z$ denote the \emph{zeta function} defined by $\zeta(\sigma) =1$ for all $\sigma$ in $\Q$. The zeta function has a unique convolutional inverse $\mu$, called the \emph{M\"obius function}. Foata proved that the M\"obius function can be expressed by the following explicit formula
\[
\mu(\sigma) =
\begin{cases}
(-1)^{\foata{\sigma}} & \text{if }\sigma \text{ is simple},\\
0 & \text{else},
\end{cases}
\]
where $\sigma\in\Q$ is \emph{simple} if all the letters of $\sigma$ are distinct, that is, $\supp{\sigma}$ is
a set, not a multiset.
This may be formulated as an
identity in a completion
$
\Z[[\Q]]:=\left\{ \sum_{\sigma \in \Q} z_\sigma \sigma: z_\sigma \in \Z \right\}
$
of the monoid algebra $\Z[\Q]$,
allowing infinite $\Z$-linear combinations of
elements of $\Q$ (see \cite[Th\'eor\`eme 2.4]{foata}):
%and implicitly \cite[Theorem 3.2]{foatathesis} \galen{these reference Foata's extension of MacMahon's theorem}):

\begin{equation}
    \label{Moebius-identity}
1 = \left( \sum_{\sigma \in \Q}  \sigma \right)
\left( \sum_{\sigma \in \Q}  \mu(\sigma) \thinspace \sigma \right)
= \left( \sum_{\sigma \in \Q}  \sigma \right)
\left( \sum_{\text{simple }\sigma \in \Q}  (-1)^{\foata{\sigma}} \sigma \right).
\end{equation}

Now introduce an $\ell \times \ell$ matrix $B:=(b_{ij})_{i,j=1,2,\ldots,\ell}$ of indeterminates,
and let $\Z[[b_{ij},t]]$ be the (usual, commutative) power series
ring in $\{b_{ij}\}_{i,j=1}^\ell$ along with one further indeterminate $t$.  One can then define a ring homomorphism
$$
\begin{array}{rcl}
\Z[[\Q]]& \overset{u_t}{\longrightarrow} & \Z[[b_{ij},t]] \\
\sigma & \longmapsto & t^{\foata{\sigma}} \cdot \bb{\sigma}
\end{array}
$$
where if
$\sigma=\begin{pmatrix}
i_1 & i_2 & \cdots & i_n\\
\sigma_1 & \sigma_2 & \cdots & \sigma_n
\end{pmatrix}
$
then $\bb{\sigma}:=\prod_{k=1}^n b_{i_k \sigma_k}$.

Applying the homomorphism $u_t$ to
both sides of
\eqref{Moebius-identity} gives a
$t$-version of MacMahon's Master Theorem.

\begin{thm}\label{mmt}
In $\Z[[b_{ij},t]]$ one has the identity
  \[
%  \begin{aligned}
    \sum_{\sigma\in\Q} t^{\fc} \bb{\sigma}
    =
    \left( \sum_{\text{simple }\sigma\in\Q} (-t)^{\fc}\bb{\sigma}
    \right)^{-1}
    =
    \left(\sum_{H\subseteq[\ell]} \sum_{\sigma\in\SH} (-t)^{\fc}\bb{\sigma}
    \right)^{-1}.
%  \end{aligned}
  \]
\end{thm}

\begin{rem}
Setting $t=1$ in Theorem~\ref{mmt} gives an identity in $\Z[[b_{ij}]]$:
\begin{equation}
    \label{Foata-version-of-MMT}
    \sum_{\sigma\in\Q} \bb{\sigma}
    =
    \left(\sum_{H\subseteq[\ell]} \sum_{\sigma\in\SH} (-1)^{\fc}\bb{\sigma}
    \right)^{-1},
  \end{equation}
which is equivalent to an identity in Foata's proof of the (commutative) MacMahon Master Theorem,
as we recall here.  Introduce two
sets of $\ell$ variables
$\mathbf{x}=(x_1,\ldots,x_\ell),
\mathbf{y}=(y_1,\ldots,y_\ell)$
related by the matrix $B$ of indeterminates
as follows: $\mathbf{y}=B\mathbf{x}$,
that is, $y_i = \sum_{j} b_{ij} x_j$.
Then MacMahon's Master Theorem is this identity
in $\Z[[b_{ij}]]$:
\begin{equation}
    \label{standard-MMT}
\sum_{\overline{a} \in \{0,1,2,\ldots\}^\ell }
\left(\text{coefficient of }
\mathbf{x}^{\overline{a}} \text{ in }
\mathbf{y}^{\overline{a}} \right)
=\det(I_\ell - B)^{-1},
\end{equation}
where
$
\mathbf{x}^{\overline{a}}:=x_1^{a_1} \cdots x_\ell^{a_\ell}.
$
It is not hard to check that the
left sides and right sides of \eqref{standard-MMT}
and
\eqref{Foata-version-of-MMT}
are the same:  the left side of \eqref{Foata-version-of-MMT}
needs to be grouped according to the multiplicity vector $\overline{a}$ giving the
support $\supp{\sigma}$, and the right side must be reinterpreted in terms of the
permutation expansion of a determinant.
\end{rem}

\begin{rem}
Theorem~\ref{mmt} is similar in spirit to Garoufalidis-Lê-Zeilberger's \emph{quantum MacMahon Master Theorem}
  \cite[Theorem 1]{garoufalidis} (see also Konvalinka-Pak \cite[Theorem 1.2]{Konvalinka}). Their quantum
  version inserts a $(-q)^{-\text{inv}(\sigma)}$ in order to produce a $q$-determinant, but
  $\text{inv}(\sigma)\not=\foata{\sigma}$.
\end{rem}

We now specialize $b_{ij}=x_j$ in
Theorem~\ref{mmt} to deduce Theorem~\ref{chains-gf-thm}, whose statement we recall here.

\vskip.1in
\noindent{\bf Theorem~\ref{chains-gf-thm}.}
\emph{
For $\ell=1,2,\ldots$, one has
$$
\sum_{\overline{a} \in \{1,2,\ldots\}^\ell}
\Poin(\poset{a},t) \cdot \mathbf{x}^{\overline{a}}
=\frac{1}{
1-\sum_{j=1}^\ell e_j(\mathbf{x}) \cdot (t-1)(2t-1)\cdots ((j-1)t-1)},
$$
where
$
e_j(\mathbf{x}):=\sum_{1 \leq i_1 <\cdots<i_j\leq \ell} x_{i_1} \cdots x_{i_j}
$
is the $j^{th}$ elementary symmetric function.
}

\vskip.1in
\begin{proof}
  Setting $b_{ij} = x_j$ in Theorem \ref{mmt} gives
  \begin{equation}
  \label{specialized-t-MMT}
    \sum_{\sigma\in\Q} t^{\fc} \prod_{k} x_{\sigma_k}
    = \left(\sum_{H\subseteq[\ell]} \sum_{\sigma\in\SH} (-t)^{\fc} \prod_{k \in H}  x_k
    \right)^{-1}.
  \end{equation}
  Let us manipulate both sides of equation \eqref{specialized-t-MMT}.  On the left side, grouping terms according to $\supp{\sigma}$ gives

 $$
 \sum_{\overline{a} \in \{0,1,2,\ldots\}^\ell}
  \mathbf{x}^{\overline{a}} \sum_{\sigma \in \mathfrak{S}_{\multi{a}}} t^{\foata{\sigma}}.
 $$
On the right side of \eqref{specialized-t-MMT}, any subset $H \subseteq [\ell]$ of cardinality $j \geq 1$ satisfies
$$
\sum_{\sigma\in\SH} (-t)^{\fc}=
\sum_{\sigma \in \mathfrak{S}_j} (-t)^{\cyc(\sigma)}
=(-t)(1-t)(2-t) \cdots (j-1-t)
$$
by \eqref{stirling-number-gf}.
Therefore grouping according
to $j=\#H$,
and noting
$
\sum_{\substack{H\subseteq[\ell]:\\ \#H=j}}
\prod_{k \in H} x_k
=e_j(\mathbf{x})
$
lets one rewrite the sum inside the parentheses on the right side of \eqref{specialized-t-MMT} as this:
$$
1+\sum_{j=1}^\ell
 (-t)(1-t)(2-t) \cdots (j-1-t)
\cdot
e_j(\mathbf{x}).
$$
So far this gives
$$
\sum_{\overline{a} \in \{0,1,2,\ldots\}^\ell}
  \mathbf{x}^{\overline{a}} \sum_{\sigma \in \mathfrak{S}_{\multi{a}}} t^{\foata{\sigma}}
  =
  \left(
  1+\sum_{j=1}^\ell
 (-t)(1-t)(2-t) \cdots (j-1-t)
\cdot
e_j(\mathbf{x}) \right)^{-1}.
$$
Now perform two more substitutions:  first replace
$t$ by $t^{-1}$, giving this
\[
\sum_{\overline{a} \in \{0,1,2,\ldots\}^\ell}
  \mathbf{x}^{\overline{a}} \sum_{\sigma \in \mathfrak{S}_{\multi{a}}} t^{-\foata{\sigma}}
   =
  \left(
  1+\sum_{j=1}^\ell
 (-t^{-1})(1-t^{-1})(2-t^{-1}) \cdots (j-1-t^{-1})
\cdot
e_j(\mathbf{x}) \right)^{-1},
\]
and then replace $x_i$ by $tx_i$ for $i=1,2,\ldots,\ell$, so that $\mathbf{x}^{\overline{a}} \mapsto t^{\total{a}} \mathbf{x}^{\overline{a}}$ and $e_j(\mathbf{x}) \mapsto t^j e_j(\mathbf{x})$, giving this
 \[
 \sum_{\overline{a} \in \{0,1,2,\ldots\}^\ell}
  \mathbf{x}^{\overline{a}} \sum_{\sigma \in \mathfrak{S}_{\multi{a}}} t^{|\overline{a}|-\foata{\sigma}}
  =\left(
  1-\sum_{j=1}^\ell
 (t-1)(2t-1) \cdots ((j-1)t-1)
 \cdot
e_j(\mathbf{x}) \right)^{-1}.
 \]
Comparison of the left side with Theorem~\ref{chains-thm} shows that this last equation is Theorem~\ref{chains-gf-thm}.
\end{proof}

\begin{rem}
  We justify here the claim from the Introduction that
  Theorem~\ref{chains-gf-thm} generalizes the formula \eqref{stirling-number-gf}:
  $$
  \Poin(\antichain_\ell,t)=1(1+t)(1+2t) \cdots (1+(\ell-1) t).
  $$
  Since $\antichain_\ell=P_{\overline{a}}$ where $\overline{a}=(1,1,\ldots,1)$, we seek to explain why the coefficient of
  $x_1\dots x_\ell$ in the power series on the right side in Theorem~\ref{chains-gf-thm} should be
  $1(1+t)(1+2t) \cdots (1+(\ell-1) t)$.
Introducing the abbreviation $\langle t \rangle_j:=(t-1)(2t-1)\cdots ((j-1)t-1)$,
the right side in Theorem~\ref{chains-gf-thm} can be rewritten and expanded as
  \begin{equation}
    \label{theorem-RHS}
\frac{1}{
  1-\sum_{j=1}^\ell e_j(\mathbf{x})  \langle t \rangle_j }
=\sum_{n\ge0} \left(\sum_{j=1}^\ell e_j(\mathbf{x}) \langle t \rangle_j \right)^n.
  \end{equation}
If we let $A_\ell$ denote the coefficient of $x_1 \cdots x_\ell$ in this series,
then it suffices to explain why
\begin{equation}
\label{desired-recursion}
A_{\ell+1} = (1+\ell t) \cdot A_\ell.
\end{equation}
For this coefficient extraction, it is safe to
replace each $e_j(\mathbf{x})=e_j(x_1, x_2,\ldots ,x_\ell)$ in  \eqref{theorem-RHS}
with an infinite variable version  $e_j(\mathbf{x})=e_j(x_1, x_2,\ldots)$.
Extracting the coefficient of $x_1 \cdots x_\ell$ on the right side in \eqref{theorem-RHS}
  shows
$$
    A_\ell= \sum_{ \substack{ \text{ordered set partitions} \\ \pi=(B_1,\ldots,B_n) \text{ of }[\ell]}} w(\pi),
\qquad \text{ where }
w(\pi) := \prod_{B_i\in \pi} \langle t \rangle_{|B_i|}.
$$
To explain \eqref{desired-recursion}, note that each ordered set partition $\hat{\pi}$ of $[\ell+1]$ can be obtained
from a unique ordered set
partition $\pi=(B_1,\dots,B_n)$ of $[\ell]$ as follows:   either $\hat{\pi}$ has
added $\ell+1$ into one of the preexisting blocks $B_i$ of $\pi$,
or $\hat{\pi}$ has a singleton block
$\{\ell+1\}$, inserted into one of the $n+1$ locations in the sequence   $(B_1,\ldots,B_n)$.
Thus having fixed an ordered set partition $\pi$ of $[\ell]$, the sum of $w(\hat{\pi})$ over $\hat{\pi}$
which correspond to $\pi$ is this sum:
\[
\begin{aligned}
 &w(\pi) \cdot (|B_1|t-1)+\cdots  +w(\pi) \cdot (|B_n|t-1) + \underbrace{w(\pi) + w(\pi)+ \cdots + w(\pi)}_{n+1\text{ times}} \\
 &=w(\pi) \left( \sum_{i=1}^n (|B_i|t-1) + n+1\right)
 =w(\pi) \left( \ell t - n + n+1\right)
 = w(\pi) \cdot (1+\ell t).
 \end{aligned}
\]
Summing this over all possible $\pi$ gives \eqref{desired-recursion}.
\end{rem}

%%%%%%%%%%%%%%%%%%%%%%%%%%%%%%%%%%%%%%%%%%%%%%%%%%%%%%%%%%%%%%%%%%%%%%%%%%%%%%%%

\section{Posets of width two and proof of Theorem~\ref{width-two-thm}}
\label{sect:width 2}

%%%%%%%%%%%%%%%%%%%%%%%%%%%%%%%%%%%%%%%%%%%%%%%%%%%%%%%%%%%%%%%%%%%%%%%%%%%%%%%%

The \emph{width} of a poset $P$ is the maximum size of an antichain in $P$.
A famous result of Dilworth from 1950 (see \cite[Ch 3, Exer 77(d)]{ec1}), asserts that the
width $d$ of $P$ is the minimum number of chains required in a
\emph{chain decomposition} $P = P_1 \cup P_2 \cup \cdots \cup P_d$,
that is, where each $P_i$ is a totally ordered subset $P_i \subseteq P$.

Consequently, a poset $P$ of width two can be decomposed into two
chains $P=P_1 \cup P_2$, possibly with some order relations between
elements of $P_1$ and $P_2$.

Recall that \edit{in} the Introduction we defined a descent-like statistic on
$\sigma=[\sigma_1,\ldots,\sigma_n]$
in $\linext{P}$,  as the cardinality
$\des_{P_1,P_2}(\sigma):=\#\Des_{P_1,P_2}(\sigma)$ of \edit{this} set,
% \edit{Recall the descent-like statistic $\des$ defined in the Introduction. For $\sigma=[\sigma_1,\ldots,\sigma_n]$ in $\linext{P}$, define $\des_{P_1,P_2}(\sigma):=\#\Des_{P_1,P_2}(\sigma)$ }
\begin{equation}
  \label{width-two-Des-defn}
\Des_{P_1,P_2}(\sigma)
:=\{ i \in [n-1]: \sigma_i \in P_2, \,\, \sigma_{i+1}\in P_1, \text{ with } \sigma_i, \sigma_{i+1}\text{ incomparable in }P\},
\end{equation}
in order to state the following result.
\vskip.1in
{\bf Theorem~\ref{width-two-thm}.}
{\it
For a width two poset decomposed into two chains as $P=P_1 \cup P_2$, one has
$$
\Poin(P,t) =
\sum_{\sigma \in \linext{P}} t^{\des_{P_1,P_2}(\sigma)}.
$$
}
\vskip.1in

To prove this, we start with the following observation.

\begin{coro}\label{same size}
For posets $P$ of width two, one has
 $$
\Poin(P,t) = \sum_{\pi \in \PTPar(P)} t^{\pairs(\pi)}
$$
where $\pairs(\pi)$ is the number of two-element blocks $B_i$ in
$\pi$.  In particular, setting $t=1$,
 $$
  \#\linext P = \# \PTPar(P) (=\#\PTP(P)).
 $$
  \end{coro}

\begin{proof}
Antichains in $P$ have at most two elements, so Proposition~\ref{cones:prop:transverse}(iii) implies that $P$-transverse permutations have only
$1$-cycles and $2$-cycles.  But then this implies that the map $\PTP(P) \rightarrow \PTPar(P)$
sending a $P$-transverse permutation $\tau$ to the set partition $\pi$ given by its cycles is a bijection,
with \edit{$\pairs(\pi)=n-\cyc(\tau)$.}
The result then follows from Theorem~\ref{tranverse-perm-cycles-thm}.
\end{proof}

\begin{ex}
\label{two-disjoint-chain-example}
  Let $P=\chain{a} \sqcup \chain{b}$ be a poset which is a disjoint union of two
  chains $\chain{a}, \chain{b}$ having $a,b$ elements respectively.  One can check that a $P$-transverse partition having $\pairs(\pi)=k$ is completely determined by the choice of a $k$ element subset $x_1 <_P \cdots <_P x_k$ from $\chain{a}$ and a $k$ element subset
  $y_1 <_P \cdots <_P y_k$ from $\chain{b}$ to constitute the two-element blocks, as follows: $\{x_1,y_1\},\ldots,\{x_k,y_k\}$.  This implies
  $$
  \Poin(\chain{a} \sqcup \chain{b} ,t)=\sum_{k=0}^{\min(a,b)} \binom{a}{k} \binom{b}{k} t^k.
  $$
  This is consistent with
  $\#\linext{\chain{a} \sqcup \chain{b} } =\binom{a+b}{a}$, since
  setting $t=1$ in the equation above  gives
  $$
  \binom{a+b}{a}=\sum_{k=0}^{\min(a,b)} \binom{a}{k} \binom{b}{k}
  $$
  which is an instance of the Chu-Vandermonde summation.
\end{ex}

In light of Corollary~\ref{same size}, to prove Theorem~\ref{width-two-thm}, one would need a bijection
from $\linext P$ to $\PTPar(P)$ (or $\PTP(P)$) that sends the statistic $\des_{(P_1,P_2)}(-)$ to the number of
pairs or number of $2$-cycles.  Unfortunately, there does not seem to be a consistent labeling of a width
two poset $P=P_1 \cup P_2$ to make the bijection $\Psi$ from Section~\ref{sect:bijection} play this role.
Nevertheless, having fixed the chain decomposition\footnote{So we assume here that
$P_1 \cap P_2=\varnothing$, but there may be order relations between elements of $P_1$ and $P_2$.}
$P=P_1 \sqcup P_2$, we provide in the proof below such a bijection
$
\Omega: \linext P \rightarrow \PTPar(P).
$

\begin{proof}[Proof of Theorem~\ref{width-two-thm}.]
  We describe $\Omega$ and $\Omega^{-1}$ recursively, via induction on $n:=\#P$.  There are two cases, based on whether $P$ has one or two minimal elements.

\vskip.1in
\noindent
{\sf Case 1.}
There is a unique minimum element $p_0\in P$.

In this case, given $\sigma=[\sigma_1,\dots,\sigma_n]$ in $\linext{P}$,
we must have $\sigma_1=p_0$, so that $\{p_0\}$ should be a singleton block of $\pi=\Omega(\sigma)$,
and one produces the remaining blocks of $\pi$ by applying $\Omega$ recursively to $[\sigma_2,\ldots,\sigma_n]$.
This is depicted schematically here:

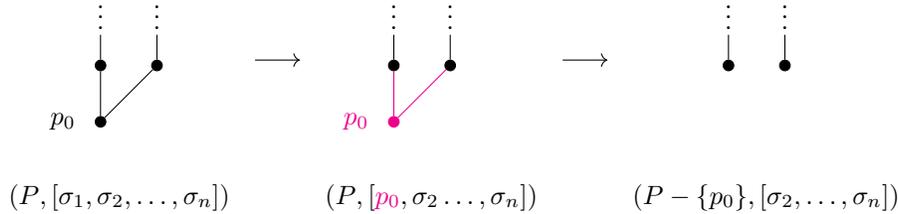
\begin{center}
\begin {tikzpicture}[]

\node[invisivertex] (P) at (-.5,0){$p_0$};
\node[point] (1) at (0,0){};
\node[point] (2) at (0,.75){};
\node[point] (3) at (.75,.75){};

\node[invisivertex] (dotsL) at (0,1.5){$\vdots$};
\node[invisivertex] (dotsR) at (.75,1.5){$\vdots$};
\node[invisivertex] (ell2) at (.25,-1){$(P,[\sigma_1,\sigma_2,\dots,\sigma_n])$};
%\node[invisivertex] (ell3) at (.5,-1.5){$\pi=\{B_1,\dots,B_k\}$};

\path[-] (1) edge [bend left =0] node[above] {} (2);
\path[-] (1) edge [bend left =0] node[above] {} (3);
\path[-] (2) edge [bend left =0] node[above] {} (dotsL);
\path[-] (3) edge [bend left =0] node[above] {} (dotsR);
\end{tikzpicture}
\begin{tikzpicture}[]
\node[invisivertex] (P1) at (0,-1.5){};
% \node[invisivertex] (P2) at (0,1){};
\node[invisivertex] (1) at (0,.5){};
\node[invisivertex] (2) at (.75,.5){};
\path[->] (1) edge [bend left =0] node[above] {} (2);
\end{tikzpicture}
\begin {tikzpicture}[]

\node[invisivertex] (P) at (-.5,0){\textcolor{magenta}{$p_0$}};
\node[magpoint] (1) at (0,0){};
\node[point] (2) at (0,.75){};
\node[point] (3) at (.75,.75){};

\node[invisivertex] (dotsL) at (0,1.5){$\vdots$};
\node[invisivertex] (dotsR) at (.75,1.5){$\vdots$};
% \node[invisivertex] (ell1) at (.5,-.5){$p=\sigma_1$};
\node[invisivertex] (ell2) at (.5,-1){$(P,[\textcolor{magenta}{p_0},\sigma_2\dots,\sigma_n])$};
%\node[invisivertex] (ell3) at (.5,-1.5){$\pi=\{B_1,\dots,B_k\}$};

\path[-] (1) edge [bend left =0,color=magenta] node[above] {} (2);
\path[-] (1) edge [bend left =0,color=magenta] node[above] {} (3);
\path[-] (2) edge [bend left =0] node[above] {} (dotsL);
\path[-] (3) edge [bend left =0] node[above] {} (dotsR);
\end{tikzpicture}
\begin{tikzpicture}[]
\node[invisivertex] (P1) at (0,-1.5){};
% \node[invisivertex] (P2) at (0,1){};
\node[invisivertex] (1) at (0,.5){};
\node[invisivertex] (2) at (.75,.5){};
\path[->] (1) edge [bend left =0] node[above] {} (2);
\end{tikzpicture}
\begin {tikzpicture}[]

\node[invisivertex] (P) at (-.5,0){};
% \node[point] (1) at (0,0){};
\node[point] (2) at (0,.75){};
\node[point] (3) at (.75,.75){};

\node[invisivertex] (dotsL) at (0,1.5){$\vdots$};
\node[invisivertex] (dotsR) at (.75,1.5){$\vdots$};
\node[invisivertex] (ell) at (.5,-1){$(P-\{p_0\},[\sigma_2,\dots,\sigma_n])$};
%\node[invisivertex] (ell3) at (.5,-1.5){$\pi=\{B_1,\dots,B_k,\{p\}\}$};

% \path[-] (1) edge [bend left =0] node[above] {} (2);
% \path[-] (1) edge [bend left =0] node[above] {} (3);
\path[-] (2) edge [bend left =0] node[above] {} (dotsL);
\path[-] (3) edge [bend left =0] node[above] {} (dotsR);
\end{tikzpicture}
\end{center}

For the inverse map $\Omega^{-1}$, given a $P$-transverse partition $\pi$, since the blocks of $\pi$ are antichains
in $P$, the unique minimum element $p_0$ of $P$  must lie in a singleton block $\{p_0\}$ in $\pi$.  So make $\Omega^{-1}(\pi)=\sigma$ have
$\sigma_1=p_0$, and construct $[\sigma_2, \ldots, \sigma_n]$ by applying $\Omega^{-1}$ recursively to the $(P-\{p_0\})$-transverse partition obtained from $\pi$
by removing the block $\{p_0\}$.

\vskip.1in
\noindent
{\sf Case 2.}
There are two minimal elements of $P$.

Label these two minimal elements $p_1, p_2$ of $P$ so that $p_i \in P_i$ for $i=1,2$.
Note that this implies that every $\sigma=[\sigma_1, \sigma_2, \ldots, \sigma_n]$
in $\linext{P}$ has either $\sigma_1=p_1$ or $\sigma_1=p_2$.
Note also that any $P$-transverse partition $\pi$ only has blocks of cardinality $1$ or $2$,
which yields two subcases for defining $\Omega$ and $\Omega^{-1}$:

\begin{itemize}
    \item The Subcase 2a for
     \begin{itemize}\item defining $\Omega$ occurs when $\sigma_1=p_1$,
                     \item defining $\Omega^{-1}$ occurs when $\{p_1\}$ appears as a singleton block within $\pi$.
    \end{itemize}
    \item   The Subcase 2b for
       \begin{itemize}\item defining $\Omega$ occurs when $\sigma_1=p_2$,
                      \item defining $\Omega^{-1}$ occurs when
    $p_1$ appears in a two-element block within $\pi$.
    \end{itemize}
\end{itemize}

\vskip.1in
\noindent
{\sf Subcase 2a.}\\
When defining $\Omega$, if $\sigma_1=p_1$, then make $\{p_1\}$ a singleton block of $\pi=\Omega(\sigma)$,
and produce the remaining blocks of $\pi$ by applying $\Omega$ recursively to $[\sigma_2,\ldots,\sigma_n]$.

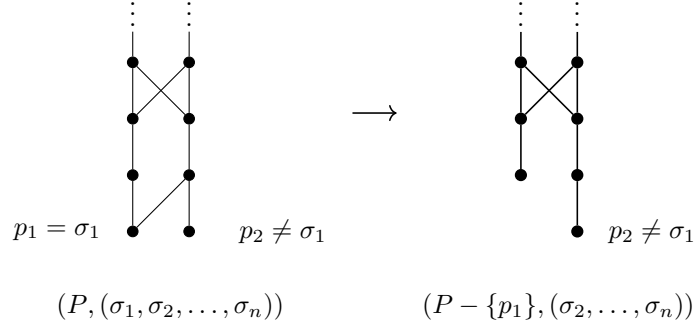
\begin{center}
    \begin{tikzpicture}[]

    \node[invisivertex] (P) at (-1,0){$p_1=\sigma_1$};
    \node[invisivertex] (P) at (2,0){$p_2 \neq \sigma_1$};
    \node[point] (1) at (0,0){};
    \node[point] (2) at (0,.75){};
    % \node[invisivertex] (2) at (0,1.25){$\vdots$};
    \node[point] (3) at (0,1.5){};
    \node[point] (4) at (0,2.25){};
    \node[point] (5) at (.75,0){};
    \node[point] (6) at (.75,.75){};
    \node[point] (7) at (.75,1.5){};
    \node[point] (8) at (.75,2.25){};

    \node[invisivertex] (dotsL) at (0,3){$\vdots$};
    \node[invisivertex] (dotsR) at (.75,3){$\vdots$};
    \node[invisivertex] (ell2) at (.5,-1){$(P,(\sigma_1,\sigma_2,\dots,\sigma_n))$};
    %\node[invisivertex] (ell3) at (.5,-1.5){$\pi=\{B_1,\dots,B_k\}$};

    \path[-] (1) edge [bend left =0] node[above] {} (2);
    \path[-] (1) edge [bend left =0] node[above] {} (6);
    \path[-] (2) edge [bend left =0] node[above] {} (3);
    \path[-] (3) edge [bend left =0] node[above] {} (4);
    \path[-] (5) edge [bend left =0] node[above] {} (6);
    \path[-] (6) edge [bend left =0] node[above] {} (7);
    \path[-] (7) edge [bend left =0] node[above] {} (8);
    \path[-] (3) edge [bend left =0] node[above] {} (8);
    \path[-] (4) edge [bend left =0] node[above] {} (7);
    \path[-] (4) edge [bend left =0] node[above] {} (dotsL);
    \path[-] (8) edge [bend left =0] node[above] {} (dotsR);
    \end{tikzpicture}
    \begin{tikzpicture}[-latex ,auto,on grid ,
    semithick ,
    state/.style ={ circle ,top color =white , bottom color = blue!10 ,
    draw, black , text=blue , minimum width =1 cm}]
    \node[invisivertex] (P1) at (0,-2){};
    \node[invisivertex] (P2) at (0,2){};
    \node[invisivertex] (1) at (0,.75){};
    \node[invisivertex] (2) at (.75,.75){};
    \path[->] (1) edge [bend left =0] node[above] {} (2);
    \end{tikzpicture}
    \begin{tikzpicture}[-latex ,auto,on grid ,
    semithick ,
    state/.style ={ circle ,top color =white , bottom color = blue!10 ,
    draw, black , text=blue , minimum width =1 cm}]

      % \node[invisivertex] (P) at (-.75,0){$p=\sigma_1$};
      \node[invisivertex] (P) at (1.75,0){$p_2 \neq \sigma_1$};
      % \node[magpoint] (1) at (0,0){};
      \node[point] (2) at (0,.75){};
      \node[point] (3) at (0,1.5){};
      % \node[invisivertex] (lablab) at (-.5,2){\textcolor{magenta}{$\sigma_{i-1}$}};
      \node[point] (4) at (0,2.25){};
      \node[point] (5) at (.75,0){};
      \node[point] (6) at (.75,.75){};
      \node[point] (7) at (.75,1.5){};
      \node[point] (8) at (.75,2.25){};

      \node[invisivertex] (dotsL) at (0,3){$\vdots$};
      \node[invisivertex] (dotsR) at (.75,3){$\vdots$};
      \node[invisivertex] (ell2) at (.5,-1){$(P-\{p_1\},(\sigma_2,\dots,\sigma_n))$};
      %\node[invisivertex] (ell3) at (.5,-1.5){$\pi=\{B_1,\dots,B_k,\{\sigma_1\}\}$};

      % \path[-] (1) edge [bend left =0] node[above] {} (2);
      % \path[-] (1) edge [bend left =0] node[above] {} (6);
      \path[-] (2) edge [bend left =0] node[above] {} (3);
      \path[-] (3) edge [bend left =0] node[above] {} (4);
      \path[-] (5) edge [bend left =0] node[above] {} (6);
      \path[-] (6) edge [bend left =0] node[above] {} (7);
      \path[-] (7) edge [bend left =0] node[above] {} (8);
      \path[-] (3) edge [bend left =0] node[above] {} (8);
      \path[-] (4) edge [bend left =0] node[above] {} (7);
      \path[-] (4) edge [bend left =0] node[above] {} (dotsL);
      \path[-] (8) edge [bend left =0] node[above] {} (dotsR);
    \end{tikzpicture}
  \end{center}

  To define $\Omega^{-1}$, if $\{p_1\}$ is a singleton block of $\pi$,
  make $\Omega^{-1}(\pi)=\sigma$ have
$\sigma_1=p_1$, and construct $[\sigma_2, \ldots, \sigma_n]$ by applying $\Omega^{-1}$ recursively to the $(P-\{p_1\})$-transverse partition obtained from $\pi$
by removing the block $\{p_1\}$.

\vskip.1in
\noindent
{\sf Subcase 2b.}\\
When defining $\Omega$, if $\sigma_1=p_2$, then $p_1$ appears elsewhere in $\sigma$, say $p_1=\sigma_{i+1}$ where $i \geq 1$.  Because $\sigma$ lies in $\linext{P}$ and $\sigma_{i+1}=p_1$ is the minimum element of $P_1$, this forces $\sigma_1,\sigma_2,\ldots,\sigma_{i}$
to all be elements of $P_2$.  In this case, add to $\pi$ the singleton blocks $\{\sigma_1\},\{\sigma_2\},\dots, \{\sigma_{i-1}\}$ along with the two-element block \edit{$\{\sigma_i,\sigma_{i+1}\}=\{\sigma_i,p_1\}$},
and compute the rest of $\Omega(\sigma)=\pi$ recursively by
replacing $(P,\sigma)$ with
  $
  (P-\{\sigma_1,\sigma_2\dots,\sigma_{i+1}\},
  [\sigma_{i+2},\sigma_{i+3},\dots,\sigma_n]).
  $
Here is the schematic picture:
  
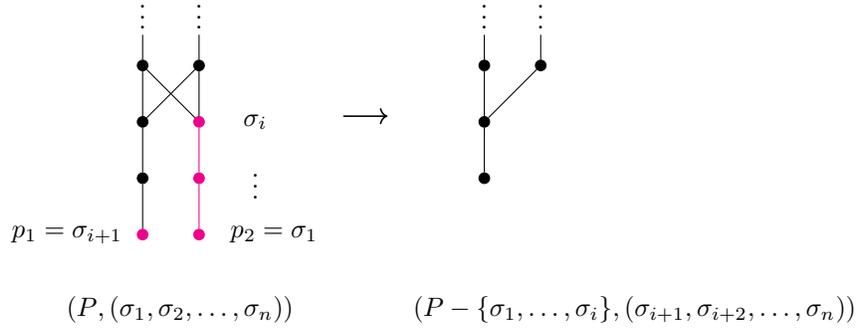
\begin{center}
    \begin{tikzpicture}[]

    \node[invisivertex] (P) at (-1,0){$p_1=\sigma_{i+1}$};
    \node[invisivertex] (P) at (1.75,0){$p_2=\sigma_1$};
    \node[invisivertex] (P) at (1.5,.75){$\vdots$};
    \node[invisivertex] (P) at (1.5,1.5){$\sigma_{i}$};
    \node[magpoint] (1) at (0,0){};
    \node[point] (2) at (0,.75){};
    % \node[invisivertex] (2) at (0,1.25){$\vdots$};
    \node[point] (3) at (0,1.5){};
    \node[point] (4) at (0,2.25){};
    \node[magpoint] (5) at (.75,0){};
    \node[magpoint] (6) at (.75,.75){};
    \node[magpoint] (7) at (.75,1.5){};
    \node[point] (8) at (.75,2.25){};

    \node[invisivertex] (dotsL) at (0,3){$\vdots$};
    \node[invisivertex] (dotsR) at (.75,3){$\vdots$};
    \node[invisivertex] (ell2) at (.5,-1){$(P,(\sigma_1,\sigma_2,\dots,\sigma_n))$};
    %\node[invisivertex] (ell3) at (.5,-1.5){$\pi=\{B_1,\dots,B_k\}$};

    \path[-] (1) edge [bend left =0] node[above] {} (2);
    % \path[-] (1) edge [bend left =0,color=magenta] node[above] {} (6);
    \path[-] (2) edge [bend left =0] node[above] {} (3);
    \path[-] (3) edge [bend left =0] node[above] {} (4);
    \path[-] (5) edge [bend left =0,color=magenta] node[above] {} (6);
    \path[-] (6) edge [bend left =0,color=magenta] node[above] {} (7);
    \path[-] (7) edge [bend left =0] node[above] {} (8);
    \path[-] (3) edge [bend left =0] node[above] {} (8);
    \path[-] (4) edge [bend left =0] node[above] {} (7);
    \path[-] (4) edge [bend left =0] node[above] {} (dotsL);
    \path[-] (8) edge [bend left =0] node[above] {} (dotsR);
    \end{tikzpicture}
    \begin{tikzpicture}[-latex ,auto,on grid ,
    semithick ,
    state/.style ={ circle ,top color =white , bottom color = blue!10 ,
    draw, black , text=blue , minimum width =1 cm}]
    \node[invisivertex] (P1) at (0,-2){};
    \node[invisivertex] (P2) at (0,2){};
    \node[invisivertex] (1) at (0,.75){};
    \node[invisivertex] (2) at (.75,.75){};
    \path[->] (1) edge [bend left =0] node[above] {} (2);
    \end{tikzpicture}
    \begin{tikzpicture}[]

    % \node[invisivertex] (P) at (-1,0){$p=\sigma_i$};
    % \node[invisivertex] (P) at (1.75,0){$q=\sigma_1$};
    % \node[invisivertex] (P) at (1.5,1){$\vdots$};
    % \node[invisivertex] (P) at (1.75,2){$\sigma_{i-1}$};
    % \node[magpoint] (1) at (0,0){};
    \node[point] (2) at (0,.75){};
    % \node[invisivertex] (2) at (0,1.25){$\vdots$};
    \node[point] (3) at (0,1.5){};
    \node[point] (4) at (0,2.25){};
    % \node[magpoint] (5) at (1,0){};
    % \node[magpoint] (6) at (1,1){};
    % \node[magpoint] (7) at (1,2){};
    \node[point] (8) at (.75,2.25){};

    \node[invisivertex] (dotsL) at (0,3){$\vdots$};
    \node[invisivertex] (dotsR) at (.75,3){$\vdots$};
    \node[invisivertex] (ell2) at (2,-1){$(P-\{\sigma_1,\dots,\sigma_i\},
    (\sigma_{i+1},\sigma_{i+2},\dots,\sigma_n))$};
    %\node[invisivertex] (ell3) at (.5,-1.5){$\pi=\{B_1,\dots,B_k,\{\sigma_1\},\{\sigma_2\},\dots, \{\sigma_{i-2}\},\{\sigma_{i-1},\sigma_i\}\}$};

    % \path[-] (1) edge [bend left =0] node[above] {} (2);
    % \path[-] (1) edge [bend left =0,color=magenta] node[above] {} (6);
    \path[-] (2) edge [bend left =0] node[above] {} (3);
    \path[-] (3) edge [bend left =0] node[above] {} (4);
    % \path[-] (5) edge [bend left =0,color=magenta] node[above] {} (6);
    % \path[-] (6) edge [bend left =0,color=magenta] node[above] {} (7);
    % \path[-] (7) edge [bend left =0] node[above] {} (8);
    \path[-] (3) edge [bend left =0] node[above] {} (8);
    % \path[-] (4) edge [bend left =0] node[above] {} (7);
    \path[-] (4) edge [bend left =0] node[above] {} (dotsL);
    \path[-] (8) edge [bend left =0] node[above] {} (dotsR);
    \end{tikzpicture}
  \end{center}

When defining $\Omega^{-1}(\pi)$, if $p_1$ appears in some two-element block of $\pi$, then it appears in some block $\{p_1,p_2'\}$ for some $p_2'$ in $P_2$.   We claim that
$\pi$ being $P$-transverse then forces any elements $p <_P p_2'$ in $P_2$ to lie in singleton blocks $\{p\}$ of $\pi$.
To see this claim, assume not, so that some such $p$ lies in a two-element block of $\pi$, necessarily of the form $\{p_1',p\}$ for some $p_1'$ in $P_1$ with $p_1 <_P p_1'$.
This leads to a contradiction of Proposition~\ref{cones:prop:transverse}(ii), since  $(p_1',p_1)$ would then be a relation in $\overline{P \cup \pi}$ via this transitive chain of relations:
$
p_1' \equiv_\pi p <_P p_2' \equiv_\pi p_1.
$

In this subcase, list the totally ordered (and possibly empty) collection of all elements $p$ in $P_2$ with $p <_P p_2'$ at the beginning of $\sigma$ as $\sigma_1,\sigma_2,\ldots,\sigma_{i-1}$, followed by
$\sigma_i \sigma_{i+1}=p_2' p_1$.
Then compute the rest of $\Omega^{-1}(\pi)=\sigma$ recursively,
by applying $\Omega^{-1}$ to the $(P-\{\sigma_1,\sigma_2,\ldots,\sigma_{i+1}\})$-transverse
partition obtained from $\pi$ by removing the singleton blocks $\{\sigma_1\},\{\sigma_2\}, \ldots, \{\sigma_{i-1}\}$ and
the two-element block $\{\sigma_i ,\sigma_{i+1}\}=\{p_2', p_1\}$.

  \bigskip

It is not hard to check that the two maps $\Omega, \Omega^{-1}$
defined recursively in this way are actually mutually inverse bijections.  By construction, $\Omega$ has the property that the
two-element blocks of $\pi=\Omega(\sigma)$ are exactly those containing $P$-incomparable pairs $\{\sigma_i,\sigma_{i+1}\}$ for which $\sigma_i \in P_1$ and $\sigma_{i+1} \in P_2$, as claimed.
\end{proof}

Recall from the Introduction that the number of usual descents of a permuation $\sigma$ is defined as
$\des(\sigma)=\#\Des(\sigma)$ where
$$
\Des(\sigma):= \{i\in[n-1] \mid \sigma_i>\sigma_{i+1}\}.
$$
This was used to define the $P$-Eulerian polynomial in equation \eqref{P-Eulerian-polynomial} as
$\sum_{\sigma \in \linext{P}} t^{\des(\sigma)}$, assuming that $P$ is naturally labeled, that is, $\linext P$ contains the identity \edit{permutation} $\sigma=[1,2,\ldots,n]$.

\begin{coro}\label{width-2:coro:descents}
  When $P$ is a width two poset having a chain decomposition $P_1\cup P_2$ with $P_1$
  an order ideal of $P$, then the Poincar\'e polynomial for $P$ coincides with the $P$-Eulerian polynomial:
  $$
  \Poin(P,t)=\sum_{\sigma \in \linext{P}} t^{\des(\sigma)}.
  $$
\end{coro}

\begin{proof}
Let $\#P_i=n_i$ for $i=1,2$, so that $n=\#P=n_1+n_2$.
One can then choose a natural labeling for $P$ by $[n]$ such that the elements of the order ideal $P_1$ are labeled by the initial segment $[n_1]=\{1,2,\ldots,n_1\}$,
and $P_2$ is labeled by $\{n_1+1,n_1+2,\ldots,n\}$.  We claim that with this
natural labeling, one has $\Des_{(P_1,P_2)}(w) = \Des(w)$.  This is because
the labeling renders one of the
conditions in the definition \eqref{width-two-Des-defn} of
$\Des_{(P_1,P_2)}(\sigma)$ superfluous:  assuming that
$\sigma_i \in P_2$ and $\sigma_{i+1} \in P_1$, then $\sigma_i, \sigma_{i+1}$ must
already be incomparable in $P$, because otherwise $\sigma_i <_P \sigma_{i+1}$
(since $\sigma$ was a linear extension of $P$) and then $P_1$ being an order ideal
would force $\sigma_i \in P_1$, a contradiction.
Now since $P_1, P_2$ are totally ordered in $P$, and $\sigma$ lies in $\linext{P}$, one has $\sigma_i \in P_2$
and $\sigma_{i+1}$ in $P_1$ if and only if
$\sigma_i>_\Z \sigma_{i+1}$. That is $\Des_{(P_1,P_2)}(w) = \Des(w)$.
\end{proof}

\begin{ex}
\label{two-row-skew-example}
An interesting family of posets to which
Corollary~\ref{width-2:coro:descents} applies
are the posets $P(\lambda/\mu)$ associated with
two-row skew Ferrers diagrams $\lambda/\mu$.
A \emph{Ferrers diagram} associated to a partition (of a number) $\lambda=(\lambda_1,\ldots,\lambda_\ell)$
has $\lambda_i$ square cells drawn left-justified in row $i$. A \emph{skew Ferrers diagram} $\lambda/\mu$ for
two partitions $\lambda, \mu$ having $\lambda_i \geq \mu_i$ is the diagram for $\lambda$ with the cells occupied by the diagram for $\mu$ removed.
There is a poset structure $P(\lambda/\mu)$ on the cells of $\lambda/\mu$
in which a cell $(i,j)$ in row $i$ and column $j$
has $(i,j) \leq_{P(\lambda/\mu)} (i',j')$
if both $i \leq i'$ and $j \leq j'$.

When $\lambda/\mu$ has only two parts, we will call it a \emph{two-row skew Ferrers diagram}.
Three examples of
such $\lambda/\mu$ and their associated $P(\lambda/\mu)$ are shown below.
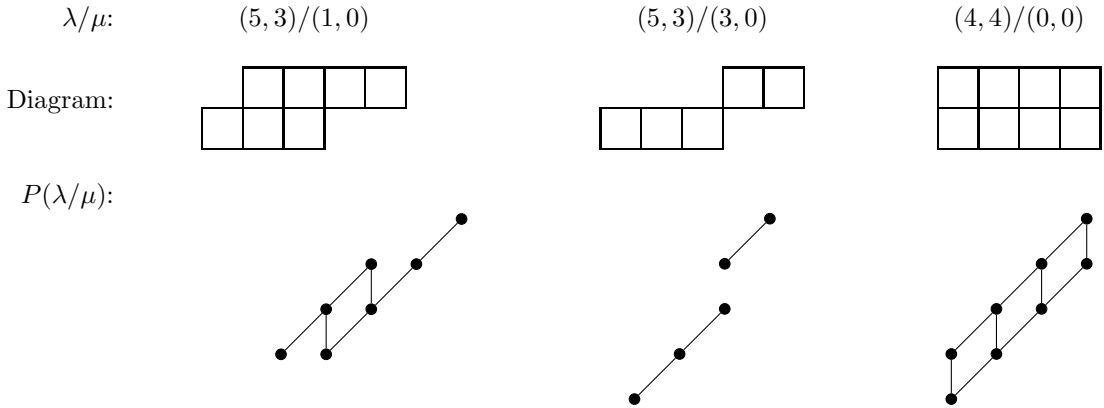
\begin{center}
  \begin{tabular}{r c  c c}
    $\lambda/\mu$: & $(5,3)/(1,0)$ & $(5,3)/(3,0)$ & $(4,4)/(0,0)$ \\
    &&&\\
    Diagram: & $\ydiagram{1+4,0+3}$ & $\ydiagram{3+2,0+3}$ & $\ydiagram{0+4,0+4}$\\
    &&&\\
    $P(\lambda/\mu)$: &&&\\
    &
    \begin{tikzpicture}[scale=.6]
    \node[invisivertex] (1) at (-2,0){};
      \node[point] (2) at (2,1){};
      \node[point] (3) at (3,2){};
      \node[point] (7) at (4,3){};
      \node[point] (4) at (1,1){};
      \node[point] (5) at (2,2){};
      \node[point] (6) at (3,3){};
    %   \node[point] (8) at (-1,4){};
      \node[point] (9) at (5,4){};

    %   \path[-] (1) edge [bend left =0] node[above] {} (2);
      \path[-] (2) edge [bend left =0] node[above] {} (3);
      \path[-] (4) edge [bend left =0] node[above] {} (5);
      \path[-] (5) edge [bend left =0] node[above] {} (6);
    %   \path[-] (1) edge [bend left =0] node[above] {} (4);
      \path[-] (2) edge [bend left =0] node[above] {} (5);
      \path[-] (3) edge [bend left =0] node[above] {} (6);
    %   \path[-] (6) edge [bend left =0] node[above] {} (8);
      \path[-] (3) edge [bend left =0] node[above] {} (7);
    %   \path[-] (7) edge [bend left =0] node[above] {} (8);
      \path[-] (7) edge [bend left =0] node[above] {} (9);
    \end{tikzpicture}
    &
    \begin{tikzpicture}[scale=.6]
    \node[invisivertex] (1) at (-2,0){};
    \node[invisivertex] (1) at (7,0){};
      \node[point] (1) at (1,0){};
      \node[point] (2) at (2,1){};
      \node[point] (3) at (3,2){};
    %   \node[point] (7) at (-1,3){};
    %   \node[point] (4) at (2,1){};
    %   \node[point] (5) at (1,2){};
      \node[point] (6) at (3,3){};
       \node[point] (8) at (4,4){};
    %   \node[point] (9) at (-2,4){};

      \path[-] (1) edge [bend left =0] node[above] {} (2);
      \path[-] (2) edge [bend left =0] node[above] {} (3);
    %   \path[-] (4) edge [bend left =0] node[above] {} (5);
    %   \path[-] (5) edge [bend left =0] node[above] {} (6);
    %   \path[-] (1) edge [bend left =0] node[above] {} (4);
    %   \path[-] (2) edge [bend left =0] node[above] {} (5);
    %   \path[-] (3) edge [bend left =0] node[above] {} (6);
      \path[-] (6) edge [bend left =0] node[above] {} (8);
    %   \path[-] (3) edge [bend left =0] node[above] {} (7);
    %   \path[-] (7) edge [bend left =0] node[above] {} (8);
    %   \path[-] (7) edge [bend left =0] node[above] {} (9);
    \end{tikzpicture}
    &
    \begin{tikzpicture}[scale=.6]
      \node[point] (1) at (1,0){};
      \node[point] (2) at (2,1){};
      \node[point] (3) at (3,2){};
      \node[point] (7) at (4,3){};
      \node[point] (4) at (1,1){};
      \node[point] (5) at (2,2){};
      \node[point] (6) at (3,3){};
      \node[point] (8) at (4,4){};

      \path[-] (1) edge [bend left =0] node[above] {} (2);
      \path[-] (2) edge [bend left =0] node[above] {} (3);
      \path[-] (4) edge [bend left =0] node[above] {} (5);
      \path[-] (5) edge [bend left =0] node[above] {} (6);
      \path[-] (1) edge [bend left =0] node[above] {} (4);
      \path[-] (2) edge [bend left =0] node[above] {} (5);
      \path[-] (3) edge [bend left =0] node[above] {} (6);
      \path[-] (6) edge [bend left =0] node[above] {} (8);
      \path[-] (3) edge [bend left =0] node[above] {} (7);
      \path[-] (7) edge [bend left =0] node[above] {} (8);
      \end{tikzpicture}
  \end{tabular}
\end{center}

The decomposition $P(\lambda/\mu)=P_1 \cup P_2$
where $P_i$ correspond to the cells in row $i$ of $\lambda/\mu$ shows that $P(\lambda/\mu)$ has width two, and furthermore $P_1$ forms an order ideal.
Therefore Corollary~\ref{width-2:coro:descents} implies that for any two-row skew Ferrers diagram $\lambda/\mu$ one has
\begin{equation}
\label{two-row-skew-diagram-poincare}
\Poin(P(\lambda/\mu),t)=
\sum_{\sigma \in \linext{P(\lambda/\mu)}} t^{\des(\sigma)}.
\end{equation}

On the other hand, there is a well-known bijection between linear extensions $\sigma$ of $P(\lambda/\mu)$ and the \emph{standard Young tableaux} $Q$ of
shape $\lambda/\mu$, which are (bijective) labelings of the cells of the diagram by $[n]$ where $n=\sum_i \lambda_i - \sum_i \mu_i$, with the numbers increasing left-to-right in rows and top-to-bottom in columns; see \cite[\S 7.10]{ec2}.  There is also a notion of \emph{descent set} $\Des(Q)$ for such tableaux, having $i \in \Des(Q)$ whenever $i+1$ labels a cell in a lower row of $Q$
than $i$. However, in general
when $\sigma$ corresponds to $Q$, one does \emph{not}
have $\des(\sigma)=\des(Q)$, so that $\Poin(P(\lambda/\mu),t)$ differs from the
generating function $\sum_Q t^{\des(Q)}$
of standard tableaux $Q$ shape $\lambda/\mu$ by $\des(Q)$.
For example, there are two standard tableaux of shape $\lambda/\mu=(2,1)/(0,0)$
$$
Q_1=\begin{ytableau} 1&2\\3  \end{ytableau} \qquad
Q_2=\begin{ytableau} 1&3\\2  \end{ytableau}
$$
both having $\des(Q_i)=1$, however $\Poin(P(\lambda/\mu),t)=1+t$.

In two special cases, however, they (essentially) coincide.
\begin{itemize}
    \item When $P_{\lambda/\mu} =\chain{a} \sqcup \chain{b}$ is a disjoint union of two chains, as in Example~\ref{two-disjoint-chain-example}, one can check that, if one (naturally) labels $\chain{a} \sqcup \chain{b}$ so
    that the elements of the order ideal $\chain{b}$ are labeled $1,2,\ldots,b$
    while $\chain{a}$ is labeled $b+1,b+2,\ldots,b+a$,
    then one \emph{does} have $\des(\sigma)=\des(Q)$, and hence
    $$
     \sum_Q t^{\des(Q)}  = \Poin( P_{\chain{a} \sqcup \chain{b}},t)= \sum_k \binom{a}{k} \binom{b}{k} t^k.
    $$
    \item When $\lambda/\mu$ is a $2 \times n$ rectangle, so that $P_{\lambda/\mu}=\chain{2} \times \chain{n}$ is a Cartesian product poset, then $\sigma$ in $\linext{P}$ and standard Young tableaux $Q$ of shape $2 \times n$ can both be
    identified with \emph{Dyck paths of semilength $n$}, that is, lattice paths from $(0,0)$ to $(2n,0)$ in $\Z^2$ taking steps northeast or \edit{southeast} and staying weakly above the $x$-axis.
    One can check that
    \begin{itemize}
        \item $\Des(\sigma)$ corresponds to
    \emph{valleys} (i.e. southwest steps followed by a northeast step), while
    \item $\Des(Q)$ correspond to \emph{peaks} (i.e. northeast steps followed by a southwest step).
    \end{itemize}
    In general, such a Dyck path has one more peak than valley \cite[Exercises 6.19(i,\thinspace ww,\thinspace aaa)]{ec2}.  Hence one has
    $$
     \Poin( \chain{2} \times \chain{n}, t) = \frac{1}{t} \sum_Q t^{\des(Q)} = \sum_{k=0}^{n-1} \frac{1}{n} \binom{n}{k} \binom{n}{k+1} t^k,
    $$
    which is the generating function for the {\it Narayana numbers} $N(n,k):=\frac{1}{n} \binom{n}{k-1} \binom{n}{k}$ (see \cite[p.2]{Branden} and \cite[Exer. 6.36(a)]{ec2}).  Upon setting $t=1$, the Naryana numbers sum to the {\it Catalan number}
    %\cite[Corollary 6.2.3]{ec2} (check this ref ???)
    $$
    \#\linext{\chain{2} \times \chain{n}}=
    \frac{1}{n+1}\binom{2n}{n}.
    $$
\end{itemize}

Note that for any (non-skew) partition $\lambda$, the celebrated \emph{hook-length formula} of
Frame, Robinson and Thrall \cite[Corollary 7.21.6]{ec2} gives a simple product formula for
$
\#\linext{P(\lambda)}
=\left[ \Poin(P(\lambda),t)\right]_{t=1}.
$

\begin{prob}\label{open-problem-perms}
  Combinatorially interpret
  $\Poin(P(\lambda),t)$ for other partitions $\lambda$, and in particular, for $m \times n$ rectangular partitions,  where $P(\lambda) = \chain{m} \times \chain{n}$ is a
  Cartesian product of chains.

\end{prob}

\noindent Below we give $\Poin(\chain{3} \times \chain{n},t)$ for $2\leq n\leq 8$.

\bigskip

%\begin{equation}
%\label{3-by-n-data}
%\begin{aligned}
%\Poin(\chain{3} \times \chain{2},t)
%&=1+3t+t^2,\\
%\Poin(\chain{3} \times \chain{3},t)
%&=1+9t+19t^2+11t^3+2t^4,\\
%\Poin(\chain{3} \times \chain{4},t)
%&=1+18t+92t^2+174t^3+133t^4+40t^5+4t^6,\\
%\Poin(\chain{3} \times \chain{5},t)
%&=1 + 30t + 280t^{2} + 1091t^{3} + 1987t^{4} + 1746t^{5} + 731t^{6} + 132t^{7} + 8t^{8},\\
%\Poin(\chain{3} \times \chain{6},t)
%&=1 + 45t + 665t^{2} + 4383t^{3} + 14603t^{4} + 25957t^{5} + 25064t^{6}+ 12965t^{7}+ 3413t^{8}+ 404t^{9}\\
%&\quad + 16t^{10},\\
%\Poin(\chain{3} \times \chain{7},t)
%&=1 + 63t + 1351t^{2} + 13475t^{3} + 71305t^{4} + 213539t^{5} + 373651t^{6} + 385578t^{7} + 232310t^{8}\\
%&\quad + 79023t^{9} + 14174t^{10} + 1168t^{11} + 32t^{12},\\
%\Poin(\chain{3} \times \chain{8},t)
%&=1 + 84t + 2464t^{2} + 34608t^{3} + 266470t^{4} + 1206826t^{5} + 3343958t^{6} + 5782699t^{7}\\
%&\quad + 6275503t^{8} + 4240489t^{9} + 1743730t^{10} + 417622t^{11} + 53884t^{12} + 3232t^{13} + 64t^{14}\\
%\end{aligned}
%\end{equation}

\begin{tabular}{|c|l|}\hline
  $n$ & $\Poin(3 \times n,t)$ \\\hline\hline
%& \\
  $2$ & $1+3t+t^2$\\
\hline
%& \\
  $3$ & $1+9t+19t^2+11t^3+2t^4$\\
\hline
%& \\
  $4$ & $1+18t+92t^2+174t^3+133t^4+40t^5+4t^6$ \\
\hline
%& \\
  $5$ & $1 + 30t + 280t^{2} + 1091t^{3} + 1987t^{4} + 1746t^{5} + 731t^{6} + 132t^{7} + 8t^{8}$\\
\hline
%& \\
  $6$ & $1 + 45t + 665t^{2} + 4383t^{3} + 14603t^{4} + 25957t^{5} + 25064t^{6}+ 12965t^{7}+ 3413t^{8}+ 404t^{9}$ \\
\hline
%& \\
  $7$ & $1 + 63t + 1351t^{2} + 13475t^{3} + 71305t^{4} + 213539t^{5} + 373651t^{6} + 385578t^{7} + 232310t^{8}$\\
&\quad $+ 79023t^{9} + 14174t^{10} + 1168t^{11} + 32t^{12}$\\
%& \\\hline
\hline
  $8$ & $1 + 84t + 2464t^{2} + 34608t^{3} + 266470t^{4} + 1206826t^{5} + 3343958t^{6} + 5782699t^{7}$\\
      &\quad $+ 6275503t^{8} + 4240489t^{9} + 1743730t^{10} + 417622t^{11} + 53884t^{12} + 3232t^{13} + 64t^{14}$\\
  \hline
\end{tabular}

\bigskip

\end{ex}

\begin{rem}
  Since equation \eqref{stirling-number-gf} shows
  that the Poincar\'e polynomial $\Poin(P,t)$ for the
  antichain poset $P=P_{(1^n)}=\antichain_n$ has {\it only real roots},
  one might wonder whether this holds for some more general
  class of posets.
  It does not hold for all posets, and not even for all disjoint unions of chains $\poset{a}$,
  since

%  nor for all skew posets $P_{\lambda/\m}$, since when $\lambda = (6,4,2)$ and
%  $\mu = (4,2)$ we have $P_{\lambda/\mu}=\mathbf{2}\sqcup\mathbf{2}\sqcup\mathbf{2}$ and
\[
  % \Poin(P_{\lambda/\mu},t) =
  \Poin(P_{(2,2,2)},t)=1+12t+43t^2+30t^3+4t^4
\]
has a pair of non-real complex roots.
It can fail even for rectangular Ferrers posets, e.g.,
$\lambda=(3,3,3)$ has
$
\Poin(P(\lambda),t)=\Poin(\chain{3} \times \chain{3},t)=1+9t+19t^2+11t^3+2t^4
$
in the above table, with two non-real complex roots.

On the other hand, computations show that $\Poin(P,t)$ \emph{is} real-rooted
for all posets $P$ of width two having at most $9$
elements.
This leads to the following question.

\begin{question}
Is $\Poin(P,t)$ real-rooted when the poset $P$ has width two?
\end{question}

\end{rem}

\section*{Acknowledgements}

The authors gratefully acknowledge Dennis Stanton for conversations about MacMahon's Master Theorem,
as well as  \edit{Anders Bj\"orner}, Jesus DeLoera, Theo Douvropoulos, Michael Falk, \edit{Ira Gessel}, Benjamin Steinberg, Volkmar Welker, Chi-Ho Yuen for enlightening discussions and references.
They thank Philip Zhang for asking them questions about real-rootedness
at the 2019 Mid-Atlantic Algebra, Geometry, and Combinatorics Workshop.
\edit{Finally, the authors thank an anonymous referee for helpful comments.}
%% if you use biblatex then this generates the bibliography
%% if you use some other method then remove this and do it your own way
%\bibliography{paper}{}
%\bibliographystyle{plain}

\end{document}